\documentclass{article}[12pt]
\usepackage{graphicx,pstricks,comment}
\usepackage[all]{xy}
\usepackage{subfigure}
\usepackage{rotating}
\usepackage{amsfonts}
\usepackage{amssymb}
\usepackage{amsmath}

\usepackage[utf8]{inputenc}

\usepackage{tikz}
\usepackage{pbox}

\usepackage{hyperref}

\newcommand{\redbluedot}{\begin{scope}[scale=.1]
\draw[color=black,fill=black] (90:1) arc (90:270:1);
\draw[color=black,fill=white] (90:1) arc (90:-90:1);
\end{scope}}
\newcommand{\reddot}{\begin{scope}[scale=.1]
\draw[color=black,fill=white] (0:1) arc (0:360:1);
\end{scope}}
\newcommand{\bluedot}{\begin{scope}[scale=.1]
\draw[color=black,fill=black] (0:1) arc (0:360:1);
\end{scope}}

\newcommand{\Zerodot}[1]{\begin{scope}[shift = {(0,0)}]#1\end{scope}}
\newcommand{\idot}[1]{\begin{scope}[shift = {(-.8,1)}]#1\end{scope}}
\newcommand{\jdot}[1]{\begin{scope}[shift = {(0,1)}]#1\end{scope}}
\newcommand{\kdot}[1]{\begin{scope}[shift = {(.8,1)}]#1\end{scope}}
\newcommand{\ijdot}[1]{\begin{scope}[shift = {(-.8,2)}]#1\end{scope}}
\newcommand{\ikdot}[1]{\begin{scope}[shift = {(0,2)}]#1\end{scope}}
\newcommand{\jkdot}[1]{\begin{scope}[shift = {(.8,2)}]#1\end{scope}}
\newcommand{\ijkdot}[1]{\begin{scope}[shift = {(0,3)}]#1\end{scope}}
\newcommand{\dotify}[8]{\Zerodot{#1},\idot{#2},\jdot{#3},\kdot{#4},\ijdot{#5},\ikdot{#6},\jkdot{#7},\ijkdot{#8}}

\newcommand{\cube}{
\draw (0,0) -- (-.8,1) -- (0,2) -- (.8,1) -- cycle;
\draw (-.8,1) -- (-.8,2) -- (0,3) -- (.8,2) -- (.8,1);
\draw (0,2) -- (0,3);}

\usepackage{indentfirst}
\usepackage{epsfig}
\usepackage{psfrag}
\usepackage{enumerate}
\usepackage{array}
\usepackage{theorem}

\newcommand{\C}{\mathbb{C}}
\newcommand{\R}{\mathbb{R}}
\newcommand{\N}{\mathbb{N}}
\newcommand{\Z}{\mathbb{Z}}

\newcommand{\Q}{\mathbb{Q}}
\newcommand{\K}{\mathbb{K}}
\newcommand{\quat}{\mathbb{H}}
\newcommand{\oct}{\mathbb{O}}

\newcommand{\la}{\langle}
\newcommand{\ra}{\rangle}

\newtheorem{thm}{Theorem}
\newtheorem{cor}{Corollary}
\newtheorem{lem}{Lemma}
\newtheorem{prop}{Proposition}
\newtheorem{dfn}{Definition}

{\theorembodyfont{\rmfamily} }
{\theorembodyfont{\rmfamily} }

\newcommand{\Pf}{{\em Proof}. }
\newcommand{\EPf}{\hfill$\Box$\vspace{.5cm}}

\title{Presentations for cusped arithmetic hyperbolic lattices}

\author{Alice Mark, Julien Paupert\footnote{Second author partially supported by National Science Foundation Grant DMS-1708463.}}

\begin{document}
\maketitle

\begin{abstract}
We present a general method to compute a presentation for any cusped arithmetic hyperbolic lattice $\Gamma$, applying a classical result of Macbeath to a suitable $\Gamma$-invariant horoball cover of the corresponding symmetric space.  As applications we compute presentations for the Picard modular groups ${\rm PU}(2,1,\mathcal{O}_d)$ for $d=1,3,7$ and the quaternion hyperbolic lattice ${\rm PU}(2,1,\mathcal{H})$ with entries in the Hurwitz integer ring $\mathcal{H}$. The implementation of the method for these groups is computer-assisted.
\end{abstract}


\section{Introduction}

Discrete subgroups and lattices in semisimple Lie groups form a rich and well-studied class of finitely generated groups acting on non-positively curved metric spaces. The case of real rank one, where the associated symmetric space is negatively curved, is of special interest. There are essentially two main families of constructions of such lattices, arithmetic on one hand and geometric on the other. Arithmetic lattices are roughly speaking obtained by taking matrices with entries lying in the integer ring of some number field; the general definition is more complicated and we will not give it here, as the arithmetic lattices that we consider in this paper are of this simplest type. By Margulis' celebrated superrigidity and arithmeticity theorems, all (irreducible) lattices in $G$ are of this arithmetic type when $G$ is a semisimple Lie group of real rank at least 2. 

The other family involves geometric constructions such as polyhedra, reflections or other types of involutions or other finite-order isometries. A prototype of this type of construction is given by Coxeter groups in the constant curvature geometries $E^n$, $S^n$ and $H^n$, which are generated by reflections across hyperplanes. These groups are classical and were classified by Coxeter in the spaces $E^n$ and $S^n$, whereas their hyperbolic counterparts (studied by Vinberg and others) are still not completely understood. However by construction these groups come equipped with data including a presentation (as an abtract Coxeter group) and a fundamental domain for their action on the symmetric space.  

Arithmetic lattices are given by a global description and their global structure is in some sense well understood by work of Siegel, Borel, Tits, Prasad and others. However concrete information such as a presentation and a fundamental domain are not readily accessible from the arithmetic construction. One can obtain geometric information such as volume by Prasad's celebrated volume formula (\cite{Pr}) but computing the constants appearing in this formula usually involves some non-trivial work (see for example \cite{Be} and \cite{Sto}). 

Very few presentations of arithmetic lattices, and of lattices in general, are known. 
Presentations can provide useful geometric and algebraic information about groups, such as explicit index of torsion-free subgroups (effective Selberg lemma, as used for example in \cite{Sto}), cohomology of the group $\Gamma$ or quotient space $X/\Gamma$, see for instance \cite{Y} (for the Picard modular groups with $d=1,3$) and of course representations of $\Gamma$, for instance if one is interested in deformations of $\Gamma$ into a larger Lie group.

Presentations for ${\rm SL}(n,\Z)$ with $n \geqslant 3$ were given by Steinberg (\cite{Ste}, following  Magnus); the case of ${\rm SL}(2,\Z)$ is classical and possibly dates to Gauss; see also Siegel \cite{Si}. In rank one, Swan gave in \cite{Sw} presentations for the Bianchi groups ${\rm PGL}(2,\mathcal{O}_d)$ (where $\mathcal{O}_d$ denotes the ring of integers of $\Q [i\sqrt{d}]$ for $d$ a positive square-free integer), following Bianchi's original construction in \cite{Bi}. These act as isometries of (real) hyperbolic 3-space, as they are lattices in ${\rm PGL}(2,\C) \simeq {\rm Isom}^+({\rm H}_\R^3)$.

Presentations for the related Picard modular groups ${\rm PU}(2,1,\mathcal{O}_d)$ were found only recently in the simplest cases of $d=3$ (\cite{FP}) and $d=1$ (\cite{FFP}). One of the reasons for this is that the associated symmetric space, complex hyperbolic 2-space ${\rm H}^2_\C$, is more complicated and in particular has non-constant (pinched) negative curvature. A particular feature of such spaces, the absence of totally geodesic real hypersurfaces, makes constructions of fundamental domains difficult as there are no obvious walls to use to bound such domains. The presentations obtained for $d=1,3$ were in fact obtained by constructing fundamental domains and using the Poincar\'e polyhedron theorem. This approach seems to become too complicated when considering more complicated groups, such as Picard modular groups with higher values of $d$, and no further such constructions have appeared. Using a similar strategy, Zhao gave in \cite{Zh} generating sets for the Picard modular groups with $d=1,2,3,7,11$ but he does not go as far as obtaining a presentation, finding a set whose translates covers the space but without control over intersections and cycles. (We will in fact use a covering argument closely related to the one he uses to cover a fundamental prism on the ideal boundary by isometric spheres, see Lemma~\ref{d=7cover}).

 In this paper we present a method for obtaining presentations for cusped hyperbolic lattices, i.e. non-cocompact lattices in semisimple Lie groups of real rank one, based on a classical result of Macbeath  (Theorem~\ref{macbeath} below) which gives a presentation for a group $\Gamma$ acting by homeomorphisms on a topological space $X$, given an open subset $V$ whose $\Gamma$-translates cover $X$. We apply this by finding a suitable horoball $V$ based at a cusp point of $\Gamma$ whose $\Gamma$-translates cover $X$, then analyzing the triple intersections and associated cyles to obtain a presentation for $\Gamma$. The main tools for this analysis come from the additional arithmetic structure that we get by assuming that $\Gamma$ is in fact an \emph{integral lattice} in the sense that it is contained in ${\rm GL}(n+1,\mathcal{O}_E)$ for some number field $E$ (or finitely generated division algebra over $\Q$). The crucial such tool that we use is the notion of \emph{level} between two $E$-rational boundary points in $\partial_\infty X$ (see Definition~\ref{level}) which gives a notion of distance between such points using only algebraic data. More importantly for us, levels measure the relative sizes of horospheres based at the correpsonding boundary points, which allows us to control whether or not such horospheres intersect at a given height (see Lemma~\ref{depth-height}). 
   
As applications of this method we compute presentations for the Picard modular groups ${\rm PU}(2,1,\mathcal{O}_d)$ with $d=3,1,7$, given in the appendix and Propositions~\ref{presentationd=1}, \ref{presentationd=7} respectively. The cases $d=2,11$ can be treated with the same method but are computationally more intensive; David Polletta has treated these in \cite{Po}. 
We also compute a presentation for the quaternion hyperbolic lattice which we call the \emph{Hurwitz modular group} ${\rm PU}(2,1,\mathcal{H})$, where $\mathcal{H}$ is the ring of Hurwitz integers $\mathcal{H}=\Z [i,j,k,\frac{1+i+j+k}{2}] \subset \quat$. This is a lattice in ${\rm PU}(2,1,\quat)$ (also denoted ${\rm PSp}(2,1)$), acting on the 8-dimensional symmetric space ${\rm H}^2_\quat$. As far as we know this is the first presentation ever found for a higher-dimensional quaternion hyperbolic lattice. (In dimension 1, ${\rm H}^1_\quat \simeq {\rm H}^4_\R$ and such groups have been studied e.g. in \cite{DVV}, see also \cite{A}, \cite{W} and \cite{Ph}). The implementation of our method for this group turned out to be computationally much more intensive than anticipated, and in fact the resulting presentation is too large for many purposes (see Section 5). It does however allow us to compute the abelianization of the lattice, and to find a nice generating set, see Theorems~\ref{hurwitzabel} and \ref{hurwitzgens}.

Using a somewhat similar principle Cartwright and Steger found presentations for some cocompact arithmetic lattices in ${\rm PU}(2,1)$ in their classification of the so-called fake projective planes, see \cite{CaS1} and \cite{CaS2}.

The paper is organized as follows. In Section 2 we discuss generalities about horoball coverings of hyperbolic spaces, levels for cusp points of integral lattices, and outline how we apply Macbeath's theorem in this context. In Section 3 we discuss horosphere intersections in more detail, in particular the quantitative relation between levels and heights of horospheres for integral lattices. In Sections 4 and 5 we apply this method to compute presentations for the Picard and Hurwitz modular groups respectively. 

We would like to thank Daniel Allcock for suggesting this method and for many helpful comments, Matthew Stover for pointing out a mistake in an earlier version of the paper and Martin Deraux for helpful discussions. We would also like to thank the referee for numerous comments and suggestions which greatly improved the paper; in particular their suggestion to automate the computational steps was critical to obtaining certifiably correct presentations for the final version of the paper.




\section{Horoball coverings and lattice presentations}

\subsection{Adapted horoball coverings and covering complex}\label{horoballs}

Let $X$ be a negatively curved symmetric space, i.e. 
a hyperbolic space ${\rm H}_\K^n$, with $\K=\R, \C, \quat$ or $\oct$ (and $n\geqslant 2$ if $\K=\R$, $n=2$ if $\K=\oct$). We refer the reader to \cite{CG} for general properties of these spaces and their isometry groups. In particular isometries of such spaces are roughly classified into the following 3 types: \emph{elliptic} (having a fixed point in $X$), \emph{parabolic} (having no fixed point in $X$ and exactly one on $\partial_\infty X$) or \emph{loxodromic} (having no fixed point in $X$ and exactly two on $\partial_\infty X$). 

Let $\Gamma$ be a lattice in ${\rm Isom}(X)$; the well-known Godement compactness criterion states that $\Gamma$ contains parabolic isometries if and only if it is non-cocompact, which we now assume. A \emph{cusp point} of $\Gamma$ is a point of $\partial_\infty X$ fixed by a parabolic element of $\Gamma$; a \emph{cusp group} of $\Gamma$ is a subgroup of the form ${\rm Stab}_\Gamma(p)$ where $p \in \partial_\infty X$ is a cusp point of $\Gamma$.

Assume that we are given a $\Gamma$-invariant covering of $X$ by (open) horoballs (see Definition~\ref{horo}), i.e. a collection $\mathcal{B}$ of horoballs such that:

$$\left\lbrace\begin{array}{l}
\gamma B \in \mathcal{B} \ {\rm for \ all} \ \gamma \in \Gamma \ {\rm and} \ B \in \mathcal{B} \\
\bigcup_{B \in \mathcal{B}} B = X
\end{array}\right.
$$

We will moreover assume that each horoball $B \in \mathcal{B}$ is based at a cusp point of $\Gamma$, and that each cusp point of $\Gamma$ is the basepoint of a unique horoball in $\mathcal{B}$ (giving a bijection between cusp points of $\Gamma$ and horoballs in $\mathcal{B}$); we will call such a covering $\mathcal{B}$ a \emph{$\Gamma$-adapted horoball covering}. Since the lattice $\Gamma$ has only finitely many cusp points modulo the action of $\Gamma$, it follows that such a horoball covering is a finite union of $\Gamma$-orbits of horoballs. 

Given a $\Gamma$-adapted horoball covering $\mathcal{B}$, the 
\emph{covering complex} $\mathcal{C}(\mathcal{B})$ associated to $\mathcal{B}$ is the simplicial 2-complex with vertex set $\mathcal{B}$, with an edge connecting each pair of vertices $B_1$ and $B_2$ such that $B_1 \cap B_2 \neq \emptyset$, and a triangle for each triple of vertices $B_1, B_2, B_3$ such that  $B_1 \cap B_2 \cap B_3 \neq \emptyset$. This the 2-skeleton of a simplicial complex sometimes called the nerve of the covering. By the above remark the quotient of the covering complex by the action of $\Gamma$ is a finite simplicial 2-complex.


We will use the following classical result of Macbeath (\cite{M}): 
\begin{thm}[\cite{M}]\label{macbeath} Let $\Gamma$ be a group acting by homeomorphisms on a topological space $X$. Let $V$ be an open subset of $X$ whose $\Gamma$-translates cover $X$.\\
(1) If $X$ is connected then the set  $E(V)=\{ \gamma \in \Gamma \ | \ V \cap \gamma V \neq \emptyset \}$ generates $\Gamma$.\\
(2) If moreover $X$ is simply-connected and $V$ path-connected, then $\Gamma$ admits a presentation with generating set $E(V)$ and relations $\gamma \cdot \gamma' = \gamma \gamma'$ for all $\gamma, \gamma' \in E(V)$ such that $V \cap \gamma V  \cap \gamma \gamma' V \neq \emptyset$.
\end{thm}
We kept the notation from~\cite{M}; to clarify the notation in part (2), let $S=\{ e_\gamma \, | \, \gamma \in E(V) \}$ be a set labelled by the generating set $E(V)$ from part (1). The claim is that $\Gamma$ has presentation $\la S \, | \, R \ra$, where $R$ consists of the relations $e_\gamma \cdot e_{\gamma'}=e_{\gamma \gamma'}$ for all $\gamma,\gamma ' \in E(V)$ such that $V \cap \gamma V  \cap \gamma \gamma' V \neq \emptyset$. Note that the latter condition implies that $\gamma \gamma' \in E(V)$.

Now if as above $\Gamma$ is a lattice in ${\rm Isom(X)}$ and $\mathcal{B}$ is a $\Gamma$-adapted horoball covering of $X$, we may as remarked above write $\mathcal{B}$ as a finite union of $\Gamma$-orbits of horoballs $B_1,...,B_k$ (say, minimally). One can then apply Macbeath's theorem with $V=B_1 \cup ... \cup B_k$, after possibly enlarging each horoball $B_i$ in order for this union to be (path-)connected. For simplicity of exposition, we henceforth asume that $\Gamma$ has a single cusp, so that the $\Gamma$-adapted horoball covering consists of a single $\Gamma$-orbit of horoballs (this is the case in all examples considered in this paper). In that case the process of obtaining a presentation fom the covering complex is closely related to a complex of groups structure on the quotient of the covering complex, the only difference being that we need to take into account non-trivial edge and face stabilizers.

\subsection{Levels and proximal cusp complex}\label{levels}

Recall that if $X$ is a 
hyperbolic space ${\rm H}_\K^n$ (with $\K=\R, \C$ or $\quat$) 
then $X$ admits the following \emph{projective model} which we briefly recall.

Consider $\K^{n,1}$, the vector space $\K^{n+1}$ endowed with a Hermitian form $\langle \cdot \, , \cdot \rangle$ of signature $(n,1)$. (When $\K=\quat$ we will use the convention that scalars act on vectors on the right, whereas matrices act on vectors on the left.)
Let $V^-=\left\lbrace Z \in \K^{n,1} | \langle Z , Z \rangle <0 \right\rbrace$, $V^0=\left\lbrace Z \in \K^{n,1} | \langle Z , Z \rangle =0 \right\rbrace$ and let $\pi: \K^{n+1}-\{0\} \longrightarrow \K{\rm P}^n$ denote projectivization.
One then defines ${\rm H}_\K^n$ to be $\pi(V^-) \subset \K{\rm P}^n$, endowed with the distance $d$ (Bergman metric) given by, for $Z,W \in V^-$:

\begin{equation}\label{dist}
  \cosh ^2 \Bigl(\frac{d(\pi(Z),\pi(W))}{2}\Bigr) = \frac{|\langle Z, W \rangle|^2}{\langle Z, Z \rangle  \langle W, W \rangle}.
\end{equation}
Note that the right-hand side is independent of the choice of lifts $Z,W$. Then ${\rm Isom}^0(X)={\rm PU}(n,1,\K)$, the (projectivization of) the matrix group preserving the Hermitian form (see \cite{CG}). Note that  ${\rm PU}(n,1,\K)$ is usually denoted ${\rm PO}(n,1)$ when $\K=\R$, ${\rm PU}(n,1)$ when $\K=\C$ and ${\rm PSp}(n,1)$ when $\K=\quat$. The boundary at infinity $\partial_\infty X$ is then identified with  $\pi(V^0) \subset \K{\rm P}^n$. We would like to measure distances between points of $\partial_\infty X$ using the Hermitian form as in (\ref{dist}); one way to do this is to use integral lifts of vectors with rational coordinates as follows.

We now assume that $\Gamma$ is an \emph{integral lattice} in the sense that it is contained in ${\rm U}(H,\mathcal{O}_E)$ for some number field $E$ (or finite-degree division algebra over $\Q$ when $\K=\quat$) with ring of integers $\mathcal{O}_E$, and Hermitian form $H=\langle \cdot \, , \cdot \rangle$ defined over $E$. 
We say that an integral vector $P_0=(p_1,...,p_{n+1}) \in \mathcal{O}_E^{n+1}$ is \emph{primitive} if it has no integral submultiple in the following sense: if $P_0 \lambda ^{-1}\in \mathcal{O}_E^{n+1}$ for some $\lambda \in \mathcal{O}_E$ then $\lambda$ is a unit in $\mathcal{O}_E$.

If $p$ is an $E$-rational point in $\K P^{n}$, i.e. the projective image of a vector $P=(p_1,...,p_{n+1}) \in E^{n+1}$, a \emph{primitive integral lift} of $p$ is any lift $P_0$ of $p$ to $\mathcal{O}_E^{n+1}$ which is a primitive integral vector.

\begin{lem}\label{uniquelift} If $\mathcal{O}_E$ is a principal ideal domain then primitive integral lifts are unique up to multiplication by a unit.
\end{lem}

\begin{lem}\label{primitivecolumn} (a) Any column-vector of a matrix $A \in {\rm U}(H,\mathcal{O}_E)$ is a primitive integral vector. (b)  If moreover $E$ is imaginary quadratic with $\mathcal{O}_E$ a principal ideal domain and one of the standard basis vectors $B_i$ is $H$-isotropic, then for any $H$-isotropic primitive integral vector $V$ and $A \in {\rm U}(H,\mathcal{O}_E)$, $AV$ is a primitive integral vector.
\end{lem}

\Pf (a) Let $A \in {\rm U}(H,\mathcal{O}_E)$ and $V$ a column-vector of $A$. Then $V$ is integral; assuming that it is not primitive, there would exist a non-unit $\lambda \in \mathcal{O}_E$ such that $V\lambda^{-1}$ is also integral. But then the matrix $A'$ obtained from $A$ by replacing the column-vector $V$ by $V\lambda^{-1}$ would also be in ${\rm GL}(n+1,\mathcal{O}_E)$, with ${\rm det}\, A'={\rm det}\, A\lambda^{-1}$, a contradiction since the latter is not an integer, as ${\rm det}\, A$ is a unit and $\lambda$ is not. \\
(b) Let $V$ be an $H$-isotropic primitive integral vector and $A \in {\rm U}(H,\mathcal{O}_E)$. If $\mathcal{O}_E$ is a principal ideal domain then ${\rm PU}(H,\mathcal{O}_E)$ has a single cusp (see \cite{Zi}), therefore there exists $M \in {\rm U}(H,\mathcal{O}_E)$ mapping $B_i$ to $AV\lambda$ for some $\lambda \in E$. Then as in (a) $\lambda$ must be a unit, hence $AV$ is a column-vector of $M\lambda^{-1} \in {\rm U}(H,\mathcal{O}_E)$ and we conclude by (a). \EPf

\begin{dfn}\label{level} Given two $E$-rational points $p,q \in \partial_\infty X$, the \emph{level between $p$ and $q$}, denoted ${\rm lev}(p,q)$, is $| \langle P_0, Q_0 \rangle |^2$ for any two primitive integral lifts $P_0,Q_0$ of $p,q$ respectively. When we are given a preferred $E$-rational point $\infty \in \partial_\infty X$, the \emph{depth} of an $E$-rational point $p \in \partial_\infty X$ is the level between $p$ and $\infty$.
\end{dfn}

By Lemma~\ref{uniquelift} this is well-defined when $\mathcal{O}_E$ is a principal ideal domain. 
The \emph{proximal cusp complex of level $n$}, denoted $\mathcal{C}_n(\Gamma)$, is the complex whose vertices are cusp points of $\Gamma$, with an edge connecting 2 vertices $p,q$ whenever ${\rm lev}(p,q) \leqslant n$, and a triangle for each triple of distinct edges.

Levels give a convenient way to distinguish orbits of edges and triangles in the covering complex, by the following observation which follows from Lemmas~\ref{uniquelift} and \ref{primitivecolumn}:

\begin{lem} If $\mathcal{O}_E$ is a principal ideal domain, for any two $E$-rational points $p,q \in \partial_\infty X$ and $\gamma \in {\rm U}(H,\mathcal{O}_E)$, $\rm{lev}(\gamma p, \gamma q)={\rm lev}(p,q)$.
\end{lem}

More importantly, levels allow us to find the optimal height $u$ of a horosphere $H_u=\partial B_u$ based at a preferred $E$-rational point $\infty \in \partial_\infty X$ such that the orbit $\Gamma B_u$ covers $X$. This relies on the following result, which is part of Corollary~\ref{cordepthradius} in Section~3:


\begin{lem}\label{depth-height} There exists a decreasing function $u: \N \longrightarrow \R$ such that, for any $E$-rational point $p \in \partial_\infty X$ with depth $n$, and any (integral) $A_p \in \Gamma$ satisfying $A_p(\infty)=p$, the set $H_u \cap A_p(H_u)$ is empty if and only if $u > u(n)$.
\end{lem}    

In fact we will see in Corollary~\ref{cordepthradius} that the function $u$ is given by $u(n)=\frac{2}{\sqrt{n}}$.

\begin{dfn} The \emph{covering depth} of $\Gamma$ is the unique $n\in \N$ such that $\frac{2}{\sqrt{n+1}} < u^{cov} \leqslant \frac{2}{\sqrt{n}}$, where $u^{cov}$ denotes the maximal height such that $\Gamma B_{u^{cov}}$ covers $X$.
\end{dfn}

In practice, for the purpose of finding a presentation of $\Gamma$, we will not need to explicitly determine the covering depth or $u^{cov}$. It will suffice to bound the covering depth from above, and use the covering of $X$ at the corresponding height to apply Macbeath's theorem.   

\subsection{Reduction modulo the vertex stabilizer $\Gamma_\infty$}

We choose a preferred cusp point $\infty \in \partial_\infty X$ of $\Gamma$ (in general we will take $\infty=\pi([1,0,...,0]^T)$ in the Siegel model, see section~\ref{horosiegel}), and consider the cusp stabilizer $\Gamma_\infty = {\rm Stab}_\Gamma(\infty)$. Since $\Gamma$ is a lattice, it is well known that $\Gamma_\infty$ acts cocompactly on all horospheres based at $\infty$. Let $n$ denote the covering depth of $\Gamma$ and $u^{cov}$ the corresponding covering height, so that the $\Gamma$-translates of the horoball $B_u$ cover $X$, and let $D_\infty \subset H_u$ be a compact fundamental domain for the action of $\Gamma_\infty$ on $H_u \simeq \partial_\infty X \setminus \{ \infty \}$. In practice we will choose $D_\infty$ to be an affinely convex polytope in Heisenberg coordinates (see section~\ref{horosiegel}).

Assume that we are given a finite presentation $\Gamma_\infty= \langle S_\infty \vert R_\infty \rangle$. Then we may reduce the procedure in Macbeath's theorem to finitely many additional generators and relations as follows. 
Let $\{p_1,...,p_k\}$ denote the $E$-rational points with depth at most $n$ in $D_\infty$, and assume for simplicity that they are ordered in such a way that the first $r$ of them form a system of representatives under the action of $\Gamma_\infty$. Assume moreover that we have found for each $i=1,...,r$ an element $A_i \in \Gamma$ such that $A_i(\infty)=p_i$ (this is possible in principle since $\Gamma$ is assumed to have a single cusp).

\medskip

{\bf Generators:} The group $\Gamma$ is generated by $\{ S_\infty,A_1,...,A_r \}$. This follows easily from part (1) of Macbeath's theorem and Lemma~\ref{depth-height}, as any $E$-rational point of $\partial_\infty X$ with depth at most $n$ is in the $\Gamma_\infty$-orbit of one of $p_1,...,p_r$. Note that with the notation from the Theorem, we are using the open set $V=B_u=B$ to cover $X$, and  $E(B)=\{ \gamma \in \Gamma \, | \, B \cap \gamma B \neq \emptyset \}= \Gamma_\infty \{ A_1,...,A_r \} \Gamma_\infty =\{ \gamma_\infty^1 A_i \gamma_\infty^2\, | \, \gamma_\infty^1,\gamma_\infty^2 \in \Gamma_\infty, i=1,...,r \}$. Indeed, by Lemma~\ref{depth-height}, $B \cap \gamma B \neq \emptyset$ if and only if $\gamma \infty$  is either $\infty$ or an $E$-rational point of depth at most $n$, which is a $\Gamma_\infty$-translate of one of $p_1,...,p_r$. 

\medskip

{\bf Relations:} We now rephrase part (2) of Macbeath's theorem in this context. Let $\gamma, \gamma' \in E(B)$ satisfy $B \cap \gamma B \cap \gamma \gamma' B \neq \emptyset$, and first assume that $\gamma \infty \neq \infty$ and $\gamma \gamma' \infty \neq \infty$. After conjugating by an element of $\Gamma_\infty$ we may assume that $\gamma=A_a$, $\gamma'=\gamma_\infty^1 A_b \gamma_\infty^2$, $\gamma \gamma'=\gamma_\infty^3 A_c \gamma_\infty^4$ for some $a,b,c \in \{ 1,...,r \}$ and $\gamma_\infty^1,...,\gamma_\infty^4 \in \Gamma_\infty$. The corresponding relation $\gamma \gamma'=\gamma \cdot \gamma'$ is then: $A_a \gamma_\infty^1 A_b \gamma_\infty^2=\gamma_\infty^3 A_c \gamma_\infty^4$. Taking the image of $\infty$ under both sides of this relation gives: $A_a(\gamma_\infty^1 p_b)=\gamma_\infty^3 p_c$. 

In practice this is how we will detect the relations, finding which points of depth at most $n$ are sent to points of depth at most $n$ by the generators $A_a$. One then recovers the relation as follows. For each triple $(a,b,c)$ for which there exist $\gamma_\infty^1,\gamma_\infty^3 \in\Gamma_\infty$ such that $A_a (\gamma_\infty^1 p_b) = \gamma_\infty^3 p_c$, we obtain a relation $R_{a,b,c}$ by identifying the element $A_c^{-1}(\gamma_\infty^3)^{-1}A_a\gamma_\infty^1A_b \in \Gamma_\infty$ as a word in the generators $S_\infty$.  
The key point is that there are only finitely many such triples $(a,b,c)$ to be checked, by Lemma~\ref{cyganbound}.

Now assume that one of $\gamma \infty$, $\gamma \gamma' \infty$ is $\infty$ but not both (as the relations in $\Gamma_\infty$ have already been considered). The corresponding relation can be obtained as above, using the point $p_\infty=\infty$ with corresponding group element $A_\infty = {\rm Id}$.

Summarizing the above discussion gives:

\begin{lem}\label{pres2} With the above notation, $\Gamma$ admits the presentation
 $\Gamma=\langle S_\infty,A_1,...,A_r \vert R_\infty, R_{a,b,c} \rangle$.
\end{lem}
 
\subsection{The method in practice}\label{method}

We now give an outline of the method we use to apply Macbeath's theorem:)
)
\begin{itemize}
\item[(1)] Find an explicit (affine) fundamental domain $D_\infty \subset \partial_\infty X \simeq H_u$ for the action of $\Gamma_\infty={\rm Stab}_\Gamma(\infty)$, and a presentation $\Gamma_\infty= \langle S_\infty \vert R_\infty \rangle$.
\item[(2)] Find the covering depth $n$ of $\Gamma$. 
  Consider the corresponding covering complex $\mathcal{C}(\Gamma B_{u^{cov}}) \simeq \mathcal{C}_n(\Gamma)$.
\item[(3)] Find all $E$-rational points $\{p_1,...,p_d\}$ in $D_\infty$ with depth at most $n$ (and denote $p_\infty=\infty$).
\item[(4)] For each of the $r$ $\Gamma_\infty$-orbits of points $p_a$, find an explicit $A_a \in \Gamma$ such that $A_a(\infty)=p_a$.
\item[(5)] For each triple $(a,b,c)$ for which there exist $\gamma_\infty^1,\gamma_\infty^3 \in\Gamma_\infty$ such that $A_a (\gamma_\infty^1 p_b) = \gamma_\infty^3 p_c$ 
we obtain a relation $R_{a,b,c}$ by identifying the element $A_c^{-1}(\gamma_\infty^3)^{-1}A_a\gamma_\infty^1A_b \in \Gamma_\infty$ as a word in the generators $S_\infty$. 

To find all such triples, we seek triple intersections of the form described above.  A necessary condition for the triple intersection $B\cap A_a B\cap \gamma_\infty^3 A_c B$ to be nonemepty is for all three pairwise intersections to be nonempty.  We get two of them, $B \cap \gamma_\infty^3 A_c B\neq\emptyset$ and $B \cap A_a B\neq\emptyset$, for free since $A_a$ and $A_c$ are generators.  For each of the finitely many $A_a$'s and $A_c$'s, we want to find all $\gamma_\infty^3$ such that the third pairwise intersection is nonempty.

The third intersection $A_a B\cap \gamma_\infty^3 A_c B$ will be nonempty if and only if the level between the corresponding points $p_a$ and $\gamma_\infty^3 p_c$ is at most the covering depth $n$. By Lemma~\ref{cyganbound}, this translates to a bound on the Cygan distance between these points, which by the triangle inequality gives a bound on the Cygan distance of each point to the reference point $(0,0)$. As there are finitely many points $\{p_1\ldots, p_r\}$, we obtain an upper bound on these distances, 
and cover a ball centered at (0,0) of the corresponding radius by translates of $D_\infty$.  The finite collection of translates required to do this is our finite list of candidates for $\gamma_\infty^3$.  For each such $\gamma_\infty^3$, we iterate over the $A_a$'s and $A_c$'s, checking the level between $p_a$ and $\gamma_\infty^3 p_c$ to find all intersections.

Finally, we take all the pairs $(a,c)$ for which the three pairwise intersections are nonempty, and compute coordinates for the point $A_a^{-1}\gamma_\infty^3 A_c p_\infty$.  If the depth of this point is greater than $n$, it will not give us a relation so we discard it.  If it is less than $n$, then it is of the form $\gamma_\infty^1 p_b$.  We can easily identify $\gamma_\infty^1$ and $b$ by figuring out which translate of $D_\infty$ the point lives in.  The last step is to identify the group element $A_b^{-1}(\gamma_\infty^1)^{-1}A_a^{-1}\gamma_\infty^3 A_c$ as a word in the generators of $\Gamma_\infty$.

\end{itemize}

Then by Macbeath's theorem and Lemma~\ref{pres2}, $\Gamma=\langle S_\infty,A_1,...,A_r \vert R_\infty, R_{a,b,c} \rangle$. Note that in step (5) we don't in fact check which triples of horoballs have nonempty triple intersection (which is harder to check); rather we write a relation for each triple satisfying the necessary condition arising from Lemma~\ref{cyganbound}.
In principle this may yield redundant relations which we later eliminate.  

Step (5) was largely carried out by a computer.  Using the bound on the level coming from the Cygan distance (see Lemma~\ref{cyganbound}), we found a Cygan ball containing all centers of horoballs that could possibly have nonempty intersections with any of $\{A_1B,\ldots, A_rB \}$.  We then covered this ball by $\Gamma_\infty$-translates of $D_\infty$, and iterated over the list of centers to generate a list of cycles.  We then compute the matrix representation of $A_c^{-1}(\gamma_\infty^3)^{-1}A_a\gamma_\infty^1A_b$ and attempt to identify it as a word in $\Gamma_\infty$.  For $d=1,3,7$ we did this very inefficiently, by iterating over words in $\Gamma_\infty$ of increasing length until we found it.  This way of searching would have taken unreasonably long for the quaternions, so instead we devised a way of guessing the element based on distance, and then making small corrections to get it exactly right. See \cite{MCode} for the details of the computations.

\smallskip

In order to avoid tedious repetition of similar arguments or straightforward computations, we will  only give one detailed proof for each step for the Picard modular groups; we will give detailed arguments for the quaternionic Hurwitz lattice. We will usually choose the most difficult case, or the most instructive if the various cases are of similar difficulty. Step (1) is routine and we just state the results, except for the Hurwitz lattice, where we find a presentation for $\Gamma_\infty$ in Lemma~\ref{inftypres} and use Philippe's fundamental domain found in \cite{Ph} (see Proposition~\ref{philippedomains}). We give a detailed argument and proof for step (2) for the Picard modular group $\Gamma(7)$ in Lemma~\ref{d=7cover}, and for step (3) for the same group and depth 2 in Lemma~\ref{depth2points}. 

There seems to be no general strategy for step (4); we find all relevant matrices in this paper by combining two tricks, which luckily cover all the cases we need. The first trick is to use stabilizers of vertical complex lines in the Heisenberg group: it is easy to find such a matrix when it stabilizes the vertical axis, then we carry over to other vertical lines by conjugating by a horizontal translation. The second trick is to hit all relevant integral points by the group elements that we already know, and see if we land in the $\Gamma_\infty$-orbit of the point we are trying to reach. 


\subsection{A toy example: $\Gamma={\rm PSL(2,\Z)}$}\label{toyexample}

In order to illustrate the method, we now go through its steps for $\Gamma={\rm PSL(2,\Z)}$ exactly as we will for the more complicated Picard and Hurwitz modular groups. The results are either well-known or elementary and we state them without proof.

{\bf Presentation and fundamental domain for the cusp stablilizer $\Gamma_\infty$:} The cusp stabilizer $\Gamma_\infty$ has presentation $\langle T \rangle$; a fundamental domain for its action on $\partial {\rm H}^2_\R \setminus \{ \infty \} \simeq \R$ is $D_\infty=[0,1]$. Concretely we use the following generator for $\Gamma_\infty$:

$$T=\left[ \begin{array}{cc}
1 & 1 \\
0 & 1
\end{array}\right]
$$

{\bf Covering depth and $\Q$-rational points in $D_\infty$:} 
The covering depth of ${\rm PSL}(2,\Z)$ is 1. The $\Q$-rational points of depth 1 in $D_\infty$ are 0 and 1, both in the same $\Gamma_\infty$-orbit. An integral lift of 0 is $p_0=[0,1]^T$; we denote $p_1=Tp_0$.

In the Picard modular and Hurwitz cases this is a more difficult estimation, owing to the fact that we cannot simply look at a 2-dimensional picture.  In the Picard modular case we can look at a 3-dimensional picture to guess at the covering depth which we may then verify, but in the Hurwitz case it is much more challenging to  use a picture for help; see Figure~\ref{cubes} for some visual intuition in that case.  

{\bf Generators:} 
The following element $A_0 \in \Gamma$ maps the point $\infty=[1,0]^T$ to $p_0$:

$$A_0=I_0=\left[ \begin{array}{cc}
0 & -1 \\
1 & 0
\end{array}\right]
$$

In the Picard modular and Hurwitz cases we must find generators mapping $\infty$ to points of depth greater than 1.  While there is no general method for doing this, in practice we were always able to do it by combining the two tricks described at the end of the previous section.

{\bf Relations:} 
We list in Table~1 the relations obtained for $\rm{PSL}(2,\Z)$ by applying generators to points of depth at most 1 as described in part (5) of section~\ref{method}. The second relation is obtained by following the corresponding cycle of points, which gives $I_0TI_0TI_0 \in \Gamma_\infty$. The latter element is computed to be $T^{-1}$, giving the relation $(I_0T)^3={\rm Id}$.

For each matrix $A$ and point $p$ for which $A.p$ has depth less than or equal to the covering depth, we obtain a relation from the following cycle of points
$$\infty \xrightarrow{A_p} p \xrightarrow{W^{-1} A} p' \xrightarrow{A_{p'}^{-1}}  \infty$$
We write $A.p = W.p'$ where $W$ is an element of $\Gamma_\infty$ written as a word in its generators.  Then $A_{p'}^{-1}W^{-1}A A_p$ is an element of $\Gamma_\infty$, where $A_p$ and $A_{p'}$ are our chosen generators taking $\infty$ to $p$ and $p'$ respectively.  We write that element as a word $W'$ in the generators of $\Gamma_\infty$, and the relation we obtain from the cycle is $A_{p'}WA A_p = W'$

\begin{table}[h]\label{SL2table}
  \begin{center}
  {\renewcommand{\arraystretch}{1.2}%
\begin{tabular}{|c|c|c|c|}
\hline 
$\mathbf{A.p}$ & $\mathbf{p'}$ & $\mathbf{A_{p'}^{-1}W^{-1}A A_p}$ & $\mathbf{W'}$  \\
\hline
$I_0 p_0$ & $\infty$ & $I_0^2$ & ${\rm Id}$\\
\hline
$I_0  p_1$ & $p_0$ & $(TI_0)^{-1}(TI_0)^{-1}I_0$ & $T$\\
\hline
\end{tabular}
}
\end{center}
\caption{Action of generators on vertices for ${\rm PSL}(2,\Z)$}
\end{table}

This step works in exactly the same way in the Picard modular and Hurwitz cases, but with more points and more cycles.

\section{Horosphere intersections}\label{horosiegel}

Our main reference for this section is \cite{KP}. We will use the \emph{Siegel model} of hyperbolic space ${\rm H}_\K^n$ (with $\K=\R, \C, \quat$), which is the projective model (as described in Section~\ref{levels}) associated to the Hermitian form on $\K^{n+1}$ given by $\langle Z,W \rangle = W^* J Z$ with:

$$ J=\left(\begin{array}{ccc}
0 & 0 & 1 \\
0 & I_{n-1} & 0 \\
1 & 0 & 0 \end{array}\right)
$$
Then hyperbolic space ${\rm H}_\K^n$ can be parametrized by $\K^{n-1} \times {\rm Im} \, \K \times \R^+$ as follows, denoting as before by $\pi$ the projectivization map: ${\rm H}_\K^n=\{\pi(\psi(\zeta,v,u)) \, | \, \zeta \in \K^{n-1}, v \in {\rm Im}\,  \K, u \in \R^+)\}$, where: 

\begin{eqnarray}\label{horocoord}
\psi(\zeta,v,u)= \left(\begin{array}{c}
(-|\zeta|^2-u+v)/2 \\
\zeta \\
1
\end{array}\right)
\end{eqnarray}
With this parametrization the boundary at infinity $\partial_\infty {\rm H}_\K^n$ corresponds to the one-point compactification:
$$\left\{\pi(\psi(\zeta,v,0)) \, | \, \zeta \in \K^{n-1}, v \in {\rm Im}\,  \K \right\} \cup \{\infty\}$$
where $\infty=\pi((1,0,...,0)^T)$. The coordinates $(\zeta,v,u) \in \K^{n-1} \times {\rm Im}\, \K \times \R^+$ are called the \emph{horospherical coordinates} of the point $\pi(\psi(\zeta,v,u)) \in {\rm H}_\K^n$. 

\begin{dfn}\label{horo} For a fixed $u_0 \in \R^+$, the level set $H_{u_0}=\{\pi(\psi(\zeta,v,u_0)) \, | \, \zeta \in \K^{n-1}, v \in {\rm Im}\,  \K \}$ is called the \emph{horosphere at height $u_0$ based at $\infty$}, and $B_{u_0}=\{\pi(\psi(\zeta,v,u)) \, | \, \zeta \in \K^{n-1}, v \in {\rm Im}\,  \K, u > u_0 \}$ is called the \emph{horoball at height $u_0$ based at $\infty$}.
\end{dfn}

The punctured boundary  $\partial_\infty {\rm H}_\K^n \setminus \{ \infty \}$ is then naturally identified to the \emph{generalized Heisenberg group} ${\rm Heis}(\K,n)$, defined as the set $\K^{n-1} \times {\rm Im} \, \K$ equipped with the group law:
$$(\zeta_1,v_1)(\zeta_2,v_2)=(\zeta_1+\zeta_2,v_1+v_2+2{\rm Im} \, (\overline{\zeta_2} \cdot \zeta_1))
$$
where $\cdot$ denotes the usual Euclidean dot-product on $\K^{n-1}$. This is the classical 3-dimensional Heisenberg group when $\K=\C$ and $n=2$. The identification of  $\partial_\infty {\rm H}_\K^n \setminus \{ \infty \}$ with ${\rm Heis}(\K,n)$ is given by the simply-transitive action of ${\rm Heis}(\K,n)$ on $\partial_\infty {\rm H}_\K^n \setminus \{ \infty \}$, where the element $(\zeta_1,v_1) \in {\rm Heis}(\K,n)$ acts on the vector $\psi(\zeta_2,v_2,0)$ by left-multiplication by the following \emph{Heisenberg translation} matrix in ${\rm U}(n,1,\K)$:
\begin{eqnarray}\label{heistrans}
T_{(\zeta_1,v_1)}=\left(\begin{array}{ccc}
1 & -\zeta_1^* & (-|\zeta_1|^2+v_1)/2 \\
0 & {\rm I}_{n-1} & \zeta_1 \\
0 & 0 & 1
\end{array}\right)
\end{eqnarray} 
Given an element $U \in {\rm U}(n,\K)$, the \emph{Heisenberg rotation} by $U$ is given by the following matrix:
\begin{eqnarray}\label{heisrot}
R_U=\left(\begin{array}{ccc}
1 & 0 & 0 \\
0 & U & 0 \\
0 & 0 & 1 
\end{array}\right)
\end{eqnarray}

There is an additional class of isometries fixing $\infty$ when $\K=\quat$, coming from the action of diagonal matrices which is non-trivial in the non-commutative case. Recall that our convention is that matrices act on vectors on the left, and scalars act on vectors on the right. Then, for any unit quaternion $q \in \quat$,  the diagonal matrix $C_q=q\,{\rm Id}$ acts by the isometry of hyperbolic space given by conjugating horospherical coordinates (the result of multiplying the vector form (\ref{horocoord}) by $q$ on the left, then normalizing by $q^{-1}$ on the right):

\begin{equation}\label{conj}
C_q:(\zeta,v,u) \longmapsto (q\zeta q^{-1},qvq^{-1},u)
\end{equation}  

For this reason, when $\K=\quat$ the relevant projectivization of ${\rm U}(n,1,\quat)$ acting on ${\rm H}_\quat^n$ is ${\rm PU}(n,1,\quat)={\rm U}(n,1,\quat)/\{\pm{\rm Id} \}$ rather than ${\rm U}(n,1,\quat)/{\rm U}(1)$.

Heisenberg translations and rotations, as well as conjugation by unit quaternions, preserve the following distance function on ${\rm Heis}(\K,n)$, called the \emph{Cygan metric}, defined for $(\zeta_1,v_1), (\zeta_2,v_2) \in {\rm Heis}(\K,n)$ by:
\begin{eqnarray}\label{cygan} 
d_C((\zeta_1,v_1), (\zeta_2,v_2)) & = & \left\vert |\zeta_1-\zeta_2|^4 +  \left\vert v_1-v_2+2{\rm Im}\, (\overline{\zeta_2} \cdot \zeta_1)\right\vert^2 \right\vert^{1/4} \\
 & = & \left\vert 2 \langle \psi(\zeta_1,v_1,0),\psi(\zeta_2,v_2,0 )   \rangle \right\vert ^{1/2} 
\end{eqnarray} 

This is in fact the restriction to $\partial_\infty {\rm H}_\K^n \setminus \{ \infty \}$ of an incomplete distance function on $\overline{{\rm H}_\K^n}\setminus \{ \infty \}$ called the \emph{extended Cygan metric} (see \cite{KP}), defined for $(\zeta_1,v_1,u_1), (\zeta_2,v_2,u_2) \in \K^{n-1} \times {\rm Im} \, \K \times \R^{\geqslant 0} \simeq \overline{{\rm H}_\K^n} \setminus \{ \infty \}$ by:
\begin{eqnarray}\label{extcygan} 
d_{XC}((\zeta_1,v_1,u_1), (\zeta_2,v_2,u_2)) & = & \left\vert \left(|\zeta_1-\zeta_2|^2+|u_1-u_2|\right)^2 +  \left\vert v_1-v_2+2{\rm Im}\, (\overline{\zeta_2}\cdot \zeta_1)\right\vert^2 \right\vert^{1/4} \\
 & = & \left\vert 2 \langle \psi(\zeta_1,v_1,u_1),\psi(\zeta_2,v_2,u_2) \rangle \right\vert ^{1/2} 
\end{eqnarray} 
We define \emph{Cygan spheres}, \emph{Cygan balls}, \emph{extended Cygan spheres} and \emph{extended Cygan balls} in the usual way relative to these distance functions. 

When $\Gamma < {\rm U}(n,1,\mathcal{O}_E)$ is an integral lattice with $\mathcal{O}_E$ a principal ideal domain as in section 2.2, the Cygan distance relates to levels and depths of integral boundary points as follows:

\begin{lem}\label{cyganbound} Let $g=(g_{i,j})$ and $h=(h_{i,j}) \in {\rm U}(n,1,\mathcal{O}_E)$ satisfy $g\infty \neq \infty, h\infty \neq \infty$, and denote $(\zeta_g,v_g),(\zeta_h,v_h)$ the horospherical coordinates of $g\infty,h\infty$   respectively. Then:
$$ d_C \left(  (\zeta_g,v_g),(\zeta_h,v_h)  \right) = \left( \frac{4{\rm lev}(g\infty,h\infty)}{{\rm depth}(g\infty){\rm depth}(h\infty)} \right) ^{1/4}.
$$
In particular, given $n \geqslant 1$:
$$ {\rm lev}(g\infty,h\infty) \leqslant n \iff d_C \left(  (\zeta_g,v_g),(\zeta_h,v_h)  \right) \leqslant \left( \frac{4n}{{\rm depth}(g\infty){\rm depth}(h\infty)} \right) ^{1/4}.
$$
\end{lem}

\Pf By Lemma~\ref{primitivecolumn}, the first column vector of $g$ (resp. $h$) is a primitive integral lift of $g\infty$ (resp. $h\infty$), and it can be written in terms of horospherical coordinates as $g_{n+1,1}\psi(\zeta_g,v_g)$ (resp. $h_{n+1,1}\psi(\zeta_h,v_h)$). Therefore, using (\ref{cygan}) we have:

\begin{eqnarray*} {\rm lev}(g\infty,h\infty) & = & \la g_{n+1,1}\psi(\zeta_g,v_g),h_{n+1,1}\psi(\zeta_h,v_h) \ra \\
& = & |g_{n+1,1}|^2 |h_{n+1,1}|^2\la\psi(\zeta_g,v_g),\psi(\zeta_h,v_h) \ra \\
& = & {\rm depth} (g\infty) {\rm depth} (h\infty) \frac{d_C \left(  (\zeta_g,v_g),(\zeta_h,v_h)  \right)^4}{4}. 
\end{eqnarray*} \EPf

When we apply Macbeath's theorem we argue that the images under $\Gamma$ of the horoball $B_u$ based at $\infty$ at a certain height $u>0$ cover $X$, or equivalently cover the horosphere $H_u=\partial B_u$.
The following result allows us to control the traces on $H_u$ of these images in terms of Cygan spheres depending only on arithmetic data. It involves \emph{Ford isometric spheres}, whose definition we first recall.

\begin{dfn}\label{isometricsphere}  The \emph{Ford isometric sphere} $I_g$ of an isometry $g \in {\rm U}(n,1,\K)$ is defined as 
$$I_g=\left\{ z=(\zeta,v,u) \in {\rm H}^n_\K \, | \, |\langle \psi (z),\psi (\infty)\rangle| =| \langle \psi (z),g^{-1}\psi (\infty)\rangle| \right\}$$
\end{dfn}  

By Proposition~4.3 of \cite{KP}, the Ford isometric sphere $I_g$ of $g=(g_{i,j}) \in {\rm U}(n,1,\K)$ is in fact the extended Cygan sphere $S$ with center $g^{-1}(\infty)$ and radius $\sqrt{2/|g_{n+1,1}|}$.

\begin{prop}\label{depthradius} Let $g=(g_{i,j}) \in {\rm U}(n,1,\K)$ satisfy $g(\infty)\neq \infty$, $S=I_{g^{-1}}$ the extended Cygan sphere with center $g(\infty)$ and radius $\sqrt{2/|g_{n+1,1}|}$, and $H_{u_0}$ the horosphere based at $\infty$ at height $u_0>0$. Then $H_{u_0} \cap g(H_{u_0})=H_{u_0} \cap S$.
\end{prop}

\Pf Using the following standard form for $g$ and $g^{-1}$, equation (1.2) of \cite{KP}, where $a,b,c,d \in \K$, $\alpha,\beta,\gamma,\delta \in \K^{n-1}$ and $A \in {\rm M}_{n-1}(K)$:

\begin{eqnarray}\label{standardform}
g=\left(\begin{array}{ccc} a & \gamma^* & b \\
\alpha & A & \beta \\
c & \delta^* & d 
\end{array}\right),
&
&
g^{-1}=\left(\begin{array}{ccc} \bar{d} & \beta^* & \bar{b} \\
\delta & A^* & \gamma \\
\bar{c} & \alpha^* & \bar{a} 
\end{array}\right),
\end{eqnarray}
from Equation~(\ref{horocoord}) and Definition~\ref{isometricsphere} we get that: $z=(\zeta,v,u) \in I_{g^{-1}}=S \iff \left| \frac{\bar{c}}{2}(-|\zeta|^2-u+v)+\alpha^*\zeta+\bar{a} \right|^2=1$.

Now fix $u_0>0$ and let $z=(\zeta,v,u_0) \in H_{u_0}$. We claim that: $g^{-1}(z) \in H_{u_0} \iff z \in S$. 

Indeed, using the above forms for $g^{-1}$ and $\psi(z)$ we have:

$$g^{-1}\psi(z)=\left( \begin{array}{c} \frac{\bar{d}}{2}(-|\zeta|^2-u_0+v)+\beta^*\zeta+\bar{b} \\
\frac{\delta}{2}(-|\zeta|^2-u_0+v)+A^*\zeta+\gamma \\
\frac{\bar{c}}{2}(-|\zeta|^2-u_0+v)+\alpha^*\zeta+\bar{a} 
\end{array}\right)=\left( \begin{array}{c} \zeta_1 \\ \zeta_2 \\ \zeta_3 \end{array}\right).
$$

The corresponding point is in the horosphere $H_{u_0}$ if and only if its $u$-coordinate in horospherical coordinates equals $u_0$. Now the $u$-coordinate of a point is recovered from any lift $(\zeta_1,\zeta_2,\zeta_3)^T$ by: 
$$u=-\zeta_1\zeta_3^{-1} - \overline{\zeta_1 \zeta_3^{-1}}- |\zeta_2 \zeta_3^{-1}|^2=-\zeta_3^{-1}\zeta_1 - \zeta_3^{-1}\overline{\zeta_1 \zeta_3^{-1}}\zeta_3 - |\zeta_2 \zeta_3^{-1}|^2=-|\zeta_3|^{-2}(\overline{\zeta_3}\zeta_1+ \overline{\zeta_1}\zeta_3 +|\zeta_2|^2).$$
(Note that we conjugated by $\zeta_3$ in the second step). Therefore: $g^{-1}(z) \in H_{u_0} \iff \overline{\zeta_3}\zeta_1+\overline{\zeta_1}\zeta_3 +|\zeta_2|^2 = -u_0 |\zeta_3|^2$. Expanding the left-hand side gives:
\begin{eqnarray*}
\left( (-|\zeta|^2-u_0-v) \frac{c}{2}+\zeta^*\alpha+a   \right) \left( \frac{\bar{d}}{2}(-|\zeta|^2-u_0+v)+\beta^*\zeta+\bar{b} \right)    \\
 + \left( (-|\zeta|^2-u_0-v)\frac{d}{2}+\zeta^*\beta+b \right) \left( \frac{\bar{c}}{2}(-|\zeta|^2-u_0+v)+\alpha^*\zeta+\bar{a}   \right)  \\
 +\left( (-|\zeta|^2-u_0-v) \frac{\delta}{2}^*+\zeta^*A+\gamma^* \right)\left( \frac{\delta}{2}(-|\zeta|^2-u_0+v)+A^*\zeta+\gamma \right) 
\end{eqnarray*}

If we further expand by distributing and collecting the terms in $(-|\zeta|^2-u_0)$, $\zeta$, $v$ (and their conjugates and products), all terms vanish by the relations $gg^{-1}=g^{-1}g={\rm Id}$ applied to the standard forms of $g,g^{-1}$ from equation~(\ref{standardform}),
except for the term in $(-|\zeta|^2-u_0)/2$ which has coefficient $(a\bar{d}+b\bar{c}+\gamma^* \delta)+(c\bar{b}+d\bar{a}+\delta^*\gamma)=2$
and the term in $\zeta\zeta^*$ which has coefficient $\alpha\beta^*+\beta\alpha^*+AA^*={\rm Id}$. 
Therefore the left-hand side simplifies to $-u_0$, whence:
 $$g^{-1}(z) \in H_{u_0} \iff -u_0=-u_0|\zeta_3|^2 \iff |\zeta_3|^2=1 \iff z \in S.$$ 
 This proves the claim and hence the Lemma. \EPf

\begin{cor}\label{cordepthradius} Let $E$ be a number field such that $\mathcal{O}_E$ is a principal ideal domain, $p \in \partial_\infty X$ an $E$-rational point with depth $n \geqslant 1$ and $g_p \in {\rm U}(H,\mathcal{O}_E)$ satisfying $g_p(\infty)=p$. Then $H_{u_0} \cap g_p(H_{u_0})=H_{u_0} \cap S$, with $S$  the extended Cygan sphere centered at $p$ with radius $\left(\frac{4}{n}\right)^{1/4}$. In particular: 
$H_{u_0} \cap g_p(H_{u_0})=\emptyset \iff u_0>u(n)=\frac{2}{\sqrt{n}}$.
\end{cor}

\Pf Since $g_p(\infty)=p$ and $e_1=(1,0,...,0)^T$ is a lift of $\infty$, the first column vector of $g_p$ is a lift $P_0$ of $p$, and since $g_p \in {\rm U}(H,\mathcal{O}_E)$ it is an integral lift. In fact by Lemma~\ref{primitivecolumn} it is a primitive lift, therefore the depth of $p$ is $|\langle P_0,e_1\rangle|^2=|g_{n+1,1}|^2$, denoting as above $g_p=(g_{i,j})$, and the result follows from Proposition~\ref{depthradius}. The second part of the statement follows by using this radius in the formula (\ref{extcygan}) for the extended Cygan metric. \EPf


  
We will also use the following observation, which is Lemma~1 of \cite{FFP}, in our covering arguments:

\begin{lem}\label{cyganconvex} Extended Cygan balls are affinely convex in horospherical coordinates.
\end{lem}  

Finally, when considering the action of a discrete subgroup $\Gamma_\infty$ of ${\rm Isom}(\partial_\infty X)$ (relative to the Cygan metric) it is convenient to consider its vertical and horizontal components defined as follows (see \cite{FP} for the case $\K=\C$ and $n=2$). The homomorphism $\Pi:{\rm Heis}(\K,n) \rightarrow \K^{n-1}$ given by projection to the first factor in the decomposition of the set ${\rm Heis}(\K,n)$ as $\K^{n-1} \times {\rm Im} \, \K$ induces a short exact sequence:

\begin{equation}
1 \longrightarrow {\rm Im} \, \K \longrightarrow {\rm Isom}({\rm Heis}(\K,n)) \xrightarrow{\Pi ^*} {\rm Isom}(\K^{n-1}) \longrightarrow 1 ,
\end{equation}  
where the isometries of $\K^{n-1}$ are relative to the Euclidean metric and and ${\rm Im} \, \K$ acts by ("vertical") translations. Denoting $\Gamma_\infty^v=\Gamma_\infty \cap {\rm Im} \, \K$ and $\Gamma_\infty^h=\Pi ^*(\Gamma_\infty)$ this gives the short exact sequence:

\begin{equation}\label{exact}
1 \longrightarrow \Gamma_\infty^v \longrightarrow \Gamma_\infty \xrightarrow{\Pi ^*} \Gamma_\infty^h \longrightarrow 1 .
\end{equation}

\section{Picard modular groups}

In this section we use the method described in Section~\ref{method} to compute presentations for the Picard modular groups $\Gamma(d)={\rm PU}(2,1,\mathcal{O}_d)$ with $d=1,3,7$. The following propositions summarize the results in this section. Recall that presentations for $\Gamma(d)={\rm PU}(2,1,\mathcal{O}_d)$ with $d=1,3$ were obtained in \cite{FP} and \cite{FFP} respectively. We only include these cases as test cases for our method; it turns out that the presentation we obtain when $d=1$ simplifies nicely (thanks to Magma \cite{Mag}) so we include it below, whereas the presentation we obtain when $d=3$ is much more complicated than the Falbel-Parker presentation, so we only include it in the appendix. We note the abelianization in each case as a corollary of the presentation.  

\begin{cor}{\bf ([FP])} The abelianization of $\Gamma(3)$ is $\Z/6\Z$.
\end{cor}

\begin{prop}\label{presentationd=1} The Picard modular group $\Gamma(1)={\rm PU}(2,1,\mathcal{O}_1)$ admits the presentation $\la I, A  \, | \, \mathcal{R}_1 \ra$, where $\mathcal{R}_1$ is the following set of 6 relations:
$$\begin{array}{l}
I^2 = {\rm Id} \\
    A^8 = {\rm Id} \\
    I  A^{-2}  I  A^2  I  A^2  I  A^{-2} = {\rm Id} \\
    (I  A^3  I  A^{-3})^3 = {\rm Id} \\
    (A^{-1}  I  A^{-2}  I  A^{-1}  I  A^3  I  A^{-1}  I 
    A^{-2}  I)^3 = {\rm Id} \\
    I  A^{-2}  I  A  I  A^3  I  A  I  A^{-2}  I 
    A^{-1}  I  A  I  A^2  I  A^{-1}  I  A  I  A^{-2}  -\\
    I  A^{-1}  I  A  I  A^{-2}  I  A^{-1}  I  A  I 
    A^{-2}  I  A^{-1}  I  A^3  I  A  I  A^3  I  A^3
     I  A^{-1}  I  A^3  I  A^3  I  A = {\rm Id} 
\end{array}$$
\end{prop}

\begin{cor}{\bf ([FFP])} The abelianization of $\Gamma(1)$ is $\Z/2\Z \times \Z/4\Z$.
\end{cor}

\begin{prop}\label{presentationd=7} The Picard modular group $\Gamma(7)={\rm PU}(2,1,\mathcal{O}_7)$ admits the presentation $\la T_1, R, I\, | \, \mathcal{R}_7 \ra$, where $\mathcal{R}_7$ is the following set of 13 relations:
$$\begin{array}{l}
R^2 = {\rm Id} \\
    I^2 = {\rm Id} \\
    (R  I)^2 = {\rm Id} \\
    R  T_1  R  T_1 = T_1  R  T_1 R  \\
    (T_1  I  T_1^{-1}  R)^4 = {\rm Id} \\
    (T_1^{-1}  I  T_1  R)^4 = {\rm Id} \\
    T_1^{-1}  I  T_1^{-1}  I  T_1  I  T_1  I  T_1^{-3}  I  T_1 
    I  T_1  I  T_1^{-1}  I  T_1^{-1} = {\rm Id} \\
    (T_1^{-1}  I  T_1  I  T_1  I  T_1^{-1}  I  T_1^{-1}  I)^2 =
    {\rm Id} \\
    (I  T_1^{-1}  R)^7 = {\rm Id} \\
    T_1^{-1}  I  T_1  I  T_1  I  T_1^{-2}  I  T_1^{-1}  I  T_1 
    I  T_1^2  I  T_1^{-1}  I  T_1^{-1}  I  T_1  I = {\rm Id} \\
    T_1^{-1}  I  T_1  I  T_1  I  R  T_1  I  R  T_1  I 
    T_1  I  T_1^{-1}  I  T_1^{-1}  I  T_1  R  T_1^{-1}  I  R 
    T_1^{-1}  I = {\rm Id} \\
    R  T_1  I  R  T_1  I  T_1  I  T_1^{-1}  I  T_1^{-1}  I
     R  T_1^{-1}  I  R  T_1^{-1}  I  T_1^{-1}  I  T_1  I  T_1 
    I  T_1^{-1} = {\rm Id} \\
    R  T_1  I  R  T_1  R  T_1^{-1}  I  T_1  I  T_1  I 
    R  T_1  I  T_1  I  T_1^{-1}  R  T_1  R  I  T_1  R 
    T_1^{-1}  I  T_1  I  T_1  I  T_1  I  T_1^{-1} = {\rm Id} 
 \end{array}$$
\end{prop}

\begin{cor} The abelianization of $\Gamma(7)$ is $\Z/2\Z$.
\end{cor}

\medskip 

The action of $\Gamma_\infty(d)={\rm Stab}_{\Gamma(d)}(\infty)$ on $\partial {\rm H}^2_\C$ is well understood for all $d$, see \cite{FP} for $d=3$, \cite{FFP} for $d=1$ and Section~5.3 of \cite{PW} for all other values (using unpublished notes of Falbel-Francsics-Parker). We will refer to these papers for presentations and fundamental domains for  $\Gamma_\infty(d)$ which we state in Lemmas~\ref{d=3Gammainfty}, \ref{d=1Gammainfty} and \ref{d=7Gammainfty}. We will denote $\tau=\frac{1+i\sqrt{d}}{2}$ when $d \equiv 3$ (mod 4), so that $\mathcal{O}_d=\Z[\tau]$.

\subsection{The Eisenstein-Picard modular group $\Gamma(3)={\rm PU}(2,1,\mathcal{O}_3)$}

{\bf Presentation and fundamental domain for the cusp stabilizer $\Gamma_\infty(3)$:}

\begin{lem}\label{d=3Gammainfty}
\begin{enumerate} 
\item The cusp stabilizer $\Gamma_\infty(3)$ admits the following presentation:

$$\Gamma_\infty(3)=\left\langle\ T_1, \, T_\tau, \, R  \ \Big\vert\ 
\begin{array}{c}  \left[ [ T_1,T_\tau ], T_1 \right], \,  \left[ [ T_1,T_\tau ], T_\tau \right], \,\left[ [ T_1,T_\tau ], R \right], \\
R^{-1}T_\tau R=T_1,  \, R^{-1}T_1 R=T_1T_\tau^{-1}, \, R^6
\end{array}\right\rangle.$$

\item Let $D_\infty(3) \subset \partial {\rm H}_\C^2$ be the affine convex hull of the points with horospherical coordinates $(0,0)$, $(1,0)$, $(\frac{\tau +1}{3},0)$, $(0,2\sqrt{3})$, $(1,2\sqrt{3})$, $(\frac{\tau +1}{3},2\sqrt{3})$. Then $D_\infty(3)$ is a fundamental domain for $\Gamma_\infty(3)$ acting on $\partial {\rm H}_\C^2 \setminus \{ \infty \}$.
\end{enumerate}
\end{lem}  
Concretely, we use the following generators for $\Gamma_\infty(3)$ (recall that $\tau=\frac{1+i\sqrt{3}}{2}$)
$$
\begin{array}{ccc} T_1=T_{(1,\sqrt{3})}=\left[ \begin{array}{ccc} 1 & -1 & \tau^2 \\ 0 & 1 & 1 \\ 0 & 0 & 1 \end{array}\right] & 
  T_\tau = T_{(\tau,\sqrt{3})}=\left[ \begin{array}{ccc} 1 & -\bar{\tau} & \tau^2 \\ 0 & 1 & \tau \\ 0 & 0 & 1 \end{array}\right] &
  R=\left[ \begin{array}{ccc} 1 & 0 & 0 \\ 0 & \tau & 0 \\ 0 & 0 & 1 \end{array}\right]
\end{array}$$

{\bf Covering depth and $\Q[i\sqrt{3}]$-rational points in $D_\infty(3)$:}

We denote $B\left((z,t),r\right)$ the open extended Cygan ball centered at $p=(z,t) \in \partial_\infty {\rm H}_\C^2$ with radius $r$ (see Equation~\ref{extcygan} for the definition of the extended Cygan metric). Recall that $u(n)=\frac{2}{\sqrt{n}}$ is the height at which balls of depth $n$ appear, in the sense of Corollary~\ref{cordepthradius}. 

\begin{center}
\begin{figure}\label{d=3pic}
\caption{Covering the prism $D_\infty(3)$ by Cygan balls of depth 1}
\includegraphics[height=8cm]{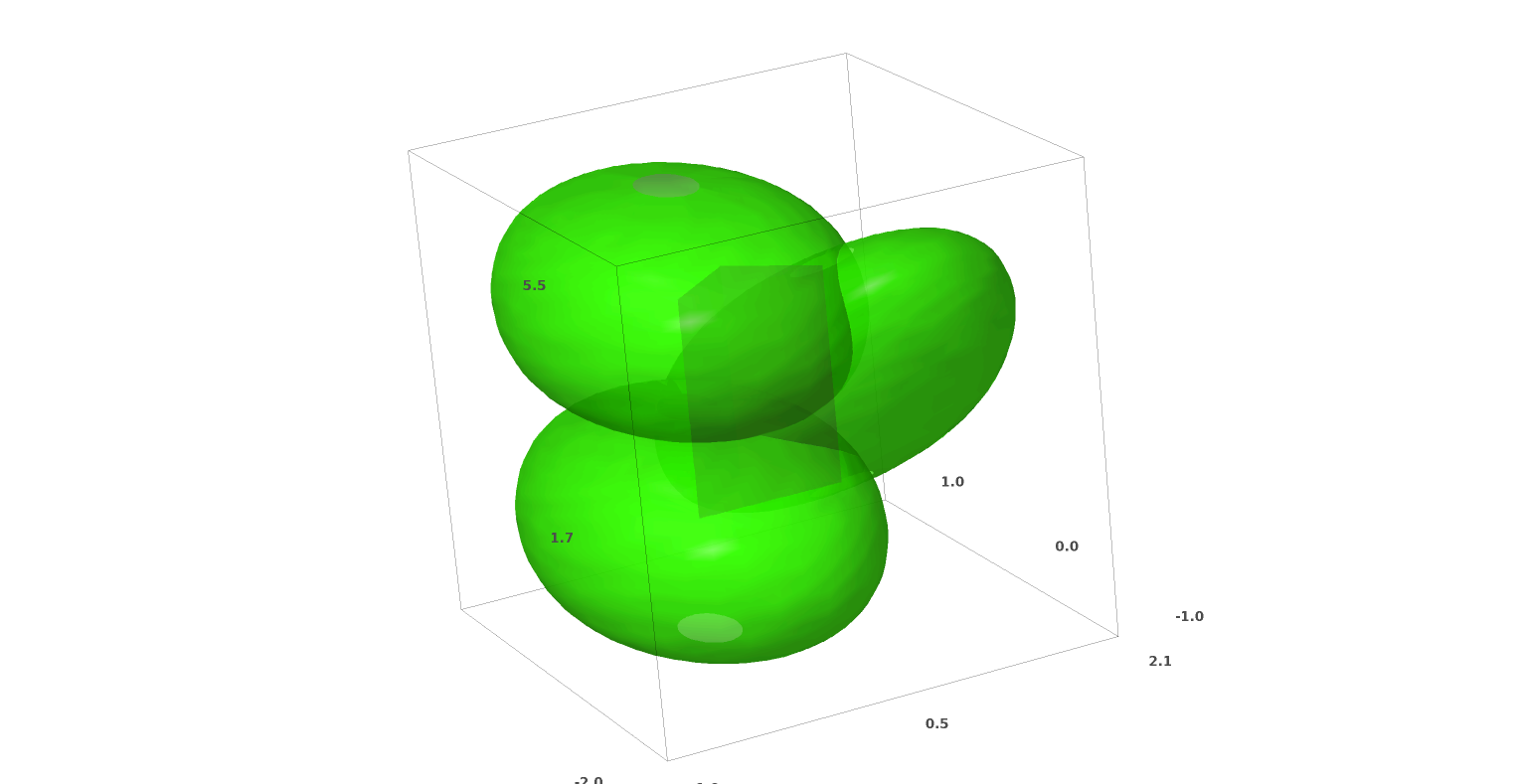}
\end{figure}
  \end{center}

\begin{lem}\label{d=3cover} Let $u=u(5)+\varepsilon=0.895$ and $H_u$ the horosphere of height $u$ based at $\infty$. Then the prism $D_\infty(3)\times \{ u \}$ is covered by the intersections with $H_u$ of the following extended Cygan balls of depth 1:  $B\left((0,0),\sqrt{2}\right)$, $B\left((0,2\sqrt{3}),\sqrt{2}\right)$ and $B\left((1,\sqrt{3}),\sqrt{2}\right)$.
\end{lem}
We omit the proof, which is similar to the proof of Lemma~\ref{d=7cover} but much simpler; see Figure~\ref{d=3pic}.

  Note that in order to cover $D_\infty(3)\times \{ u \}$ we only need balls of depth 1, in particular none of depths 3 or 4 even though they are present at the height $u=u(5)+\varepsilon$ which we consider. It seems however necessary to pass to this height, as we have observed experimentally that at height $u=u(4)+\varepsilon$ balls of depth at most 3 do not cover $D_\infty(3)\times \{ u \}$ (there are no $\Q[i\sqrt{3}]$-rational points of depth 2), i.e. that the covering depth is more than 3.
  
\begin{cor} The covering depth of $\Gamma(3)$ is at most 4.
\end{cor}  

By inspection, we see that the $\Q[i\sqrt{3}]$-rational points in $D_\infty(3)$ with depth at most 4 are, in horospherical coordinates:

\begin{itemize}
\item Depth 1: $(0,0)$, $(0,2\sqrt{3})$ and $(1,\sqrt{3})$, all in the same $\Gamma_\infty(3)$-orbit;
\item Depth 3: $(0,\frac{2}{3}\sqrt{3})$, $(1,\frac{5}{3}\sqrt{3})$ in one $\Gamma_\infty(3)$-orbit, andF $(0,\frac{4}{3}\sqrt{3})$, $(1,\frac{1}{3}\sqrt{3})$ in the other;
  \item Depth 4:   $(0,\sqrt{3})$, $(1,0)$ and $(1,2\sqrt{3})$, all in the same $\Gamma_\infty(3)$-orbit.
  \end{itemize}

Integral lifts of representatives of $\Gamma_\infty(3)$-orbits of these points are:

$$
\begin{array}{cccc} p_0=\left[ \begin{array}{c} 0 \\ 0 \\ 1 \end{array}\right] & 
  p_{3,1}=\left[ \begin{array}{c} -1 \\ 0 \\ i\sqrt{3} \end{array}\right] &
  
  p_{3,2}=\left[ \begin{array}{c} -2 \\ 0 \\ i\sqrt{3} \end{array}\right]  &
  p_{4}=\left[ \begin{array}{c} i\sqrt{3} \\ 0 \\ 2  \end{array}\right] 
  \end{array}
$$

{\bf Generators:} The following elements $A_\alpha \in \Gamma(3)$ map the point $\infty=[1,0,0]^T$ to the corresponding $p_\alpha$ as above (for $\alpha=0; 3,1; 3,2; 4$) :

$$
\begin{array}{cc} A_0=I_0=\left[ \begin{array}{ccc} 0 & 0 & 1 \\ 0 & -1 & 0 \\ 1 & 0 & 0 \end{array}\right] & 
  A_{3,1}=A_3=\left[ \begin{array}{ccc} -1 & 0 & i\sqrt{3} \\ 0 & 1 & 0 \\ i\sqrt{3} & 0 & 2 \end{array}\right] \\
  \\
  A_{3,2}=A_3^{-1} 
  &  A_{4}=\left[ \begin{array}{ccc} -i\sqrt{3} & 0 & -2\\ 0 & -1 & 0 \\ -2 & 0 & i\sqrt{3}  \end{array}\right] 
  \end{array}
$$



{\bf Relations:}  We obtain a complete set of relations between these generators by applying generators to points of depth at most 4 as described in part (5) of section~\ref{method}. The detailed steps of the computation can be found on the companion Sagemath Jupyter notebook in \cite{MCode}. The direct output is a presentation with 8 generators and 583 relations, which simplifies, thanks to Magma \cite{Mag}, to the presentation in the appendix. More specifically, we obtain that particular simplification by specifying a subset of the generators which must be preserved, using in this case the command {\tt Simplify (G: Preserve:=[1,4,5]);}.

To illustrate the steps involved we compute by hand the cycles and relations corresponding to triples $(a,b,c)$ for which the point $A_a(p_b)$ has depth at most 4, hence is of the form $\gamma_\infty (p_c)$ for some $\gamma_\infty \in \Gamma_\infty$ and $p_c \in D_\infty$ with depth at most 4. In the notation of part (5) of section~\ref{method}, this corresponds to $\gamma_\infty^1={\rm Id}$ and  $\gamma_\infty^3=\gamma_\infty$. The results are listed in Table~2, with the same notation as in the example of section~\ref{toyexample}.

\begin{table}[h]\label{d=3table}
  \begin{center}
  {\renewcommand{\arraystretch}{1.2}%
\begin{tabular}{|c|c|c|c|}
\hline 
$\mathbf{A.p}$ & $\mathbf{p'}$ & $\mathbf{A_{p'}^{-1}W^{-1}A A_p}$ & $\mathbf{W'}$  \\
\hline

$I_0 p_0$ & $\infty$ & $I_0^2$ & ${\rm Id}$\\
\hline

$A_4 p_4$ & $\infty$ & $A_4^2$ & ${\rm Id}$\\
\hline

$R p_0$ & $p_0$ & $I_0RI_0$ & $R$\\
\hline

$R p_{3,1}$ & $p_{3,1}$ & $A_3^{-1}RA_3$ & $R$\\
\hline

$R p_{3,2}$ & $p_{3,2}$ & $A_3RA_3^{-1}$& $R$\\
\hline

$I_0 (0,2\sqrt{3})$ & $p_{3,2}$ & $A_3T_vI_0T_vI_0$ & $R^3$\\
\hline

$I_0 (1,\sqrt{3})$ & $p_0$ & $I_0T_\tau ^{-1}T_vI_0T_1I_0$ & $T_v^{-1}T_1T_\tau ^{-1}$\\
\hline

$I_0 p_{3,1}$ & $p_0$ & $I_0T_vI_0A_3$ & $R^3T_v^{-1}$\\
\hline

$I_0 p_{3,2}$ & $p_4$ & $A_4^{-1}T_vI_0A_3^{-1}$ & $R^3$\\
\hline

$I_0 p_4$ & $p_{3,1}$ & $A_3^{-1}T_vI_0A_4$ & $R^3$\\
\hline

$A_3  p_0$ & $p_4$ & $A_4^{-1}A_3I_0$ & $R^3T_v$ \\
\hline

$A_3  p_4$ & $p_0$ & $I_0T_v^{-1}A_3A_4$ &  $R^3 $ \\
\hline

$A_3  p_{3,1}$ & $p_{3,2}$ & $A_3^3$ &  $R^3 $ \\
\hline

$A_4  p_{0}$ & $p_{3,2}$ & $A_3A_4I_0$ &  $R^3T_v$ \\
\hline

$A_4  p_{3,1}$ & $p_0$ & $I_0T_v^{-1}A_4A_3$ &  $R^3$ \\
\hline

\end{tabular}
}
\end{center}
\caption{Action of generators on vertices for $d=3$}
\end{table}

\subsection{The Gauss-Picard modular group $\Gamma(1)={\rm PU}(2,1,\mathcal{O}_1)$}

{\bf Presentation and fundamental domain for the cusp stabilizer $\Gamma_\infty(1)$:}

\begin{lem}\label{d=1Gammainfty} 
\begin{enumerate}
\item The cusp stabilizer $\Gamma_\infty(1)$ admits the following presentation:

$$\Gamma_\infty(1)= \left\langle\ T_2, \, T_\tau, \, T_v, \, R  \ \Big\vert\ 
\begin{array}{c} [ T_\tau,T_2 ]=T_v^4, \, [ T_v, T_2 ], \, [ T_v,T_\tau ], \, [ T_v,R ], \, R^4,  \\

  R T_2 R^{-1}=T_\tau^2T_2^{-1}T_v^4, \, R T_\tau R^{-1}=T_\tau T_2^{-1}T_v^2
  \end{array}\right\rangle.$$

\item Let $D_\infty(1) \subset \partial {\rm H}_\C^2$ be the affine convex hull of the points with horospherical coordinates $(0,0)$, $(1,0)$, $(\tau,0)$, $(0,2)$, $(1,2)$, $(\tau,2)$. Then $D_\infty(1)$ is a fundamental domain for $\Gamma_\infty(1)$ acting on $\partial {\rm H}_\C^2 \setminus \{ \infty \}$.
\end{enumerate}
\end{lem}  
Concretely, we use the following generators for $\Gamma_\infty(1)$ (denoting $\tau=1+i$):

$$
\begin{array}{cc} T_2=T_{(2,0)}=\left[ \begin{array}{ccc} 1 & -2 & -2 \\ 0 & 1 & 2 \\ 0 & 0 & 1 \end{array}\right] & 
  T_\tau = T_{(\tau,0)}=\left[ \begin{array}{ccc} 1 & -\bar{\tau} & -1 \\ 0 & 1 & \tau \\ 0 & 0 & 1 \end{array}\right] \\
  \\
  T_v=T_{(0,2)}=\left[ \begin{array}{ccc} 1 & 0 & i \\ 0 & 1 & 0 \\ 0 & 0 & 1 \end{array}\right] & 
  R=\left[ \begin{array}{ccc} 1 & 0 & 0 \\ 0 & i & 0 \\ 0 & 0 & 1 \end{array}\right]
\end{array}$$

{\bf Covering depth and $\Q[i]$-rational points in $D_\infty(1)$:}

We denote $B\left((z,t),r\right)$ the open extended Cygan ball centered at $p=(z,t) \in \partial_\infty {\rm H}_\C^2$ with radius $r$ (see Equation~\ref{extcygan} for the definition of the extended Cygan metric). Recall that $u(n)=\frac{2}{\sqrt{n}}$ is the height at which balls of depth $n$ appear, in the sense of Corollary~\ref{cordepthradius}. 

\begin{lem}\label{d=1cover} Let $u=u(5)+\varepsilon=0.895$ and $H_u$ the horosphere of height $u$ based at $\infty$. Then the prism $D_\infty(1)\times \{ u \}$ is covered by the intersections with $H_u$ of the following extended Cygan balls:
  \begin{itemize}
\item (depth 1)  $B\left((0,0),\sqrt{2}\right)$, $B\left((0,2),\sqrt{2}\right)$, $B\left((\tau,0),\sqrt{2}\right)$, $B\left((\tau,2),\sqrt{2}\right)$, $B\left((1,\sqrt{7}),\sqrt{2}\right)$,
\item (depth 2) $B\left((1,1),\sqrt[4]{2}\right)$.

  \end{itemize}
\end{lem}

\begin{center}
\begin{figure}\label{d=1pic}
\caption{Covering the prism $D_\infty (1)$ by Cygan balls of depth 1 and 2}
\includegraphics[height=8cm]{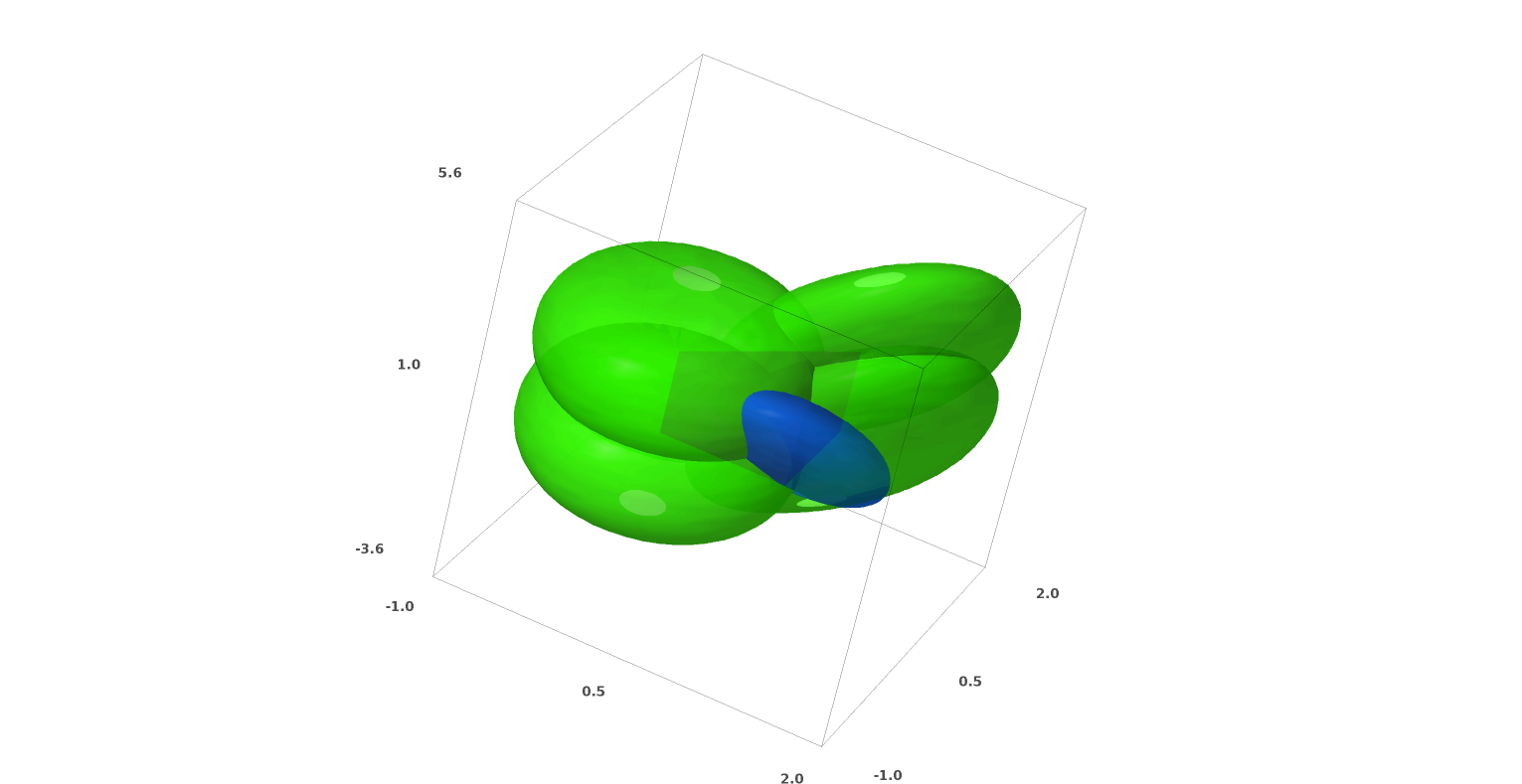}
\end{figure}
  \end{center}

  We omit the proof, which is similar to the proof of Lemma~\ref{d=7cover} but simpler; see Figure~\ref{d=1pic}.

  Note that in order to cover  $D_\infty(1)\times \{ u \}$ we only need balls of depth at most 2, in particular none of depth 4 even though they are present at the height $u=u(5)+\varepsilon$ which we consider. It is however necessary to pass to this height, as we have observed experimentally that at height $u=u(4)+\varepsilon$ balls of depth at most 2 do not cover $D_\infty(1)\times \{ u \}$ (there are no $\Q[i]$-rational points of depth 3).

\begin{cor} The covering depth of $\Gamma(1)$ is at most 4.
\end{cor}

By inspection, we see that the $\Q[i]$-rational points in $D_\infty(1)$ with depth at most 4 are, in horospherical coordinates:

\begin{itemize}
\item Depth 1: $(0,0)$, $(0,2)$, $(\tau,0)$ and $(\tau,2)$, all in the same $\Gamma_\infty(1)$-orbit;
\item Depth 2: $(1,1)$
  \item Depth 4:   $(0,1)$, $(\tau,1)$ in one $\Gamma_\infty(1)$-orbit and $(1,0)$, $(1,2)$ in the other.
  \end{itemize}

Integral lifts of representatives of $\Gamma_\infty(1)$-orbits of these points are:

$$
\begin{array}{cccc} p_0=\left[ \begin{array}{c} 0 \\ 0 \\ 1 \end{array}\right] & 
  p_{2}=\left[ \begin{array}{c} -1 \\ 1+i \\ 1+i \end{array}\right] &
  p_{4,1}=\left[ \begin{array}{c} i \\ 0 \\ 2 \end{array}\right]  &
  p_{4,2}=\left[ \begin{array}{c} -1 \\ 2 \\ 2  \end{array}\right] 
  \end{array}
$$

{\bf Generators:}  The following elements $A_\alpha \in \Gamma(1)$ map the point $\infty=[1,0,0]^T$ to the corresponding $p_\alpha$ as above (for $\alpha=0; 2; 4,1; 4,2$) :

$$
\begin{array}{cc} A_0=I_0=\left[ \begin{array}{ccc} 0 & 0 & 1 \\ 0 & -1 & 0 \\ 1 & 0 & 0 \end{array}\right] & 
  A_{2}=\left[ \begin{array}{ccc} -1 & -2 & 2+i \\ 1+i & 2+i & -1-i \\ 1+i & 1+i & -i \end{array}\right] \\
  \\
  A_{4,1}=\left[ \begin{array}{ccc} i & 0 & 1 \\ 0 & -1 & 0 \\ 2 & 0 & -i \end{array}\right] & 
    A_{4,2}=\left[ \begin{array}{ccc} -1 & -2 & 2\\ 2 & 3 & -2 \\ 2 & 2 & -1  \end{array}\right] 
  \end{array}
$$

{\bf Relations:}  We obtain a complete set of relations between these generators by applying generators to points of depth at most 4 as described in part (5) of section~\ref{method}. The detailed steps of the computation can be found on the companion Sagemath Jupyter notebook in \cite{MCode}. The direct output is a presentation with 8 generators and 247 relations, which simplifies, thanks to Magma \cite{Mag}, to the presentation in Proposition~\ref{presentationd=1}. More specifically, we obtain that particular simplification by specifying a subset of the generators which must be preserved, using in this case the command {\tt Simplify (G: Preserve:=[5,6]);}.

To illustrate the steps involved we compute by hand the cycles and relations corresponding to triples $(a,b,c)$ for which the point $A_a(p_b)$ has depth at most 4, hence is of the form $\gamma_\infty (p_c)$ for some $\gamma_\infty \in \Gamma_\infty$ and $p_c \in D_\infty$ with depth at most 4. In the notation of part (5) of section~\ref{method}, this corresponds to $\gamma_\infty^1={\rm Id}$ and  $\gamma_\infty^3=\gamma_\infty$. The results are listed in Table~3, with the same notation as in the example of section~\ref{toyexample}.

\begin{table}[h]\label{d=1table}
  \begin{center}
  {\renewcommand{\arraystretch}{1.2}%
\begin{tabular}{|c|c|c|c|}
\hline 
$\mathbf{A.p}$ & $\mathbf{p'}$ & $\mathbf{A_{p'}^{-1}W^{-1}A A_p}$ & $\mathbf{W'}$  \\
\hline

$I_0 p_0$ & $\infty$ & $I_0^2$ & ${\rm Id}$\\
\hline

$A_{4,1} p_{4,1}$ & $\infty$ & $A_{4,1}^2$ & ${\rm Id}$\\
\hline

$A_{4,2} p_{4,2}$ & $\infty$ & $A_{4,2}^2$ & ${\rm Id}$\\
\hline

$R p_0$ & $p_0$ & $I_0RI_0$ & $R$\\
\hline

$R p_{4,1}$ & $p_{4,1}$ & $A_{4,1}^{-1}RA_{4,1}$ & $R$\\
\hline

$I_0 (0,2)$ & $p_0$ & $I_0T_vI_0T_vI_0$& $RT_v^{-1}$\\
\hline

$I_0 (\tau,0)$ & $(\tau,0)$ & $(T_\tau I_0)^{-1}I_0T_\tau I_0$ & $R^3T_\tau$\\
\hline

$I_0 (\tau,2)$ & $p_2$ & $(T_vT_\tau I_0)^{-1}I_0T_\tau T_2^{-1}A_2$ & $R^2$\\
\hline

$I_0 p_2$ & $p_0$ & $I_0T_\tau ^{-1}T_vI_0A_2$ & $RT_2^{-1}T_v^{-1}$\\
\hline

$I_0 p_{4,1}$ & $p_0$ & $I_0T_v ^2I_0A_{4,1}$ & $RT_v^{-1}$\\
\hline

$I_0 p_{4,2}$ & $p_0$ & $I_0T_2^{-1}I_0A_{4,2}$ & $T_2^{-1}$\\

\hline
$A_2 p_0$ & $p_0$ & $I_0T_v ^2I_0A_{4,1}$ & $RT_v^{-1} $\\

\hline
$A_2 p_2$ & $\infty$ & $A_2^2$ & $R^{-1}T_2^{-1}T_\tau T_v^{-3}$\\

\hline
$A_2 p_{4,1}$ & $p_2$  & $A_2^{-1}T_v^{-4}T_2^{-1}T_\tau A_2A_{4,1}$ & $R^2T_2^{-1}$\\

\hline
$A_2 p_{4,2}$ & $p_2$ & $A_2^{-1}T_v ^{-1}A_2A_{4,2}$ & $T_2T_\tau ^{-1}R $\\

\hline
$A_{4,1} p_0$ & $p_0$ & $I_0T_v^{-1}A_{4,1}I_0$ & $R^{-1}T_v^2 $\\

\hline
$A_{4,1} p_2$ & $p_2$ & $A_2^{-1}A_{4,1}A_2$ & $R^2T_\tau ^{-1}T_2^{-1}T_v^{-2}$\\

\hline
$A_{4,2} p_0$ & $p_0$ & $I_0T_2^{-1}A_{4,2}I_0$ & $T_2^{-1} $\\

\hline
$A_{4,2} p_2$ & $(\tau,0)$ & $(T\tau I_0)^{-1}T_v A_{4,2}A_2$ & $R^{-1}T_\tau T_2^{-1}$\\

\hline
\end{tabular}
}
\end{center}
\caption{Action of generators on vertices for $d=1$}
\end{table}

\subsection{The Picard modular group $\Gamma(7)={\rm PU}(2,1,\mathcal{O}_7)$}\label{Pic_mod_Gamma(7)}

{\bf Presentation and fundamental domain for the cusp stabilizer $\Gamma_\infty(7)$:}

\begin{lem}\label{d=7Gammainfty}
\begin{enumerate}
\item The cusp stabilizer $\Gamma_\infty(7)$ admits the following presentation:

$$\Gamma_\infty(7)= \left\langle\ T_1, \, T_\tau, \, T_v, \, R \ \Big\vert\ 
  \begin{array}{c} [ T_\tau,T_1 ]=T_v, \, [ T_v,T_1 ], \, [ T_v, T_\tau ], \, [ T_v, R ], \, (RT_\tau)^2, \, (RT_1)^2=T_v, \,  R^2
\end{array}\right\rangle$$

\item Let $D_\infty(7) \subset \partial {\rm H}_\C^2$ be the affine convex hull of the points with horospherical coordinates $(0,0)$, $(1,0)$, $(\tau,0)$, $(0,2\sqrt{7})$, $(1,2\sqrt{7})$, $(\tau,2\sqrt{7})$. Then $D_\infty(7)$ is a fundamental domain for $\Gamma_\infty(7)$ acting on $\partial {\rm H}_\C^2 \setminus \{ \infty \}$.
\end{enumerate}
\end{lem} 
Concretely, we use the following generators for $\Gamma_\infty(7)$ (denoting $\tau=\frac{1+i\sqrt{7}}{2}$):

$$
\begin{array}{cc} T_1=T_{(1,\sqrt{7})}=\left[ \begin{array}{ccc} 1 & -1 & \tau-1 \\ 0 & 1 & 1 \\ 0 & 0 & 1 \end{array}\right] & 
  T_\tau = T_{(\tau,0)}=\left[ \begin{array}{ccc} 1 & -\bar{\tau} & -1 \\ 0 & 1 & \tau \\ 0 & 0 & 1 \end{array}\right] \\
  \\
  T_v=T_{(0,2\sqrt{7})}=\left[ \begin{array}{ccc} 1 & 0 & i\sqrt{7} \\ 0 & 1 & 0 \\ 0 & 0 & 1 \end{array}\right] & 
  R=\left[ \begin{array}{ccc} 1 & 0 & 0 \\ 0 & -1 & 0 \\ 0 & 0 & 1 \end{array}\right]
\end{array}$$

{\bf Covering depth and $\Q[i\sqrt{7}]$-rational points in $D_\infty(7)$:} \\

We denote $B\left((z,t),r\right)$ the open extended Cygan ball centered at $p=(z,t) \in \partial_\infty {\rm H}_\C^2$ with radius $r$ (see Equation~\ref{extcygan} for the definition of the extended Cygan metric). Recall that $u(n)=\frac{2}{\sqrt{n}}$ is the height at which balls of depth $n$ appear, in the sense of Corollary~\ref{cordepthradius}. 

\begin{lem}\label{d=7cover} Let $u=u(8)+\varepsilon=0.70711$ and $H_u$ the horosphere of height $u$ based at $\infty$. Then the prism $D_\infty(7)\times \{ u \}$ is covered by the intersections with $H_u$ of the following extended Cygan balls:
  \begin{itemize}
\item (depth 1)  $B\left((0,0),\sqrt{2}\right)$, $B\left((0,2\sqrt{7}),\sqrt{2}\right)$, $B\left((\tau,0),\sqrt{2}\right)$, $B\left((\tau,2\sqrt{7}),\sqrt{2}\right)$, $B\left((1,\sqrt{7}),\sqrt{2}\right)$,
\item (depth 2) $B\left((\tau/2,3\sqrt{7}/2),\sqrt[4]{2}\right)$, $B\left(((\tau+1)/2,\sqrt{7}),\sqrt[4]{2}\right)$
\item (depth 4) $B\left((\tau,\sqrt{7}),1\right)$.
  \end{itemize}
\end{lem}

\begin{center}
\begin{figure}\label{d=7pic}
\caption{Covering the prism $D_\infty (7)$ by Cygan balls of depth 1, 2 and 4}
\includegraphics[height=8cm]{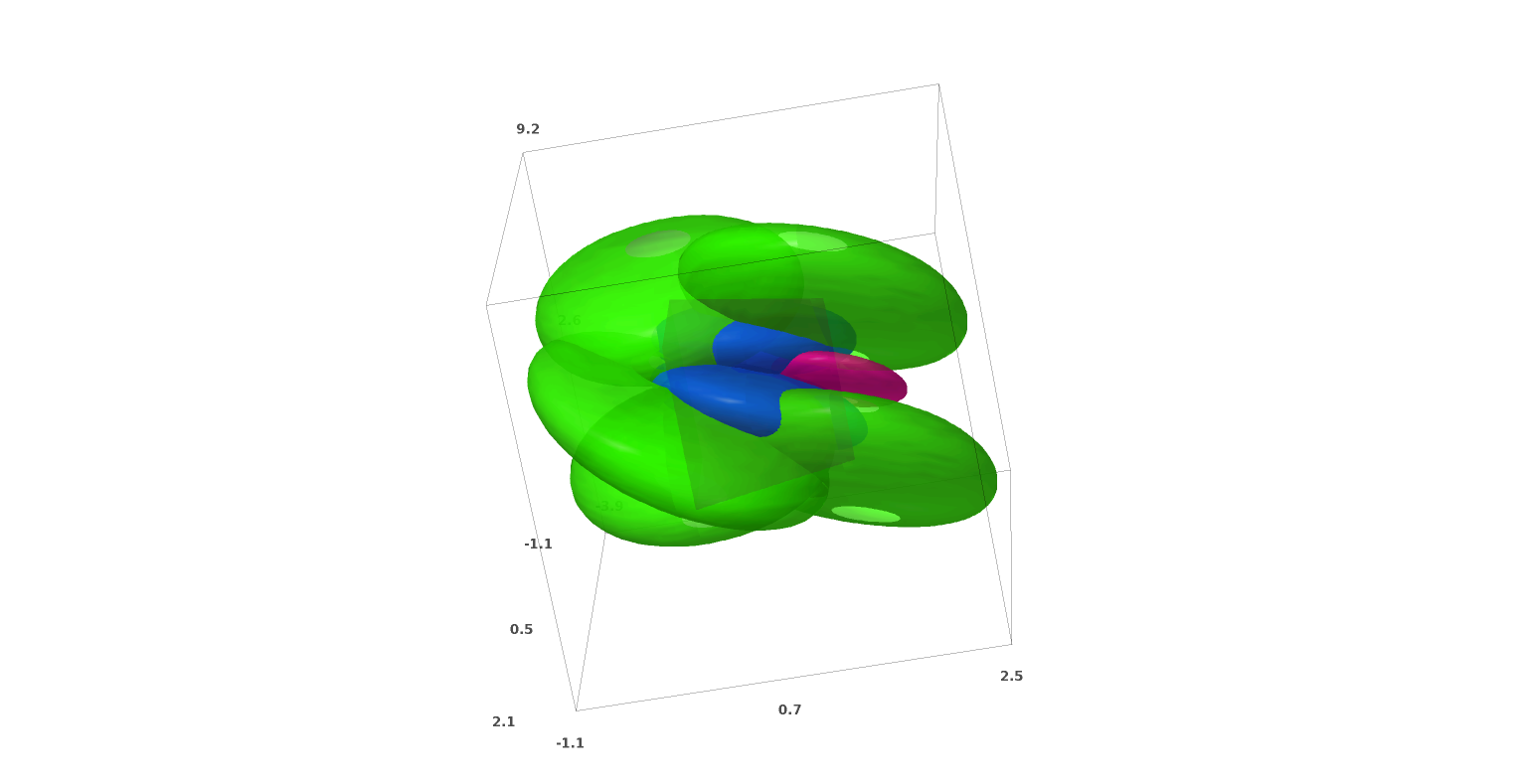}
\end{figure}
  \end{center}
 
\Pf Figure~\ref{d=7pic} shows the prism $D_\infty(7)$ and the relevant Cygan balls. We prove the result by dissecting the prism $D_\infty(7)\times \{ u \}$ into affine polyhedra, each of which lies in one of the extended Cygan balls. This is reminiscent of the proof of Proposition~5.2 of \cite{Zh}. Consider the following points of $\partial_\infty {\rm H}_\C^2$, in horospherical coordinates (see Figure~\ref{d=7coverpic}):

$$\begin{array}{ccccc}
  q_1=(0.65\tau,2.8) & q_2=(\tau,2) & q_3=(\tau,3.3) & q_4=(0.35+0.65\tau,3.4) & q_5=(0.7+0.3\tau,4.2) \\
  q_6=(1,4.3) & q_7=((\tau +1)/2,1.5) & q_8=(\tau,1.5) & q_9=((\tau +1)/2,0) & q_{10}=(1,1) \\
  q_{11}=(\tau,4) & q_{12}=((\tau +1)/2,2\sqrt{7}) & q_{13}=(\tau/2,4) & q_{14}=(\tau/2,2\sqrt{7}) & q_{15}=(0,3.5) \\
  q_{16}=(0.3\tau,3) & q_{17}=(\tau/2,1) & q_{18}=(0,1.7) & q_{19}=(\tau/2,0) &
\end{array}
$$

Denoting ${\rm Hull}(S)$ the affine hull (in horospherical coordinates) of a subset $S \subset H_u \simeq \partial_\infty {\rm H}_\C^2 \times \{ u \}$, we claim that the following affinely convex pieces of $D_\infty(7)\times \{ u \}$ are each contained in the corresponding (open) extended Cygan sphere:

\begin{itemize}
\item $D_1={\rm Hull}\left((0,0),(1,0),q_9,q_{10},q_{17},q_{18},q_{19}\right) \subset B\left((0,0),\sqrt{2}\right)$

  \item $D_2={\rm Hull}\left((0,2\sqrt{7}),(1,2\sqrt{7}),q_6,q_{12},q_{13},q_{14},q_{15} \right) \subset B\left((0,2\sqrt{7}),\sqrt{2}\right)$

  \item  $D_3={\rm Hull}\left( (\tau,0),q_7,q_8,q_9,q_{17},q_{19} \right) \subset B\left((\tau,0),\sqrt{2}\right)$ 

  \item $D_4={\rm Hull}\left((\tau,2\sqrt{7}),q_{11},q_{12},q_{13},q_{14} \right) \subset B\left((\tau,2\sqrt{7}),\sqrt{2}\right)$

  \item $D_5={\rm Hull}\left((1,\sqrt{7}),q_6,q_7,q_{10},q_{15},q_{16},q_{17},q_{18}   \right) \subset B\left((1,\sqrt{7}),\sqrt{2}\right)$

    \item $D_6={\rm Hull}\left(q_1,q_3,q_4,q_5,q_{11},q_{12},q_{13},q_{15},q_{16}  \right) \subset B\left((\tau/2,3\sqrt{7}/2),\sqrt[4]{2}\right)$

    \item $D_7={\rm Hull}\left(q_1,q_2,q_4,q_5,q_6,q_7,q_8,q_{16},q_{17} \right) \subset B\left(((\tau+1)/2,\sqrt{7}),\sqrt[4]{2}\right)$

      \item $D_8={\rm Hull}\left(q_1,q_2,q_3,q_4 \right) \subset B\left((\tau,\sqrt{7}),1\right)$
\end{itemize}  

To verify each of these claims, we check numerically that each of the vertices indeed belongs to the ball in question using Equation~(\ref{extcygan}), then extend to the whole affine covex hull by Lemma~\ref{cyganconvex}. For example, the point $q_1=(0.65\tau,2.8)$ indeed belongs to $B\left((\tau/2,3\sqrt{7}/2),\sqrt[4]{2}\right), B\left(((\tau+1)/2,\sqrt{7}),\sqrt[4]{2}\right)$ and $B\left((\tau,\sqrt{7}),1\right)$ because:

$$\begin{array}{c}
  d_{XC}\left((0.65\tau,2.8,u), (\tau/2,3\sqrt{7}/2,0) \right)\simeq 1.179 < \sqrt[4]{2}\simeq 1.189\\
  
  d_{XC}\left((0.65\tau,2.8,u), ((\tau+1)/2,\sqrt{7},0) \right)\simeq 1.172 < \sqrt[4]{2}\simeq 1.189\\

  d_{XC}\left((0.65\tau,2.8,u), (\tau,\sqrt{7},0) \right)\simeq 0.982 < 1
  
\end{array}$$
The result then follows as the prism $D_\infty(7)\times \{ u \}$ is the union of the affinely convex pieces $D_1,...,D_8$. Indeed, each of these pieces has all of its vertices on the boundary of the prism, the union of the boundaries of the pieces covers the boundary of the prism, and each piece shares a codimension-1 face with its neighbors, see Figure~\ref{d=7coverpic}. \EPf

\begin{cor} The covering depth of $\Gamma(7)$ is at most 7.
\end{cor}  
Note that in the above covering argument we have only needed balls of depth at most 4 (in particular none of depth 7) even though they are present at the height $u=u_8+\varepsilon$ which we consider. It is however necessary to pass to this height, as we have observed experimentally that at height $u=u_7+\varepsilon$ balls of depth at most 4 do not cover $D_\infty(7)\times \{ u \}$ (there are no $\Q[i\sqrt{7}]$-rational points of depths 5 or 6).

\begin{center}
\begin{figure}\label{d=7coverpic}

\begin{tikzpicture}[scale = 2]
\node[circle,fill=green,inner sep=0pt,minimum size=5pt,label=below left:{$(0,0)$}] at (0,0) (01) {};
\node[circle,fill=green,inner sep=0pt,minimum size=5pt] at (-.5,0) (02) {};
\node[circle,fill=green,inner sep=0pt,minimum size=5pt] at (0,-.5) (03) {};

\node[circle,fill=green,inner sep=0pt,minimum size=5pt,label=below right:{$(\tau,0)$}] at (2,0) (tau1) {};
\node[circle,fill=green,inner sep=0pt,minimum size=5pt] at (2.5,0) (tau4) {};
\node[circle,fill=green,inner sep=0pt,minimum size=5pt] at (2,-.5) (tau3) {};

\node[circle,fill=green,inner sep=0pt,minimum size=5pt,label=left:{$(1,0)$}] at (-1.5,0) (12) {};
\node[circle,fill=green,inner sep=0pt,minimum size=5pt,label=right:{$(1,0)$}] at (4.5,0) (14) {};
\node[circle,fill=green,inner sep=0pt,minimum size=5pt,label=below:{$(1,0)$}] at (.25,-1.468) (13) {};

\node[circle,fill=green,inner sep=0pt,minimum size=5pt,label=above left:{$(0,2\sqrt{7})$}] at (0,5.292) (v1) {};
\node[circle,fill=green,inner sep=0pt,minimum size=5pt] at (-.5,5.292) (v2) {};
\node[circle,fill=green,inner sep=0pt,minimum size=5pt] at (0,5.792) (v5) {};

\node[circle,fill=green,inner sep=0pt,minimum size=5pt,label=above right:{$(\tau,2\sqrt{7})$}] at (2,5.292) (taupv1) {};
\node[circle,fill=green,inner sep=0pt,minimum size=5pt] at (2.5,5.292) (taupv4) {};
\node[circle,fill=green,inner sep=0pt,minimum size=5pt] at (2,5.792) (taupv5) {};

\node[circle,fill=green,inner sep=0pt,minimum size=5pt,label=above:{$(1,2\sqrt{7})$}] at (-1.5,5.292) (1pv2) {};
\node[circle,fill=green,inner sep=0pt,minimum size=5pt,label=above:{$(1,2\sqrt{7})$}] at (4.5,5.292) (1pv4) {};
\node[circle,fill=green,inner sep=0pt,minimum size=5pt,label=above:{$(1,2\sqrt{7})$}] at (.25,6.76) (1pv5) {};

\draw (01) -- (tau1) -- (taupv1) -- (v1) -- (01);
\draw (02) -- (v2) -- (1pv2) -- (12) -- (02);
\draw (03) -- (tau3) -- (13) -- (03);
\draw (tau4) -- (taupv4) -- (1pv4) -- (14) -- (tau4);
\draw (taupv5) -- (1pv5) -- (v5) -- (taupv5);



\node[circle,fill=black,inner sep=0pt,minimum size=5pt,label=above left:{$q_1$}] at (1.3,2.8) (q1) {};

\node[circle,fill=black,inner sep=0pt,minimum size=5pt,label=right:{$q_2$}] at (2,2) (q21) {};
\node[circle,fill=black,inner sep=0pt,minimum size=5pt] at (2.5,2) (q24) {};

\node[circle,fill=black,inner sep=0pt,minimum size=5pt,label=right:{$q_3$}] at (2,3.3) (q31) {};
\node[circle,fill=black,inner sep=0pt,minimum size=5pt] at (2.5,3.3) (q34) {};

\node[circle,fill=black,inner sep=0pt,minimum size=5pt,label=right:{$q_4$}] at (2.85,3.4) (q4) {};

\node[circle,fill=black,inner sep=0pt,minimum size=5pt,label=below right:{$q_5$}] at (3.9,4.2) (q5) {};

\node[circle,fill=black,inner sep=0pt,minimum size=5pt,label=right:{$q_6$}] at (4.5,4.3) (q64) {};
\node[circle,fill=black,inner sep=0pt,minimum size=5pt,label=left:{$q_6$}] at (-1.5,4.3) (q62) {};

\node[circle,fill=black,inner sep=0pt,minimum size=5pt,label=right:{$q_7$}] at (3.5,1.5) (q7) {};

\node[circle,fill=black,inner sep=0pt,minimum size=5pt,label=right:{$q_8$}] at (2,1.5) (q81) {};
\node[circle,fill=black,inner sep=0pt,minimum size=5pt] at (2.5,1.5) (q84){};

\node[circle,fill=black,inner sep=0pt,minimum size=5pt,label=below:{$q_9$}] at (3.5,0) (q94) {};
\node[circle,fill=black,inner sep=0pt,minimum size=5pt,label=below right:{$q_9$}] at (1.125,-.984) (q93) {};

\node[circle,fill=black,inner sep=0pt,minimum size=5pt,label=right:{$q_{10}$}] at (4.5,1) (q104) {};
\node[circle,fill=black,inner sep=0pt,minimum size=5pt,label=left:{$q_{10}$}] at (-1.5,1) (q102) {};

\node[circle,fill=black,inner sep=0pt,minimum size=5pt,label=right:{$q_{11}$}] at (2,4) (q111) {};
\node[circle,fill=black,inner sep=0pt,minimum size=5pt] at (2.5,4) (q114) {};

\node[circle,fill=black,inner sep=0pt,minimum size=5pt,label=above:{$q_{12}$}] at (3.5,5.292) (q124){};
\node[circle,fill=black,inner sep=0pt,minimum size=5pt,label=above right:{$q_{12}$}] at (1.125,6.276) (q125){};

\node[circle,fill=black,inner sep=0pt,minimum size=5pt,label=above left:{$q_{13}$}] at (1,4) (q13) {};

\node[circle,fill=black,inner sep=0pt,minimum size=5pt,label=above:{$q_{14}$}] at (1,5.292) (q141) {};
\node[circle,fill=black,inner sep=0pt,minimum size=5pt] at (1,5.792) (q145) {};

\node[circle,fill=black,inner sep=0pt,minimum size=5pt,label=left:{$q_{15}$}] at (0,3.5) (q151) {};
\node[circle,fill=black,inner sep=0pt,minimum size=5pt] at (-.5,3.5) (q152) {};

\node[circle,fill=black,inner sep=0pt,minimum size=5pt,label=below left:{$q_{16}$}] at (.6,3) (q16) {};

\node[circle,fill=black,inner sep=0pt,minimum size=5pt,label=below right:{$q_{17}$}] at (1,1) (q17) {};

\node[circle,fill=black,inner sep=0pt,minimum size=5pt,label=left:{$q_{18}$}] at (0,1.7) (q181) {};
\node[circle,fill=black,inner sep=0pt,minimum size=5pt] at (-.5,1.7) (q182) {};

\node[circle,fill=black,inner sep=0pt,minimum size=5pt,label=below:{$q_{19}$}] at (1,0) (q191) {};
\node[circle,fill=black,inner sep=0pt,minimum size=5pt] at (1,-.5) (q193) {};

\draw (q111) -- (q13) -- (q151) -- (q16) -- (q1) -- (q31);
\draw (q141) -- (q13);
\draw (q1) -- (q21);
\draw (q16) -- (q17) -- (q81);
\draw (q181) -- (q17);
\draw (q191) -- (q17);

\draw (q62) -- (q152);
\draw (q102) -- (q182);

\draw (q193) -- (q93);

\draw (q104) -- (q94) -- (q7) -- (q64) -- (q5) -- (q124) -- (q114);
\draw (q84) -- (q7);
\draw (q24) -- (q4) -- (q34);
\draw (q4) -- (q5);

\draw (q145) -- (q125);

\node at (-1,.5) {1};
\node at (.5,.5) {1};
\node at (.5,-.9) {1};
\node at (4.2,.3) {1};
\node at (-1,4.5) {2};
\node at (.5,4.5) {2};
\node at (.5,6.2) {2};
\node at (4,4.7) {2};
\node at (1.5,.5) {3};
\node at (3,.5) {3};
\node at (1.5,-.65) {3};
\node at (1.5,4.7) {4};
\node at (1.5,5.9) {4};
\node at (2.8,4.8) {4};
\node at (-1,2.6) {5};
\node at (.4,2) {5};
\node at (4,1.7) {5};
\node at (1,3.5) {6};
\node at (3,4) {6};
\node at (1.4,2) {7};
\node at (3.4,2.7) {7};
\node at (1.8,2.64) {8};
\node at (2.6,2.9) {8};

\end{tikzpicture}
\caption{Affine cell decomposition of the prism $D_\infty(7)$ }
\end{figure}
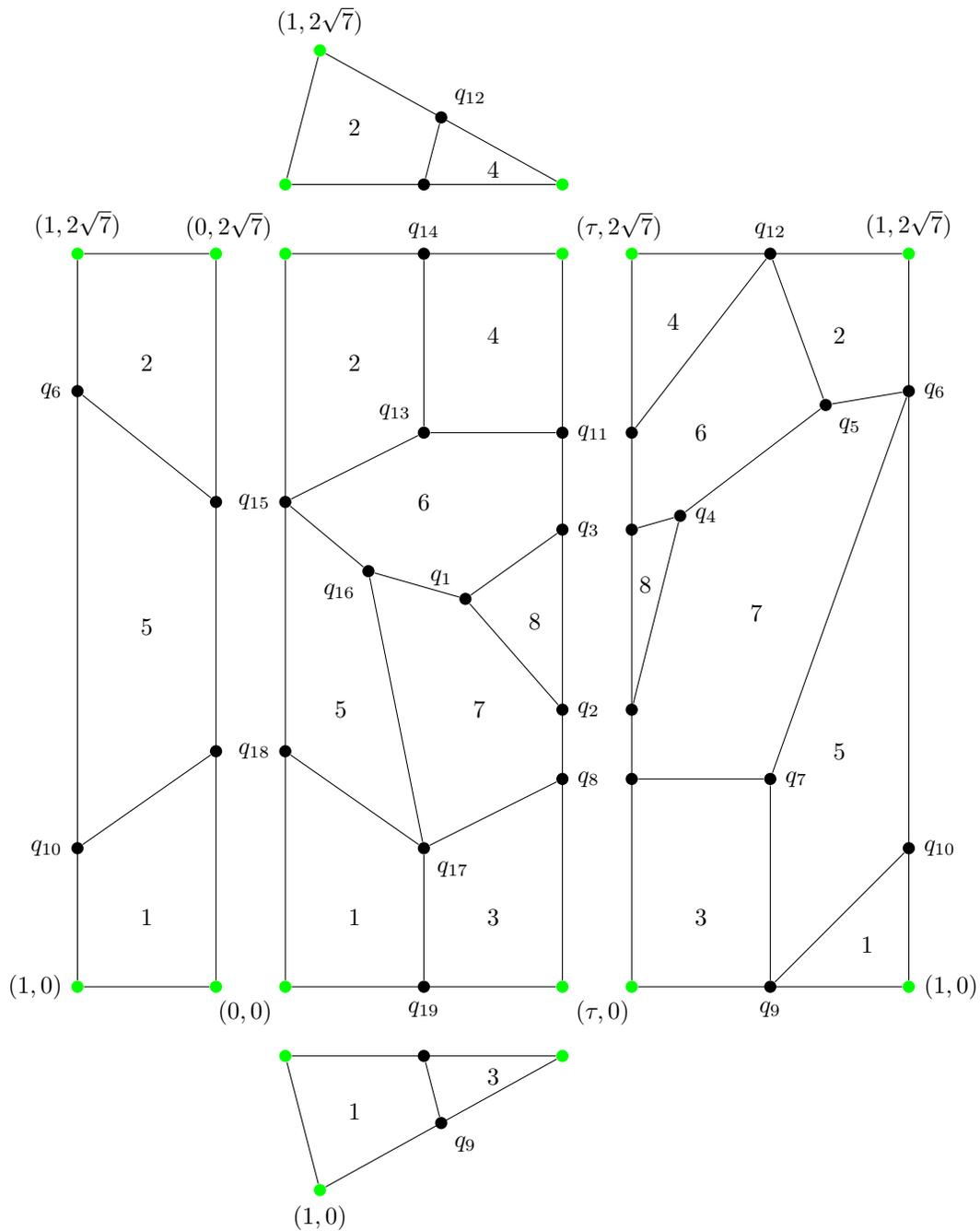
\end{center}

By inspection, we see that the $\Q[i\sqrt{7}]$-rational points in $D_\infty(7)$ with depth at most 7 are the following, in horospherical coordinates. We give in Lemma~\ref{depth2points} below a detailed justifcation for the points of depth 2.

\begin{itemize}
\item Depth 1: $(0,0)$, $(0,2\sqrt{7})$, $(\tau,0)$, $(\tau,2\sqrt{7})$ and $(1,\sqrt{7})$, all in the same $\Gamma_\infty(7)$-orbit;
\item Depth 2: $(\frac{\tau}{2},\frac{3}{2}\sqrt{7})$ in one $\Gamma_\infty(7)$-orbit and $(\frac{\tau+1}{2},\sqrt{7})$ in the other;
\item Depth 4:   $(0,\sqrt{7})$, $(\tau,\sqrt{7})$, $(1,0)$, $(1,2\sqrt{7})$ in one $\Gamma_\infty(7)$-orbit, $(\frac{\tau}{2},\frac{1}{2}\sqrt{7})$ in a second, $(\frac{\tau+1}{2},0)$, $(\frac{\tau+1}{2},2\sqrt{7})$ in a third, $(\frac{\tau+1}{4}, \frac{5\sqrt{7}}{4})$ in a fourth, and $(\frac{\tau+2}{4},\frac{3\sqrt{7}}{2})$ in a fifth.
  \item Depth 7: for each $k=1,...,6$, there is a $\Gamma_\infty(7)$-orbit containing the 3 points $(0,\frac{2k}{7}\sqrt{7})$, $(\tau,\frac{2k}{7}\sqrt{7})$, $(1,\frac{2k+7}{7}\sqrt{7})$ (with $2k+7$ taken mod 7).  
  \end{itemize}

Integral lifts of representatives of $\Gamma_\infty(7)$-orbits of these points are:

$$
\begin{array}{cccccc} p_0=\left[ \begin{array}{c} 0 \\ 0 \\ 1 \end{array}\right] & 
  p_{2,1}=\left[ \begin{array}{c} -\tau-2 \\ -1 \\ \tau-1 \end{array}\right] &
  p_{2,2}=\left[ \begin{array}{c} 2 \\ 1-\tau \\ -\tau \end{array}\right] \\
\\
  p_{4,1}=\left[ \begin{array}{c} i\sqrt{7} \\ 0 \\ 2 \end{array}\right]  &
  p_{4,2}=\left[ \begin{array}{c} \tau-1 \\ \tau \\ 2  \end{array}\right] &
  p_{4,3}=\left[ \begin{array}{c} -1 \\ \tau+1 \\ 2  \end{array}\right] &
  p_{4,4}=\left[ \begin{array}{c} -2\tau-1 \\ -1 \\ \tau-2  \end{array}\right] &
  p_{4,5}=\left[ \begin{array}{c}  2\tau-4 \\ \tau \\ \tau+1  \end{array}\right] \\
\\
  p_{7,1}=\left[ \begin{array}{c} -1 \\ 0 \\ i\sqrt{7} \end{array}\right]  &
  p_{7,2}=\left[ \begin{array}{c} -2 \\ 0 \\ i\sqrt{7} \end{array}\right]  &
  p_{7,3}=\left[ \begin{array}{c} -3 \\ 0 \\ i\sqrt{7} \end{array}\right]  &
  p_{7,4}=\left[ \begin{array}{c} -4 \\ 0 \\ i\sqrt{7} \end{array}\right]  &
  p_{7,5}=\left[ \begin{array}{c} -5 \\ 0 \\ i\sqrt{7} \end{array}\right]  &
  p_{7,6}=\left[ \begin{array}{c} -6 \\ 0 \\ i\sqrt{7} \end{array}\right]    
  \end{array}
$$

\begin{lem}\label{depth2points} The $\Q[i\sqrt{7}]$-rational points in $D_\infty(7)$ with depth 2 are exactly $(\frac{\tau}{2},\frac{3}{2}\sqrt{7})$ and $(\frac{\tau+1}{2},\sqrt{7})$.
\end{lem}

\Pf We illustrate the general procedure for finding points of depth $n$, then specialize to the present  case.

The depths at which there will be $E$-rational points are the natural numbers $n$ such that $|z|^2 = n$ has a solution $z\in\mathcal{O}_d$.  For $d=7$, the first few $n$'s are $1,2,4,7,\ldots$.  
Begin by assuming we have found all points of depth less than $n$.  In this case those are just the points at depth 1.

\begin{enumerate}

\item Find all $q = x+y\tau \in\mathcal{O}_d$ (up to multiplication by a unit) with $|q|^2 = n$.  From the geometry of numbers, we know that there are only finitely many such $q$.  In this case, there are two possibilities: either $q=\tau:=\frac{1+i\sqrt{7}}{2}$ or $q= \tau-1$.

\item Consider the standard lift  of a point of $\partial H_{\mathbb{K}}^2 \setminus \{ \infty \}$:
$$P = \left(\begin{array}{c}
\frac{-|z|^2+it}{2}\\z\\1
\end{array}\right)$$
Note that $t = b\sqrt{7}$ for some $b\in \Q$.  The calculation is more transparent if we rewrite the first coordinate:
$$\frac{-|z|^2+it}{2} = \frac{-|z|^2-b}{2}+b\tau$$

For $P$ to be the vector form of an $E$-rational point $p$ of depth $2$, it must satisfy the following:
\begin{enumerate}
\item $P$ does not have depth $1$ (in general, $P$ does not have depth less than $n$).  In other words, the coordinates of $P$ are not in $\mathcal{O}_7$.
\item $Pq$'s coordinates must be in $\mathcal{O}_7$ for some $q$ from step (1). 
\end{enumerate}
Next, we do some calculations to make sure (b) is satisfied.

\item Find all $z$ in the projection to $\mathbb{C}$ of $D_\infty$ such that $zq\in\mathcal{O}_7$.  If $q = \tau$ we can have $z = 0,1,\tau,\text{ or }\frac{1+\tau}{2}$.  If $q = \tau-1$, we can have $z = 0,1,\tau,\frac{\tau}{2}$.

\item For each possible $z$, find $|z|^2$ and compute $\left(\frac{-|z|^2-b}{2}+b\tau\right)q$.  Use this to list all $b$'s such that $(z,b)\in D_{\infty}$ and $\left(\frac{-|z|^2-b}{2}+b\tau\right)q\in\mathcal{O}_7$.
$$
\begin{array}{cccccc}
z & |z|^2 & q & \left(\frac{-|z|^2-b}{2}+b\tau\right)q & b\text{'s} & \text{point(s) in horo. coords}\\
\hline
0 & 0 & \tau & -2b+\frac{b\tau}{2} & 0,2 & (0,0), (0,2\sqrt{7})\\
&  & \tau-1 & -\frac{3b}{2}-\frac{b}{2}\tau & 0,2 & (0,0), (0,2\sqrt{7})\\
\hline
1 & 1 & \tau & 2b + \frac{b-1}{2}\tau & 1 & (1,\sqrt{7})\\
&  & \tau-1 & \frac{1-3b}{2} - \frac{1+b}{2}\tau & 1 & (1,\sqrt{7})\\
\hline
\tau & 2 & \tau & -2b-\frac{b-2}{2} & 0,2 & (0,0), (0,2\sqrt{7})\\
&  & \tau-1 & \frac{2-3b}{2} - \frac{2+b}{2}\tau &  0,2 & (0,0), (0,2\sqrt{7})\\ 
\hline
\frac{1+\tau}{2} & 1 & \tau & 2b + \frac{b-1}{2}\tau & 1 & \left(\frac{1+\tau}{2}, \sqrt{7}\right)\\
\hline
\frac{\tau}{2} & \frac{1}{2} & \tau -1 & \frac{1-6b}{4} - \frac{1+2b}{4}\tau & \frac{3}{2} & \left(\frac{\tau}{2},\frac{3\sqrt{7}}{2}\right)
\end{array}
$$

\item Get rid of all the ones that are level 1 (or from any previous level).  What you are left with are the level $2$ points: $\left\{\left(\frac{1+\tau}{2}, \sqrt{7}\right), \left(\frac{\tau}{2},\frac{3\sqrt{7}}{2}\right)\right\}$. 
\end{enumerate}
\EPf

{\bf Generators:}  The following elements $A_\alpha \in \Gamma(7)$ map the point $\infty=[1,0,0]^T$ to the corresponding $p_\alpha$ as above (for $\alpha=0; 2,1; 2,2; 4,1;...; 4,5; 7,1;...;7,6$) :

$$
\begin{array}{ccc} A_0=I_0=\left[ \begin{array}{ccc} 0 & 0 & 1 \\ 0 & -1 & 0 \\ 1 & 0 & 0 \end{array}\right] &
  
  A_{2,1}=\left[ \begin{array}{ccc} -\tau-2 & -5 & 5 \\ -1 & -1 & 0 \\ \tau-1 & i\sqrt{7} & -i\sqrt{7} \end{array}\right] &
  
  A_{2,2}=\left[ \begin{array}{ccc} 2 & -\tau & \tau -1 \\ 1-\tau & -2 & \tau \\ -\tau & \tau -1 & 2 \end{array}\right] \\
  \\
  A_{4,1}=\left[ \begin{array}{ccc} i\sqrt{7} & 0 & 4 \\ 0 & 1 & 0 \\ 2 & 0 & -i\sqrt{7} \end{array}\right] &
  
  A_{4,2}=\left[ \begin{array}{ccc}\tau -1 & \tau & 2 \\ \tau & 2 & 1-\tau \\ 2 & 1-\tau & -\tau \end{array}\right] &
  
  A_{4,3}=\left[ \begin{array}{ccc} -1 & \tau -2 & 2 \\ \tau + 1 & 3 & -\tau -1 \\ 2 & 2-\tau & -1  \end{array}\right] \\
  \\
    A_{4,4}=\left[ \begin{array}{ccc} -2\tau-1 & -4 & 4 \\ -1 & -1 & 0 \\ \tau-2 & i\sqrt{7} & -i\sqrt{7} \end{array}\right] &
  
  A_{4,5}=\left[ \begin{array}{ccc}2\tau -4 & 3+3\tau & 4-3\tau \\ \tau & 2 & 1-\tau \\ 1+\tau & 3-2\tau & -2-\tau \end{array}\right]  \\
  \\
  A_{7,1}=\left[ \begin{array}{ccc} -1 & 0 & 0 \\ 0 & 1 & 0 \\ i\sqrt{7} & 0 & -1  \end{array}\right] &
  A_{7,2}=\left[ \begin{array}{ccc} -2 & 0 & i\sqrt{7} \\ 0 & 1 & 0 \\ i\sqrt{7} & 0 & 3  \end{array}\right] &
  A_{7,3}=\left[ \begin{array}{ccc} -3 & 0 & i\sqrt{7} \\ 0 & 1 & 0 \\ i\sqrt{7} & 0 & 2  \end{array}\right] \\
  \\
  A_{7,4}=\left[ \begin{array}{ccc} -4 & 0 & 3i\sqrt{7} \\ 0 & 1 & 0 \\ i\sqrt{7} & 0 & 5  \end{array}\right] &
  A_{7,5}=\left[ \begin{array}{ccc} -5 & 0 & 3i\sqrt{7} \\ 0 & 1 & 0 \\ i\sqrt{7} & 0 & 4  \end{array}\right] &
  A_{7,6}=\left[ \begin{array}{ccc} -6 & 0 & i\sqrt{7} \\ 0 & 1 & 0 \\ i\sqrt{7} & 0 & 1  \end{array}\right] 
 
  \end{array}
$$

{\bf Relations:}  We obtain a complete set of relations between these generators by applying generators to points of depth at most 4 as described in part (5) of section~\ref{method}. The detailed steps of the computation can be found on the companion Sagemath Jupyter notebook in \cite{MCode}. The direct output is a presentation with 18 generators and 406 relations, which simplifies, thanks to Magma \cite{Mag}, to the presentation in Proposition~\ref{presentationd=7}. More specifically, we obtain that particular simplification by specifying a subset of the generators which must be preserved, using in this case the command {\tt Simplify (G: Preserve:=[1,4,5]);}.

To illustrate the steps involved we compute by hand the cycles and relations corresponding to triples $(a,b,c)$ for which the point $A_a(p_b)$ has depth at most 7, hence is of the form $\gamma_\infty (p_c)$ for some $\gamma_\infty \in \Gamma_\infty$ and $p_c \in D_\infty$ with depth at most 7. In the notation of part (5) of section~\ref{method}, this corresponds to $\gamma_\infty^1={\rm Id}$ and  $\gamma_\infty^3=\gamma_\infty$. The results are listed in Table~4, with the same notation as in the example of section~\ref{toyexample}.

\begin{table}\label{d=7table1}
  \begin{center}
    \vspace{-1cm}
      {\renewcommand{\arraystretch}{1.2}%
\begin{tabular}{|c|c|c|c|}
  \hline
  
$\mathbf{A.p}$ & $\mathbf{p'}$ & $\mathbf{A_{p'}^{-1}W^{-1}A A_p}$ & $\mathbf{W'}$  \\
\hline

$I_0 p_0$ & $\infty$ & $I_0^2$ & ${\rm Id}$\\
\hline

$A_{4,1} p_{4,1}$ & $\infty$ & $A_{4,1}^2$ & ${\rm Id}$\\
\hline

$A_{4,2} p_{4,2}$ & $\infty$ & $A_{4,2}^2$ & ${\rm Id}$\\
\hline

$A_{4,3} p_{4,3}$ & $\infty$ & $A_{4,3}^2$ & ${\rm Id}$\\
\hline

$R p_0$ & $p_0$ & $I_0RI_0$ & $R$\\
\hline

$I_0 p_{2,2}$ & $p_{4,3}$ & $A_{4,3}^{-1}T_1I_0A_{2,2}$ & $T_1$\\
\hline

$I_0 p_{4,1}$ & $p_{7,5}$ & $A_{7,5}^{-1}T_vI_0A_{4,1}$ & $R$\\
\hline

$I_0 p_{4,2}$ & $p_{2,1}$ & $A_{2,1}^{-1}T_vT_1I_0A_{4,2}$  & $RT_v^{-1}T_\tau^{-1}T_1^{-2}$\\
\hline

$I_0 p_{4,3}$ & $p_0$ & $A_{4,3}^{-1}I_0T_\tau T_1T_v^{-1}I_0$ & $T_1T_\tau$\\
\hline

$I_0 p_{7,1}$ & $p_0$ & $A_{7,1}^{-1}I_0T_v^{-1}I_0$ & $R$\\
\hline

$I_0 p_{7,2}$ & $p_{4,1}$ & $A_{7,2}^{-1}I_0T_v^{-1}A_{4,1}$ & ${\rm Id}$\\
\hline

$I_0 (0,2\sqrt{7})$ & $p_{7,6}$ & $A_{7,6}^{-1}T_v I_0T_vI_0$ & ${\rm Id}$\\
\hline

$I_0 (1,\sqrt{7})$ & $p_{2,1}$ & $A_{2,1}^{-1}T_v I_0T_1I_0$ & $T_1^2T_v^{-1}$\\
\hline

$I_0 (\tau,0)$ & $(\tau,0)$ & $(T_\tau I_0)^{-1}I_0T_\tau I_0$ & $R T_\tau$\\
\hline

$A_{2,1}p_0$ & $p_{7,5}$ & $A_{7,5}^{-1}A_{2,1}I_0$ & $RT_vT_1^{-1}$\\
\hline

$A_{2,2}p_0$ & $p_{4,2}$ & $A_{4,2}^{-1}A_{2,2}I_0$ & ${\rm Id}$\\
\hline

$A_{4,1}p_0$ & $p_{7,4}$ & $A_{7,4}^{-1}A_{4,1}I_0$ & ${\rm Id}$\\
\hline

$A_{4,1}p_{4,4}$ & $p_{2,2}$ & $A_{2,2}^{-1}T_\tau A_{4,1}A_{4,4}$ & $T_1^{-1}$\\
\hline

$A_{4,1}p_{7,3}$ & $p_0$ & $I_0T_v^{-1}A_{4,1}A_{7,3}$ & $R$\\
\hline

$A_{4,2}p_0$ & $p_{2,2}$ & $A_{2,2}^{-1}A_{4,2}I_0$ & ${\rm Id}$\\
\hline

$A_{4,2}p_{7,2}$ & $p_{2,2}$ & $A_{2,2}^{-1}T_1T_\tau ^{-1}A_{4,2}A_{7,2}$ & $T_\tau T_1^{-2}T_v  $\\
\hline

$A_{4,3}p_0$ & $p_0$ & $I_0T_1^{-1}T_\tau ^{-1}T_v A_{4,3}I_0$ & $T_1^{-1}T_\tau ^{-1}T_v$\\
\hline

$A_{4,4} p_0$ & $p_{7,4}$ & $A_{7,4}^{-1}A_{4,4}I_0$ & $RT_1^{-1}$\\
\hline

$A_{4,5} p_{2,2}$ & $p_{7,5}$ & $A_{7,5}^{-1}A_{4,5}A_{2,2}$ & $RT_v$\\
\hline

$A_{7,1}p_0$ & $p_0$ & $I_0A_{7,1}I_0$ & $RT_v^{-1}$\\
\hline

$A_{7,2}p_{4,1}$ & $p_0$ & $I_0A_{7,2}A_{4,1}$ & $T_v^{-1}$\\
\hline

$A_{7,2}p_{7,3}$ & $\infty$ & $A_{7,2}A_{7,3}$ & $R$\\
\hline

$A_{7,2}p_{7,4}$ & $p_{7,1}$ & $A_{7,1}^{-1}A_{7,2}A_{7,4}$ & $R$\\
\hline

$A_{7,3}p_0$ & $p_{4,1}$ & $A_{4,1}^{-1}A_{7,3}I_0$ & $RT_v$\\
\hline

$A_{7,3}p_{4,2}$ & $p_{2,1}$ & $A_{2,1}^{-1}T_1A_{7,3}A_{4,2}$ & $T_\tau ^{-1}T_1^2T_v$\\
\hline

$A_{7,3}p_{7,1}$ & $p_{7,4}$ & $A_{7,4}^{-1}A_{7,3}A_{7,1}$ & ${\rm Id}$\\
\hline

$A_{7,3}p_{7,2}$ & $\infty$ & $A_{7,3}A_{7,2}$ & $R$\\
\hline

$A_{7,4}p_0$ & $p_{4,1}$ & $A_{4,1}^{-1}A_{7,4}I_0$ & ${\rm Id}$\\
\hline

$A_{7,5}p_{2,2}$ & $p_{4,3}$ & $A_{4,3}^{-1}T_1T_v^{-1}A_{7,5}A_{2,2}$ & $RT_\tau ^{-2}T_1$\\
\hline

$A_{7,5}p_{4,1}$ & $p_0$ & $I_0T_v^{-1}A_{7,5}A_{4,1}$ & $R$\\
\hline

$A_{7,5} p_{4,4}$ & $p_0$ & $I_0T_v^{-1}T_1 A_{7,5}A_{4,4}$ & ${\rm Id}$\\
\hline

$A_{7,5}p_{7,3}$ & $p_{7,6}$ & $A_{7,6}^{-1}A_{7,5}A_{7,3}$ & ${\rm Id}$\\
\hline

$A_{7,5}p_{7,4}$ & $\infty$ & $A_{7,5}A_{7,4}$ & $RT_v$\\
\hline

$A_{7,6}p_0$ & $p_0$ & $I_0T_v^{-1}A_{7,6}I_0$ & $T_v$\\
\hline

$A_{7,6}p_{7,1}$ & $\infty$ & $A_{7,6}A_{7,1}$ & $RT_v$\\
\hline

$A_{7,6}p_{7,2}$ & $p_{7,5}$ & $A_{7,5}^{-1}A_{7,6}A_{7,2}$ & $R$\\
\hline

\end{tabular}
}
\end{center}
\caption{Action of generators on vertices for $d=7$}
\end{table}

\section{The Hurwitz quaternion modular group ${\rm PU}(2,1,\mathcal{H})$}

In this section we use the method described in Section~\ref{method} to compute a presentation for the Hurwitz modular group $\Gamma(\mathcal{H})={\rm PU}(2,1,\mathcal{H})$ (also denoted  ${\rm PSp}(2,1,\mathcal{H})$).
Recall that the Hurwitz integer ring is $\mathcal{H}=\Z[i,j,k,\sigma] \subset \quat$, with $\sigma=\frac{1+i+j+k}{2}$. The resulting presentation is unfortunately too large for Magma to handle directly (see comments at the end of this section), however it still allows us to obtain the following results. 

\begin{thm}\label{hurwitzabel} The abelianization of $\Gamma(\mathcal{H})$ is $\Z/3\Z$.
\end{thm}

\begin{thm}\label{hurwitzgens} $\Gamma(\mathcal{H})$ is generated by $\Gamma_\infty(\mathcal{H})$ and $I_0$.
\end{thm}

Theorem~\ref{hurwitzgens} was in fact stated in \cite{Ph} (Th. 4.4.2), but the very short proof given there is inadequate (it would indeed apply to any group containing $I_0$ and a nontrivial stabilizer of $\infty$). However, it is shown in \cite{Ph} that  $\Gamma_\infty(\mathcal{H})$ can be generated by only 3 elements, which combined with the above shows that $\Gamma(\mathcal{H})$ can be generated by 4 elements. 

We will use Theorem~\ref{hurwitzgens} in the proof of Theorem~\ref{hurwitzabel}, so for future reference we note that  
Theorem~\ref{hurwitzgens} follows from the result of computational simplification of a partial presentation by Magma. More specifically, we enter the Magma command {\tt Simplify(G: Preserve:=[1,2,3,4,5,6,7,8,9,10,11,12,13,14,15,16,17]);} applied to the group $G$ presented by all generators and the first thousand relations in our main presentation (see file {\tt QuaternionsTruncated1000.txt} at \cite{MCode}). This returns a presentation of $G$ with 17 generators, corresponding to 16 generators for $\Gamma_\infty(\mathcal{H})$ and $I_0$. Theorem~\ref{hurwitzgens} follows since $\Gamma(\mathcal{H})$ is a quotient of $G$.

\vspace{.2cm}

{\bf Presentation and fundamental domain for $\Gamma_\infty(\mathcal{H})$:}

In this section we study the action of  the cusp group $\Gamma_\infty(\mathcal{H})={\rm Stab}_{\Gamma(\mathcal{H})}(\infty)$ on $\partial_\infty {\rm H}_\quat ^2 \setminus \{ \infty \} \simeq \quat \times {\rm Im} \, \quat$.
The referee informed us that a fundamental domain for the action of $\Gamma_\infty(\mathcal{H})$ on $\partial_\infty {\rm H}_\quat ^2 \setminus \{ \infty \}$ had been found by Philippe in her thesis \cite{Ph}. In a previous version of this paper we only used a \emph{coarse fundamental domain}, i.e. a domain whose translates cover the space; using her fundamental domain simplified the arguments and computations in this section, so we thank the referee for drawing our attention to it. In fact Philippe gives in \cite{Ph} (Theorem 4.3.11, Proposition 4.3.22 and pp. 103--104) two fundamental domains, related by a sequence of cut-and paste operations. We denote these two domains by $D^1_\infty(\mathcal{H})$ and $D^2_\infty(\mathcal{H})$ and will in fact use both, as they have different geometric/combinatorial advantages. More specifially, $D^2_\infty(\mathcal{H})$ is simpler combinatorially as it is a product of a 4-simplex and a 3-cube, whereas $D^1_\infty(\mathcal{H})$ is a union of two such objects. However each of the two isometric pieces of  $D^1_\infty(\mathcal{H})$ is contained in a smaller Cygan ball centered at 0, which makes it easier to use in the covering argument below.

Recall that $\partial_\infty {\rm H}_\quat ^2 \setminus \{ \infty \} \simeq \quat \times {\rm Im} \, \quat$; we will refer to $\quat$ as the horizontal direction and ${\rm Im} \, \quat$ as the vertical direction. The vertices of the horizontal component of these domains consist in the following points in $\quat$ (following the notation in \cite{Ph}):

$$\begin{array}{ccccc}\label{basepoints}
  p_0=0, & p_3=i, & c_1=(1+i/3+j+k)/2, & c_2=(-1+i/3-j-k)/2, & \\
    \\
   q_0=i/2,  & q_{14}=(i+k)/2, & q_{18}=(i+j)/2, & q_{16}=(i-k)/2, & q_{20}=(i-j)/ 2.  
\end{array}
$$

\begin{prop}[\cite{Ph}]\label{philippedomains} Consider the subsets $S={\rm Hull}(p_0,p_3,c_1,q_{14},q_{18})$, $C_1={\rm Hull}(p_0,q_0,c_1,q_{14},q_{18})$ and $C_2={\rm Hull}(p_0,q_0,c_2,q_{16},q_{20})$ of $\quat$. Then $D^1_\infty(\mathcal{H})=(C_1 \cup C_2) \times [-1,1]^3$ and $D^2_\infty(\mathcal{H})=S \times [-1,1]^3$ are fundamental domains for the action of $\Gamma_\infty(\mathcal{H})$ on $\partial_\infty {\rm H}_\quat ^2 \setminus \{ \infty \}\simeq \quat \times {\rm Im} \, \quat$.
\end{prop}  

We now determine a presentation for $\Gamma_{\infty}(\mathcal{H})$ by more algebraic methods, understanding this group as a sequence of normal extensions of simpler subgroups. The presentation we give is highly redundant in terms of generators, but we will need this "geometrically complete" set of generators when we are required to algorithmically identify given elements of  $\Gamma_{\infty}(\mathcal{H})$ in terms of the generators. (Also note that \cite{Ph} does not identify a presentation corresponding to the above fundamental domains).

We use the notation of Equations~(\ref{heistrans}), (\ref{heisrot}) and  (\ref{conj}), namely $T_{(\zeta,v)}$, $R_w$, $C_w$ respectively denote Heisenberg translation by $(\zeta,v)$, Heisenberg rotation by $w$ and conjugation by $w$. Note that for any unit quaternion $w \in \mathcal{H}$ and purely imaginary  $q \in \mathcal{H}$, $\Gamma_\infty(\mathcal{H})$ contains the following Heisenberg translations:

\begin{eqnarray}
T_w=T_{(w,i+j+k)}=\left(\begin{array}{ccc}
1 & -\bar{w} & \frac{-1+i+j+k}{2} \\
0 & 1 & w \\
0 & 0 & 1 
\end{array}\right),
& & 
T_{v_q} = T_{(0,2q)}=\left(\begin{array}{ccc}
1&0 & q\\
0& 1 & 0\\ 
0 & 0 & 1
\end{array}\right)
\end{eqnarray}

\begin{lem}\label{inftypres}
$\Gamma_{\infty}(\mathcal{H})$ admits the following presentation:

$$\Gamma_{\infty}(\mathcal{H})=\left\langle \mathcal{S}_\infty \, | \, \mathcal{R}_\infty \right\rangle=\left\langle R_i, R_{\sigma}, C_i, C_{\sigma}, T_1, T_i, T_j, T_k, T_{\sigma}, T_{v_i}, T_{v_j}, T_{v_k} \, | \, \mathcal{A}_1, \mathcal{A}_2, \mathcal{A}_3, \mathcal{A}_4, \mathcal{A}_5\right\rangle$$
where the $\mathcal{A}_n$ are the following sets of relations: 

$$\mathcal{A}_1 = \{R_i^4, R_i^2R_{\sigma}^{-3}, (R_iR_{\sigma})^3, (R_i^{-1}R_{\sigma})^2R_i^{-1}R_{\sigma}^{-2}\}$$

$$\mathcal{A}_2 = \{C_i^2, C_{\sigma}^3, (C_iC_{\sigma})^3\}$$ 

$$\mathcal{A}_3 = \left\{
\begin{array}{c}\left[T_{v_w},T\right], \left[T_1,T_w\right] = T_{v_w}^{-2}, \left[T_1,T_\sigma\right] = T_{v_i}^{-1}T_{v_j}^{-1}T_{v_k}^{-1}, \vspace{.1cm} \\ \left[T_w,T_{\hat{w}}\right] = T_{v_{w\hat{w}}}^2, \left[T_w,T_{\sigma}\right] = T_{v_w}T_{v_{\hat{w}}}^{-1}T_{v_{w\hat{w}}}, T_1T_iT_jT_k = T_\sigma^2T_{v_i}T_{v_j}^{-1}T_{v_k}\end{array} \right\}$$
where $T$ runs over $T_1,T_i,T_j,T_k,T_{\sigma},T_{v_i},T_{v_j},T_{v_k}$, $w$ runs over  $i,j,k$, and $\hat{i} = j, \hat{j} = k, \hat{k} = i$.

$$\mathcal{A}_4 = \left\{ C_iR_iC_i = R_i, C_iR_{\sigma}C_i = R_iR_{\sigma}^4R_i, C_{\sigma}R_iC_{\sigma}^{-1} = R_{\sigma}R_iR_{\sigma}^{-1}, C_{\sigma}R_{\sigma}C_{\sigma}^{-1} = R_{\sigma}\right\}$$

$$\mathcal{A}_5 = \{GTG^{-1} = E_{G,T} \},$$
where $G$ runs over $R_i.R_{\sigma},C_i,C_{\sigma}$, $T$ runs over $T_1,T_i,T_j,T_k,T_{\sigma},T_{v_i},T_{v_j},T_{v_k}$, and $E_{G,T}$ is the entry in the $G$ column and $T$ row of the table below:
$$
\begin{array}{c|cccc}
& R_i & R_{\sigma} & C_i & C_{\sigma}\\\hline
T_1 & T_i & T_{\sigma} & T_1T_{v_j}^{-1}T_{v_k}^{-1} & T_1\\
T_i & T_1^{-1}T_{v_i}T_{v_j}T_{v_k} & T_1^{-1}T_k^{-1}T_{\sigma}T_{v_j}T_{v_k}^2 & T_iT_{v_j}^{-1}T_{v_k}^{-1} & T_j\\
T_j & T_k & T_1^{-1}T_i^{-1}T_{\sigma}T_{v_i}^2T_{v_k} & T_j^{-1}T_{v_i} & T_k\\
T_k & T_j^{-1}T_{v_i}T_{v_j}T_{v_k} & T_1^{-1}T_j^{-1}T_{\sigma}T_{v_i}T_{v_j}^2 & T_k^{-1}T_{v_i} & T_i\\
T_{\sigma} & T_{\sigma}^{-1}T_iT_kT_{v_j}T_{v_k}^{-1} & T_1^{-1}T_{\sigma} & T_j^{-1}T_k^{-1}T_{\sigma}T_{v_j} & T_{\sigma}\\
T_{v_i} & T_{v_i} & T_{v_i} & T_{v_i} & T_{v_j}\\
T_{v_j} & T_{v_j} & T_{v_j} & T_{v_j}^{-1} & T_{v_k}\\
T_{v_k} & T_{v_k} & T_{v_k} & T_{v_k}^{-1} & T_{v_i}
\end{array}
$$

\end{lem}

\Pf To obtain a presentation for $\Gamma_{\infty}(\mathcal{H})$, we identify 3 of its subgroups, observe that some of them normalize each other, and build up the presentation via a sequence of extensions using the following procedure.  Suppose that $G$ is a group with subgroups $N$ and $K$ where $N$ is normal in $G$ and $G$ is an extension of $N$ by $K$.  Suppose also that we know presentations for $N$ and $K$, $N = \langle S_N\, | \, R_N\rangle$ and $K = \langle S_K\, |\, R_K\rangle$.  Then $G$ admits the presentation
$\langle S_N\cup S_K \, | \, R_N\cup R_K\cup R \rangle$, where the set $R$ consists of relations of the form $knk^{-1} = n'$ where $k$ runs over all elements of $S_K$, $n$ runs over the elements of $S_N$, and $n'\in N$ is expressed as a word in the generators $S_N$.

The three subgroups of $\Gamma_{\infty}(\mathcal{H})$ we identify are the rotation, conjugation, and translation subgroups.  

The rotation subgroup consists of all Heisenberg rotations $R_w$ where $w$ is a Hurwitz integral unit quaternion.  It is isomorphic to the binary tetrahedral group, which has order 24 (see e.g. \cite{CoS}).  It admits the presentation
$$\left\langle R_i, R_\sigma\, |  \, \mathcal{A}_1\right\rangle$$

The conjugation subgroup consists of all conjugations by unit quaternions.  Elements of this group also correspond to Hurwitz unit integral quaternions, only $C_w$ acts the same as $C_{-w}$.  Thus, this group is isomorphic to the quotient of the binary tetrahedral group by $-1$, which is the tetrahedral group (or the alternating group on 4 elements).  It admits the presentation
$$\left\langle C_i, C_\sigma\, | \, \mathcal{A}_2\right\rangle$$

The translation subgroup consists of all Heisenberg translations.  It admits the presentation
$$\left\langle T_1, T_i, T_j, T_k, T_\sigma, T_{v_i},T_{v_j},T_{v_k} \, | \, \mathcal{A}_3\right\rangle$$

The rotation subgroup is normalized by the conjugation subgroup.  The extension of the rotation subgroup by the conjugation subgroup is a finite group of order 288.  We obtain the four relations in $\mathcal{A}_4$ by conjugating $R_i$ and $R_{\sigma}$ by $C_i$ and $C_{\sigma}$.

The translation subgroup is normalized by the rotation-conjugation subgroup.  The conjugates of the translation generators by the rotation and conjugation generators are listed in the table.  The translations $T$ are along the left side, the rotations/conjugations $G$ are along the top, and the table entry $E_{G,T}$ is the element $GTG^{-1}$ written as a word in the translation generators.  $\mathcal{A}_5$ contains these relations. \EPf

\begin{cor}\label{inftyabel} The abelianization of $\Gamma_{\infty}(\mathcal{H})$ is $(\Z/3\Z)^2$, generated by (the images of) $R_\sigma$ and $C_\sigma$.
\end{cor}

{\bf Covering depth and $\Q[i,j,k]$-rational points in $D_\infty(\mathcal{H})$:}

We denote as before $B\left((\zeta,v),r\right)$ the open extended Cygan ball centered at $p=(\zeta,v) \in \partial_\infty {\rm H}_\quat^2$ with radius $r$ (see Equation~\ref{extcygan} for the definition of the extended Cygan metric $d_{XC}$). Recall that $u(n)=\frac{2}{\sqrt{n}}$ is the height at which balls of depth $n$ appear, in the sense of Corollary~\ref{cordepthradius}. 

\begin{lem}\label{hurwitzcover} Let $u=u(5)+\varepsilon=0.89443$ and $H_u$ the horosphere of height $u$ based at $\infty$. Then the prism $D^1_\infty(\mathcal{H})\times \{ u \}$ is covered by the intersections with $H_u$ of the following 17 extended Cygan balls of depth 1:  $B\left((0,0),\sqrt{2}\right)$, $B\left((\sigma,\pm i \pm j \pm k),\sqrt{2}\right)$, $B\left((\sigma-1-j-k,\pm i \pm j \pm k),\sqrt{2}\right)$.
\end{lem}

\Pf Recall from Proposition~\ref{philippedomains}, using the notation from \cite{Ph}, that $D^1_\infty(\mathcal{H})=(C_1 \cup C_2) \times [-1,1]^3  \subset \mathbb{H} \times {\rm Im} \, \quat$, where $C_1$ and $C_2$ are the 4-simplices $C_1={\rm Hull}\left(p_0,q_0,q_{14},q_{18},c_1\right) \subset \mathbb{H}$ and $C_2={\rm Hull}\left(p_0,q_0, q_{16},q_{20},c_2\right) \subset \mathbb{H}$.\\

{\bf Claim:} $C_1 \times [0,1]^3$ is contained in $B_0 \, \cup \, B_\sigma$, where $B_0=B\left((0,0),\sqrt{2}\right)$ and $B_\sigma=B\left((\sigma,i+j+k),\sqrt{2}\right)$. \\

Using coordinates $(x_1+x_2i+x_3j+x_4k,t_1i+t_2j+t_3k)$ on $\quat \times {\rm Im}\, \quat$, we separate $C_1 \times [0,1]^3$ into 2 pieces by the piecewise linear hypersurface $A=\{(x_1+x_2i+x_3j+x_4k,t_1i+t_2j+t_3k) \in \quat \times {\rm Im}\, \quat \, | \, T= L(X)\}$, where we denote $T=t_1+t_2+t_3 \, , \, X=x_1+x_2+x_3+x_4$, and $L$ is the piecewise linear function defined by:

$$ L(X) \ = \ \left\{\begin{array}{rcl}
2.7 & {\rm if} & X \in [0,1/2] \\
-1.4 X + 3.4 & {\rm if} & X \in [1/2,1]  \\
-0.9 X + 2.9 & {\rm if} & X \in [1, 5/3]. 
\end{array}\right.
$$.

Heuristically, the coordinate $X$ (resp. $T$) measures the ($l^1$-) distance from the origin $(0,0)$ in the horizontal (resp. vertical) direction, and the choice of the hypersurface $A$ was inspired by the position of the 40 vertices of $C_1 \times [0,1]^3$ relative to the Cygan balls $B\left((0,0),\sqrt{2}\right)$ and $B\left((\sigma,i+j+k),\sqrt{2}\right)$, see Figure~\ref{cubes} and below. 

 The claim is then verified by showing that $C_1 \times [0,1]^3 \cap A^- \subset B\left((0,0),\sqrt{2}\right)$ and $C_1 \times [0,1]^3 \cap A^+ \subset B\left((\sigma,i+j+k),\sqrt{2}\right)$, where $A^{\pm}$ denote the 2 half-spaces bounded by $A$, with $(0,0) \in A^-$ and $(\sigma,i+j+k) \in A^+$. This is done as in the proof of Lemma~\ref{d=7cover},  by checking the vertices numerically using Equation~(\ref{extcygan}), then extending the result to their affinely convex hull by Lemma~\ref{cyganconvex}. 
 
 Since $A$ is only piecewise linear, we first subdivide $A^+$ and $A^-$ into affinely convex pieces as follows. First we add the base vertices $s_1,...,s_5 \in \quat$ defined as the respective intersections of the affine segments $[p_0,q_{18}],[p_0,c_1], [p_0,q_{14}]$ with the hyperplane $\{ X=1/2 \}$, and $[q_0,c_1],[p_0,c_1]$ with the hyperplane $\{ X=1\}$ (see Figure~\ref{cubes}). Explicitly:
 
 $$s_1=\frac{i+j}{4}, \ s_2=\frac{3+i+3j+3k}{20}, \ s_3=\frac{i+k}{4}, \  s_4=\frac{3+5i+3j+3k}{14}, \ s_5=\frac{3+i+3j+3k}{10}.
 $$

 Since $C_1={\rm Hull}\left(p_0,q_0,q_{14},q_{18},c_1\right)$ with $X(p_0)=0$, $X(q_0)=1/2$, $X(q_{14})=X(q_{18})=1$ and $X(c_1)=5/3$, this ensures that: 
 $$\begin{array}{l}
 C_1\cap \{ X \leqslant 1/2 \} ={\rm Hull}\left(p_0,q_0,s_1,s_2,s_3\right), \\
 C_1\cap \{ 1/2 \leqslant X \leqslant 1 \} ={\rm Hull}\left(q_0,q_{14},q_{18}, s_1,s_2,s_3,s_4,s_5\right), \\ 
 C_1\cap \{ 1 \leqslant X \leqslant 5/3 \} ={\rm Hull}\left(p_0,q_0,s_1,s_2,s_3\right).
 \end{array}$$ 
 
 This produces a total of 10 vertices in the horizontal factor $\quat$, above each of which lies a 3-cube in the vertical factor, spanned by $0$, $i$, $j$, $k$, $i+j$, $i+k$, $j+k$, $i+j+k$.  We check numerically, using Equation~(\ref{extcygan}) that each of these 80 vertices lies in $B_0,B_\sigma$, or both as indicated in Figure~\ref{cubes}. Finally, we check that each point of intersection of an edge of a cube with the hypersurface $A$ lies in both $B_0$ and $B_\sigma$ (of course, the values appearing in the definition of the function $L$ above were chosen to satisfy this property). Explicitly, these intersection points comprise:
 \begin{itemize}
 \item 3 points at height $T=2.7$ on the top 3 edges of the cube above $p_0$,
 \item 12 points at height $T=2.7$ on the top 3 edges of the cubes above $q_0,s_1,s_2,s_3$,
 \item the 12 vertices at height $T=2$ (with vertical coordinates $i+j, i+k, j+k$) above $q_{18},s_4,s_5,q_{14}$,
 \item and 6 points at height $T=1.4$ on the middle 6 edges of the cube above $c_1$.
 \end{itemize}
 
For example, the relevant point above $c_1=(1+i/3+j+k)/2$ on the vertical edge $[j,i+j]$ has coordinates $(c_1,0.4i+j)=(1/2,1/6,1/2,1/2,0.4,1,0)$ and satisfies:

$$d_{XC}\left((c_1,0.4i+j,u),(0,0,0)\right) \simeq  1.410 < \sqrt{2},
$$

$$d_{XC}\left((c_1,0.4i+j,u),(\sigma,i+j+k,0)\right) \simeq 1.394 < \sqrt{2}. 
$$

All other computations are similar, and the claim follows. 
 \begin{figure}
 \centering
 \begin{tikzpicture}[scale=0.8]\label{cubes}

\begin{scope}[scale = .6]
\node[label=below:{$p_0$}] at (0,0) (p0){};
\cube
\dotify{\reddot}{\reddot}{\reddot}{\reddot}{\reddot}{\reddot}{\reddot}{\redbluedot}
\end{scope}

\begin{scope}[shift = {(4,5))}, scale = .6]
\node[label=below:{$q_0$}] at (0,0) (q0){};
\cube
\dotify{\reddot}{\reddot}{\redbluedot}{\reddot}{\redbluedot}{\reddot}{\redbluedot}{\bluedot}
\end{scope}

\begin{scope}[shift = {(8,5))}, scale = .6]
\node[label=below:{$q_{18}$}] at (0,0) (q18){};
\cube
\dotify{\redbluedot}{\reddot}{\redbluedot}{\redbluedot}{\redbluedot}{\redbluedot}{\redbluedot}{\bluedot}
\end{scope}

\begin{scope}[shift = {(8,-5))}, scale = .6]
\node[label=below:{$q_{14}$}] at (0,0) (q14){};
\cube
\dotify{\redbluedot}{\redbluedot}{\redbluedot}{\reddot}{\redbluedot}{\redbluedot}{\redbluedot}{\bluedot}
\end{scope}

\begin{scope}[shift = {(12,0))},scale = .6]
\node[label=below:{$c_1$}] at (0,0) (c1){};
\cube
\dotify{\reddot}{\redbluedot}{\reddot}{\redbluedot}{\bluedot}{\bluedot}{\bluedot}{\bluedot}
\end{scope}

\begin{scope}[shift = {(4,0))}, scale = .6]
\node[label=below:{$s_2$}] at (0,0) (s2){};
\cube
\dotify{\reddot}{\reddot}{\reddot}{\reddot}{\redbluedot}{\redbluedot}{\redbluedot}{\redbluedot}
\end{scope}

\begin{scope}[shift = {(4,2.5))}, scale = .6]
\node[label=below:{$s_1$}] at (0,0) (s1){};
\cube
\dotify{\reddot}{\reddot}{\redbluedot}{\redbluedot}{\redbluedot}{\redbluedot}{\redbluedot}{\bluedot}
\end{scope}

\begin{scope}[shift = {(4,-2.5))}, scale = .6]
\node[label=below:{$s_3$}] at (0,0) (s3){};
\cube
\dotify{\reddot}{\redbluedot}{\redbluedot}{\reddot}{\redbluedot}{\redbluedot}{\redbluedot}{\bluedot}
\end{scope}

\begin{scope}[shift = {(8,2.5))}, scale = .6]
\node[label=below:{$s_4$}] at (0,0) (s4){};
\cube
\dotify{\reddot}{\redbluedot}{\redbluedot}{\redbluedot}{\redbluedot}{\redbluedot}{\redbluedot}{\bluedot}
\end{scope}

\begin{scope}[shift = {(8,0))}, scale = .6]
\node[label=below:{$s_5$}] at (0,0) (s5){};
\cube
\dotify{\reddot}{\redbluedot}{\reddot}{\redbluedot}{\redbluedot}{\redbluedot}{\redbluedot}{\bluedot}
\end{scope}

\draw[dotted] (p0) -- (q0)
(p0) -- (q18)
(p0) -- (q14)
(p0) -- (c1)
(q0) -- (c1)
(q18) -- (c1)
(q14) -- (c1);

\end{tikzpicture}
\caption{Vertices of the prism $C_1 \times [0,1]^3$ ($p_0,q_0,q_{14},q_{18},c_1$) and its intersections with the level sets $\{X=1/2\}$ ($s_1,s_2,s_3$) and $\{X=1\}$ ($s_4,s_5$). The vertices of the horizontal base are labelled by name, above each lies a vertical 3-cube spanned by $0$, $i$, $j$, $k$, $i+j$, $i+k$, $j+k$, $i+j+k$, ordered here from bottom to top and left to right. A vertex is colored white if it belongs to $B\left((0,0),\sqrt{2}\right)$, black if  it belongs to $B\left((\sigma,i+j+k),\sqrt{2}\right)$, and black and white if it belongs to both. Cubes in a common column belong to a common level set of $X=x_1+x_2+x_3+x_4$.}
\end{figure}

We now complete the proof of Lemma~\ref{hurwitzcover} from the claim.
Replacing $C_1$ by $C_2$ by negating the first, third and fourth horizontal coordinates gives that the other subprism $C_2 \times [0,1]^3$ is contained in $B\left((0,0),\sqrt{2}\right) \, \cup \, B\left((\sigma-1-j-k,i+j+k),\sqrt{2}\right)$, since that transformation is an (affine) isometry of the Cygan metric.

Finally, by likewise negating some or all of the vertical coordinates we get that the 14 other subprisms $C_1 \times [-1,0] \times [0,1]^2$, $C_2 \times [-1,0] \times [0,1]^2$....  are contained in the extended Cygan balls obtained by switching the corresponding signs of the vertical coordinates of the center. This results in the 17 Cygan balls in the statement of Lemma~\ref{hurwitzcover}. \EPf

\begin{cor} The covering depth of $\Gamma(\mathcal{H})$ is at most 4.
\end{cor}  

We now use the other fundamental domain to list the points up to depth 4, namely $D^2_\infty(\mathcal{H}) = S\times (-1,1]^3$, where $S$ is the simplex
$S={\rm Hull}\left(0,\frac{i+j}{2},\frac{i+k}{2},\frac{1}{2}(1+\frac{i}{3}+j+k),i \right) \subset \mathbb{H}$.

\begin{lem}
The $\Q[i,j,k]$-rational points in $D^2_\infty(\mathcal{H})$ with depth at most 4 are, in horospherical coordinates:

\begin{itemize}
\item Depth 1: $(0,0)$;
\item Depth 2: $(0,i+j)$;
\item Depth 3: $(0,\frac{\pm 2i \mp 2j \pm 2k}{3})$, $(\frac{1+3i+j+k}{6},\frac{3i-j-k}{3})$, $(\frac{1+3i+j+k}{6},\frac{-i+3j+k}{3})$, $(\frac{1+i+j+3k}{6},\frac{3i-j-k}{3})$ (5 points total);
\item Depth 4: $(0,i)$, $(0,i+j+k)$, $(\frac{i+j}{2},\frac{\pm i \pm j \pm k}{2})$ (10 points total);
\end{itemize}
and all of their $\Gamma_\infty$-translates.

\end{lem}

\Pf We use the general procedure outlined in Lemma \ref{depth2points}.

Let $\Gamma_\infty^{fin}\subseteq \Gamma_\infty$ be the finite subgroup of order 288 described in Lemma \ref{inftypres}.  $\Gamma_\infty^{fin}$ describes two different actions by units $u\in\mathcal{H}$ on $\partial_\infty H_\mathbb{H}^2$.  One is the action by left multiplication on the first coordinate
$$R_u:p = (z,t) \mapsto (uz,t)$$
The other is by conjugation
$$C_u:p = (z,t)\mapsto (uzu^{-1},utu^{-1})$$

Consider the standard lift $P$ of a point $p\in \partial H_{\mathbb{K}}^2 \setminus \{ \infty \}$:
$$P = \left(\begin{array}{c}
\frac{-|z|^2+it_1+jt_2+kt_3}{2}\\z\\1
\end{array}\right)$$

For $P$ to be the vector form of a $\Q(i,j,k)$-rational point $p$ of depth $n$, it must satisfy the following:
\begin{enumerate}
\item $P$ does not have depth less than $n$ 
\item $Pq$'s coordinates must be in $\mathcal{H}$ for some $q$ with $|q|^2 = n$.  Write $z = z_1+iz_2+jz_3+kz_4$ and $q=q_1+iq_2+jq_3+kq_4+$.  Then we have
\begin{equation}\label{coordszq}
\begin{array}{c}
z_1q_1-z_2q_2-z_3q_3-z_4q_4,\\
z_1q_2+z_2q_1+z_3q_4-z_4q_3,\\
z_1q_3+z_3q_1-z_2q_4+z_4q_2,\\
z_1q_4+z_4q_1+z_2q_3-z_3q_2
\end{array} \in \Z,\text{ or}\in\Z+\frac{1}{2}
\end{equation}
and
\begin{equation}\label{coordsvertq}
\begin{array}{c}
\dfrac{1}{2}(-|z|^2q_1-t_1q_2-t_2q_3-t_3q_4),\\
\dfrac{1}{2}(-|z|^2q_2+t_1q_1+t_2q_4-t_3q_3),\\
\dfrac{1}{2}(-|z|^2q_3+t_2q_1-t_1q_4+t_3q_2),\\
\dfrac{1}{2}(-|z|^2q_4+t_3q_1+t_1q_3-t_2q_2)
\end{array} \in \Z,\text{ or}\in\Z+\frac{1}{2}
\end{equation}

\end{enumerate}
To list the points of depth $n$, we first assume we have already listed the points of depth less than $n$, and then we follow these steps:
\begin{enumerate}
    \item Find all $q\in\mathcal{H}$ with $|q|^2 = n$.  From the geometry of numbers, we know that there will only be finitely many $q$'s.  Even better, we only need to consider one representative for every $\Gamma_\infty^{fin}$ orbit of $q$'s under the $\Gamma_\infty^{fin}$ action on $\mathcal{H}$.  The reason for this is that if $u$ is a unit, then
    \begin{itemize}
        \item if $Pqu$ has coordinates in $\mathcal{H}$ then so does $Pq$.
        \item if $Puq$ has coordinates in $\mathcal{H}$ then $C_u^{-1}Pq = u^{-1}Puq$ also has coordinates in $\mathcal{H}$.
        \item if $Puqu^{-1}$ has coordinates in $\mathcal{H}$ then $C_u^{-1}Pq$ also has coordinates in $\mathcal{H}$.
    \end{itemize}
    Thus, we will always find a point in the $\Gamma_\infty^{fin}$-orbit of $P$ using our chosen $q$.
    
    \item For each $q$, find all potential first coordinates $z$, which are solutions to \eqref{coordszq}.  Like in the previous step, we only need one $z$ representing each $\Gamma_\infty^{fin}$-orbit.  We pick $z$ inside the fundamental domain $C$, and if necessary we modify our choice of $q$ by a unit so that $zq\in\mathcal{H}$.
    
    \item For each $z$, find all solutions to \eqref{coordsvertq} with $-1\leq t_1, t_2, t_3<1$.  Keep only those solutions that are not on the list of points of depth less than $n$.
\end{enumerate}

In practice, here's how that goes:

{\bf Depth 1:} In step 1, $q$ must be a unit, and there is a single orbit of $q$'s represented by $1$.  In step 2, we must have $z\in\mathcal{H}\cap D_\infty^h(\mathcal{H})$, and so $z=0$ is the only solution to \eqref{coordszq}.  

It follows that $|z|^2 = 0$, and so in step 3, the solutions to \eqref{coordsvertq} satisfying the inequalities are $(t_1,t_2,t_3)\in\Z^3$ with $0\leq t_1,t_2,t_3\leq 2$ and $t_1,t_2,t_3$ all even.  This gives us exactly the depth 1 point listed above.

{\bf Depth 2:}  In step 1, there is one orbit of $q$'s, represented by $1+i$.  In step 2, we find $z\in C$ by solving \eqref{coordszq} which specializes to 
$$z_1-z_2, z_1+z_2, z_3+z_4, z_4-z_3 \in \Z, \text{ or}\in\Z+\frac{1}{2}$$
Using some basic algebra, we conclude that 
$$2z_1,2z_2,2z_3,2z_4\in \Z, \text{ or}\in\Z+\frac{1}{2}$$
The possible values of $z$ in $C$ are then represented by $0,\frac{i}{2},\frac{i+j}{2}$.  

The values of $|z|^2$ corresponding to these $z$'s are $0,\frac{1}{4},\frac{1}{2}$ respectively.  In step 3, \eqref{coordsvertq} specializes to 
$$\dfrac{1}{2}(-|z|^2-t_1), \dfrac{1}{2}(-|z|^2+t_1), \dfrac{1}{2}(-t_2+t_3), \dfrac{1}{2}(t_3-t_2)
 \in \Z,\text{ or}\in\Z+\frac{1}{2}$$
With $z=0$, we find that taking two of $t_1,t_2,t_3$ equal to 1 and the other equal to 0 gives a solution.  These solutions are all new points of depth 2, and are $\Gamma_\infty^{fin}$-equivalent (by conjugation).  This is the depth 2 point listed above.

For the other two $z$'s, with $|z|^2 = \frac{1}{2}$ and $\frac{1}{4}$, there are no solutions to \eqref{coordsvertq}.

{\bf Depth 3:}
In step 1, there is one orbit of $q$'s, represented by $1-j-k$.  In step 2, we find $z$ by solving \eqref{coordszq} which specializes to 
$$z_1+z_3+z_4, z_2-z_3+z_4, -z_1+z_3+z_4, -z_1+z_4-z_2 \in \Z, \text{ or}\in\Z+\frac{1}{2}$$
Using some basic algebra, we conclude that 
$$3z_1,3z_2,3z_3,3z_4\in \Z, \text{ or}\in\Z+\frac{1}{2}$$
The possible values of $z$ in $C$ are then represented by 
$$0,\frac{i}{3},\frac{i+j}{3},\frac{1+3i+j+k}{3}$$

The values of $|z|^2$ corresponding to these $z$'s are $0, \frac{1}{9}, \frac{2}{9},\text{ and }\frac{1}{3}$ respectively.  In step 3, \eqref{coordsvertq} specializes to 
$$\dfrac{1}{2}(-|z|^2+t_2+t_3), \dfrac{1}{2}(t_1-t_2+t_3), \dfrac{1}{2}(|z|^2+t_2+t_1), \dfrac{1}{2}(|z|^2+t_3-t_1) \in \Z, \text{ or}\in\Z+\frac{1}{2}$$
Using some basic algebra, we find that 
$$\frac{3|z|^2}{2}, \frac{3t_1}{2}, \frac{3t_2}{2}, \frac{3t_3}{2}\in \Z, \text{ or}\in\Z+\frac{1}{2}$$
meaning that $3|z|^2\in \Z$ and its parity is the same as that of $3t_1,3t_2,3t_3$.  Therefore, we see that $|z|^2=0$ or $\frac{1}{3}$. We then solve \eqref{coordsvertq}, and in addition to finding the depth 1 point we already have, we find five new solutions at depth 3.

{\bf Depth 4:}

In step 1, there is one orbit of $q$'s, represented by $2$.  In step 2, we find $z$ by solving \eqref{coordszq} which specializes to 
$$2z_1,2z_2,2z_3,2z_4\in \Z, \text{ or}\in\Z+\frac{1}{2}$$
The possible values of $z$ in $C$ are then represented by $$0,\frac{1}{2}, \frac{i+j}{2}$$

The values of $|z|^2$ corresponding to these $z$'s are $0, \frac{1}{4}, \text{ and }\frac{1}{2}$ respectively.  In step 3, \eqref{coordsvertq} specializes to 
$$|z|^2,t_1,t_2,t_3 \in \Z, \text{ or}\in\Z+\frac{1}{2}$$
Thus, we must have $z=0$ or $z=\frac{i+j}{2}$.  When we solve for $t_1,t_2,t_3$, we recover all the depth 1 and 2 solutions we already found, as well as ten new points of depth 4.  \EPf




Integral lifts of representatives of $\Gamma_\infty(\mathcal{H})$-orbits of these points are:

$$
\begin{array}{cccccc} p_0=\left[ \begin{array}{c} 0 \\ 0 \\ 1 \end{array}\right] & 
  p_{2}=\left[ \begin{array}{c} -1 \\ 0 \\ i+j \end{array}\right] &
  p_{3,1}=\left[ \begin{array}{c} 1 \\ 0 \\ -i+j-k\end{array}\right] &
  p_{3,2}=\left[ \begin{array}{c} 1 \\ 0 \\ i-j+k\end{array}\right] \\
  \\
  p_{3,3}=\left[ \begin{array}{c} 1 \\ \frac{1-i-j+k}{2} \\ \frac{-1-3i+j+k}{2} \end{array}\right] &
  p_{3,4}=\left[ \begin{array}{c} 1 \\ -k \\  \frac{-1-i-3j-k}{2} \end{array}\right] &
  p_{3,5}=\left[ \begin{array}{c} 1 \\ \frac{1-i+j-k}{2}\\\frac{-1-i-j-3k}{2} \end{array}\right] &
  p_{4,1}=\left[ \begin{array}{c} i \\ 0 \\ 2 \end{array}\right] \\
  \\
  p_{4,2}=\left[ \begin{array}{c} i+j+k\\ 0 \\ 2  \end{array}\right] &
  p_{4,3}=\left[ \begin{array}{c} 1 \\ 1-i \\ -1-i-j-k  \end{array}\right] &
  p_{4,4}=\left[ \begin{array}{c} 1 \\ -i-k \\ -1+i-j-k  \end{array}\right] &
  p_{4,5}=\left[ \begin{array}{c}  1 \\ -i+k \\ -1-i+j-k  \end{array}\right] \\
\\
  p_{4,6}=\left[ \begin{array}{c} 1 \\ -1-i \\ -1+i+j-k  \end{array}\right] &
   p_{4,7}=\left[ \begin{array}{c} 1 \\ 1-j \\ -1-i-j+k  \end{array}\right] &
    p_{4,8}=\left[ \begin{array}{c} 1 \\ -j-k \\ -1+i-j+k  \end{array}\right] &
     p_{4,9}=\left[ \begin{array}{c} 1 \\ -j+k \\ -1-i+j+k  \end{array}\right] \\
     \\
      p_{4,10}=\left[ \begin{array}{c} 1 \\ -1-j \\ -1+i+j+k  \end{array}\right] 
  \end{array}
$$
{\bf Generators:}  The following elements $A_\alpha \in \Gamma(\mathcal{H})$ map the point $\infty=[1,0,0]^T$ to the corresponding $p_\alpha$ as above (for $\alpha=0; 2; 3,1;...;3,5; 4,1;...; 4,10$) :

$$
\begin{array}{ccc} A_0=I_0=\left[ \begin{array}{ccc} 0 & 0 & 1 \\ 0 & -1 & 0 \\ 1 & 0 & 0 \end{array}\right] & 
  A_2 = \left[ \begin{array}{ccc} -1 & 0 & i+j \\ 0 & 1 & 0 \\ i+j & 0 & 1 \end{array}\right] \\\\

A_{3,1} = \left[ \begin{array}{ccc} 1 & 0 & 0 \\ 0 & 1 & 0 \\ -i+j-k & 0 & 1 \end{array}\right]&
  A_{3,2}=\left[ \begin{array}{ccc} 1 & 0 & 0 \\ 0 & 1 & 0 \\ i-j+k & 0 & 1 \end{array}\right] \\\\

A_{3,3} = \left[ \begin{array}{ccc} 1 & 0 & 0 \\ \frac{1-i-j+k}{2} & 1 & 0 \\ \frac{-1-3i+j+k}{2} & \frac{-1-i-j+k}{2} & 1 \end{array}\right] &
  A_{3,4} = \left[ \begin{array}{ccc} 1 & 0 & 0 \\ -k & 1 & 0 \\ \frac{-1+i-3j-k}{2} & -k & 1 \end{array}\right]\\\\

A_{3,5}=\left[ \begin{array}{ccc} 1 & 0 & 0 \\ \frac{1-i+j-k}{2} & 1 & 0 \\ \frac{-1-i-j-3k}{2} & \frac{-1-i+j-k}{2} & 1 \end{array}\right] & 
  A_{4,1} = \left[ \begin{array}{ccc} i & 0 & 1 \\ 0 & -1 & 0 \\ 2 & 0 & -i \end{array}\right] \\\\

A_{4,2} = \left[ \begin{array}{ccc} i+j+k & 0 & 2 \\ 0 & -1 & 0 \\ 2 & 0 & -i-j-k \end{array}\right]&
  A_{4,3}=\left[ \begin{array}{ccc} 1 & 0 & 0 \\ 1-i & 1 & 0 \\ -1-i-j-k & -1-i & 1 \end{array}\right] \\\\ 
  
A_{4,4} = \left[ \begin{array}{ccc} 1 & 0 & 0 \\ -i-k & 1 & 0 \\ -1+i-j-k & -i-k & 1 \end{array}\right] &
  A_{4,5} = \left[ \begin{array}{ccc} 1 & 0 & 0 \\ -i+k & 1 & 0 \\ -1-i+j-k & -i+k & 1 \end{array}\right]\\\\

A_{4,6}=\left[ \begin{array}{ccc} 1 & 0 & 0 \\ -1-i & 1 & 0 \\ -1+i+j-k & 1-i & 1 \end{array}\right] & 
  A_{4,7} = \left[ \begin{array}{ccc} 1 & 0 & 0 \\ 1-j & 1 & 0 \\ -1-i-j+k & -1-j & 1 \end{array}\right] \\\\

A_{4,8} = \left[ \begin{array}{ccc} 1 & 0 & 0 \\ -j-k & 1 & 0 \\ -1+i-j+k & -j-k & 1 \end{array}\right]&
  A_{4,9}=\left[ \begin{array}{ccc} 1 & 0 & 0 \\ -j+k & 1 & 0 \\ -1-i+j+k & -j+k & 1 \end{array}\right] \\\\ 
  
A_{4,10} = \left[ \begin{array}{ccc} 1 & 0 & 0 \\ -1-j & 1 & 0 \\ -1+i+j+k & 1-j & 1 \end{array}\right] &\\

\end{array}$$

\medskip

{\bf Relations:} 
We find a complete set of relations by applying generators to points of depth at most 4 as described in part (5) of section~\ref{method}. The detailed steps of the computation can be found in the companion files in \cite{MCode}. The direct output is a presentation with 33 relations and 968,480 relations, which is unfortunately too large for Magma to handle directly (the text file for this presentation, {\tt quaternions.m}, is slightly over 200 MB and can also be found at \cite{MCode}). 

To illustrate the steps involved we compute by hand some of the cycles and relations corresponding to triples $(a,b,c)$ for which the point $A_a(p_b)$ has depth at most 4, hence is of the form $\gamma_\infty (p_c)$ for some $\gamma_\infty \in \Gamma_\infty$ and $p_c \in D_\infty$ with depth at most 4. In the notation of part (5) of section~\ref{method}, this corresponds to $\gamma_\infty^1={\rm Id}$ and  $\gamma_\infty^3=\gamma_\infty$. The results are listed in Tables~5--9, with the same notation as in the example of section~\ref{toyexample}.

We would also like to make explicit certain relations in the group, which essentially follow from these tables.

\begin{prop}\label{hurwitzrels} The following relations hold among the generators of $\Gamma_\infty(\mathcal{H})$ and $I=I_0$:
\begin{equation}
{\renewcommand{\arraystretch}{1.2}
\begin{array}{c}
I^2, [\ I,R_i], \ [I,R_\sigma], \ [I,C_i], \ [I,C_\sigma],  \
(T_{v_i}I)^3=R_iC_i, \
 (T_{v_i}T_{v_j}I)^4=R_i^2, \ (T_{v_i}T_{v_j}T_{v_k}I)^6, \\

 [IT_{v_i}^{-2}IR_iT_{v_i}^{-1}C_i,R_\sigma], \
     [C_\sigma^{-1}(IT_{v_i}^{-1}T_{v_j}^{-1}T_{v_k}^{-1})^3R_i^2C_\sigma,R_i], \
[C_\sigma^{-1}(IT_{v_i}^{-1}T_{v_j}^{-1}T_{v_k}^{-1})^3R_i^2C_\sigma,R_\sigma], \\

[C_\sigma^{-1}(IT_{v_i}^{-1}T_{v_j}^{-1}T_{v_k}^{-1})^3R_i^2C_\sigma,C_\sigma], \
 (IT_{v_i}^{-2}IR_iT_{v_i}^{-1}C_i)^2, \ 
\ IT_{v_i}^{-1}(IT_{v_j})^2=(T_{v_j}IT_{v_i}^{-1})^2R_iC_i,
 \\ 
 
 [(T_{v_i}T_{v_j}I)^2,C_\sigma^{-1}(IT_{v_i}^{-1}T_{v_j}^{-1}T_{v_k}^{-1})^3R_i^2C_\sigma], \ 
IT_1T_\sigma^{-1}IT_\sigma I=T_1 T_{v_i}^{-1}T_{v_j}^{-1}T_{v_k}^{-1}C_\sigma^{-1}
\end{array}
}
\end{equation}
\end{prop}

It is straightforward to check that each of these relations hold. However we find it instructive to illustrate how some of  them can be deduced from the cycle relations in the tables. We now give more details on how to obtain the relations $(T_{v_i}T_{v_j}T_{v_k}I)^6$ and $[(T_{v_i}T_{v_j}I)^2,C_\sigma^{-1}(IT_{v_i}^{-1}T_{v_j}^{-1}T_{v_k}^{-1})^3R_i^2C_\sigma]$.  The rest of the details are either straighforward or similar to these two.

Substituting for $A_2$ in (2) using (1) (both from Table~6), we obtain
$$ A_{4,2} = I T_{v_i}^{-1} T_{v_j}^{-1} I T_{v_i}^{-1} T_{v_j}^{-1} T_{v_k} T_{v_i} T_{v_j} I T_{v_i} T_{v_j} I C_{\sigma}^{-1} R_i C_{\sigma} T_{v_i} T_{v_j} T_{v_k}^{-1} C_{\sigma}^{-1} C_i C_{\sigma} $$
Use the relations $C_{\sigma}T_{v_w}C_{\sigma}^{-1} = T_{v_{\hat{w}}}$ and $[R_i,T_{v_w}] = 1$ from $\Gamma_{\infty}(\mathcal{H})$ and the relations (5) and (6) from Table 9 to get
$$A_{4,2} = C_{\sigma}^{-1} I T_{v_j}^{-1} T_{v_k}^{-1} I T_{v_i} I T_{v_j} T_{v_k} I T_{v_j} T_{v_k} T_{v_i}^{-1} R_i C_i C_{\sigma} $$
Use the relation $(I T_{v_i})^3 = R_iC_i$ 
to make the subsitution $IT_{v_i} I = T_{v_i}^{-1} I T_{v_i}^{-1} R_i C_i$. 
$$A_{4,2} = C_{\sigma}^{-1} I T_{v_j}^{-1} T_{v_k}^{-1} T_{v_i}^{-1} I T_{v_i}^{-1} R_i C_i T_{v_j} T_{v_k} I T_{v_j} T_{v_k} T_{v_i}^{-1} R_i C_i C_{\sigma} $$
 Use the relations $[C_i,T_{v_i}] = 1$ and $C_i T_{v_w} C_i^{-1} = T_{v_w}^{-1}$ if $w\neq i$ from $\Gamma_{\infty}(\mathcal{H})$ to get
 $$
\begin{array}{r @{\,\,=\,\,}l}
  A_{4,2} & C_{\sigma}^{-1} I T_{v_j}^{-1} T_{v_k}^{-1} T_{v_i}^{-1} I T_{v_i}^{-1} T_{v_j}^{-1} T_{v_k}^{-1} I T_{v_j}^{-1} T_{v_k}^{-1} T_{v_i}^{-1} R_i^2 C_{\sigma} \\
  &  C_{\sigma}^{-1} (I T_{v_i}^{-1} T_{v_j}^{-1} T_{v_k}^{-1})^3 R_i^2 C_{\sigma}
  \end{array}
$$

The relation $(T_{v_i}T_{v_j}T_{v_k} I)^6$ comes from substituting this expression for $A_{4,2}$ in the relation (4) from Table 8.  The relation $[(T_{v_i}T_{v_j}I)^2,C_\sigma^{-1}(IT_{v_i}^{-1}T_{v_j}^{-1}T_{v_k}^{-1})^3R_i^2C_\sigma]$ comes from first observing that (3) in Table 8 and (2) in Table 6 are conjugate to each other by $A_2$.  Since both are equal to $A_{4,2}$, we obtain $[A_2,A_{4,2}] = 1$.  Then to get the relation that appears in the presentation, we substitute for $A_2$ using (1) from Table 6 and for $A_{4,2}$ using the expression obtained above.

\begin{cor} The abelianization of $\Gamma(\mathcal{H})$ is $\Z/3\Z$, generated by (the image of) $R_\sigma$, or trivial.
\end{cor}

\Pf By Theorem~\ref{hurwitzgens}, $\Gamma=\Gamma(\mathcal{H})$ is generated by $\Gamma_{\infty}(\mathcal{H})$ and $I_0$, hence by $\mathcal{S}_\infty$ and $I_0$ (where $\mathcal{S}_\infty$ denotes our generating set for $\Gamma_\infty(\mathcal{H})$, see Lemma~\ref{inftypres}). Therefore $\Gamma$ is a quotient of the abstract group $\Delta$ with generators $\mathcal{S}_\infty$ and $I_0$ and relations $\mathcal{R}_\infty$ together with the relations listed in Proposition~\ref{hurwitzrels} (with $\mathcal{S}_\infty$ the set of relations for $\Gamma_\infty(\mathcal{H})$, see Lemma~\ref{inftypres}). By Corollary~\ref{inftyabel}, the abelianization of $\Delta$ is a quotient of $(\Z/3\Z)^2$, generated by (the images of) $R_\sigma$ and $C_\sigma$. By inspecting the relations listed in Proposition~\ref{hurwitzrels}, we see that the abelianization of $\Delta$ is $\Z/3\Z$, generated by (the image of) $R_\sigma$.  \EPf

We then obtain Theorem~\ref{hurwitzabel} by combining these results with a computational check using the full presentation. More specifically, while the full presentation is too large for Magma to handle at once, we found that it could handle files approximately 10 times smaller in a few hours. To check that the abelianization remains unchanged after adding the roughly $10^6$ relations of the mondo presentation, we subdivide the file with all relations into 10 pieces (each with approximately $10^5$ relations), see files {\tt mondo1.m},...,{\tt mondo10.m} at \cite{MCode} (each of these also contains the relations from $\Gamma_\infty(\mathcal{H})$). Recall from Corollary~\ref{inftyabel} that the  abelianization of $\Gamma_\infty(\mathcal{H})$ is $(\Z/3\Z)^2$. We add to each file enough relations to bring the abelianization down to $\Z/3\Z$ - in practice we simply add the relations $C_\sigma=1=I_0$ (which hold in the abelianization by the above, but of course not in the group). This produces the files {\tt mondo1modified.m},...,{\tt mondo10modified.m}, for which Magma is able to compute that all 10 abelianizations are $\Z/3\Z$,  proving Theorem~\ref{hurwitzabel}.

\begin{table}[h]\label{hurwitztableI}
  \begin{center}
  {\renewcommand{\arraystretch}{1.2}%
\begin{tabular}{|c|c|c|c|}
\hline 
$\mathbf{A.p}$ & $\mathbf{p'}$ & $\mathbf{A_{p'}^{-1}W^{-1}A A_p}$ & $\mathbf{W'}$  \\
\hline

$I p_0$ & $\infty$ & $I^2$ & ${\rm Id}$ \\
\hline

$I (\sigma,i+j+k)$ & $p_0$ & $IT_1T_\sigma ^{-1}IT_\sigma I$ & $T_1T_{v_i}^{-1}T_{v_j}^{-1}T_{v_k}^{-1}C_\sigma ^{-1}$\\
\hline

$I p_2$ & $p_0$ & $IT_{v_i}T_{v_j}IA_2$ & $T_{v_i}^{-1}T_{v_j}^{-1}R_i^2$\\
\hline

$I p_{3,1}$ & $p_0$ & $IT_{v_i}T_{v_j}T_{V_k}IA_{3,1}$ & ${\rm Id}$\\
\hline

$I p_{3,2}$ & $p_0$ & $IT_{v_i}^{-1}T_{v_j}T_{V_k}IA_{3,2}$ & ${\rm Id}$\\
\hline

$I p_{3,3}$ & $p_0$ & $IT_\sigma T_i^{-1}T_{v_j}T_{v_k}IA_{3,3}$ & ${\rm Id}$\\
\hline

$I p_{4,1}$ & $p_0$ & $IT_{v_i}^2IA_{4,1}$ & $T_{v_i}^{-1}R_iC_i$\\
\hline

$I p_{4,2}$ & $p_{3,1}$ & $A_{3,1}^{-1}T_{v_i}T_{v_j}T_{v_k}IA_{4,2}$ & $R_i^1T_{v_i}^{-1}T_{v_j}^{-1}T_{v_k}^{-1}$\\
\hline

$I p_{4,3}$ & $p_0$ & $IT_{v_i}^2T_{v_j}T_{v_k}^2T_k^{-1}T_i^{-1}IA_{4,3}$ & ${\rm Id}$\\
\hline

$I p_{4,4}$ & $p_0$ & $IT_k^{-1}T_1^{-1}T_{v_j}^2T_{v_k}IA_{4,4}$ & ${\rm Id}$\\
\hline





\end{tabular}
}
\end{center}
\caption{Cycles coming from $I$}
\end{table}

\begin{table}[h]\label{hurwitztableA2}
  \begin{center}
  {\renewcommand{\arraystretch}{1.2}%
\begin{tabular}{|c|c|c|c|}
\hline 
$\mathbf{A.p}$ & $\mathbf{p'}$ & $\mathbf{A_{p'}^{-1}W^{-1}A A_p}$ & $\mathbf{W'}$  \\
\hline

$A_2 p_0$ & $p_0$ & $I_0T_{v_j}^{-1}T_{v_i}^{-1}A_2I_0$ & \hfill $A_2=(T_{v_i}T_{v_j}I_0)^2$ \hfill (1)\\
\hline

$A_2 p_2$ & $\infty$ & $A_2^2$ & $R_i^2$\\
\hline

$A_2 p_{3,1}$ & $p_0$ & $I_0T_{v_i}^{-1}T_{v_j}^{-1}T_{v_k}^{-1}A_2A_{3,1}$ & $R_\sigma^{-1}R_i^{-1}R_\sigma C_\sigma^{-1}C_i C_\sigma T_{v_k}$\\
\hline


$A_2 p_{3,3}$ & $p_{3,3}$ & $A_{3,3}^{-1}A_2^{-1}R_i^{-1}T_{v_i}T_{v_j}C_\sigma^{-1}C_iC_\sigma A_{3,3}$ & $C_\sigma^{-1}C_iC_\sigma R_i$\\
\hline

$A_2 p_{4,1}$ & $p_2$ & $A_2^{-1}T_{v_i}^{-1}A_2A_{4,1}$ & $C_\sigma R_i^{-1}C_\sigma ^{-1}  T_{v_i}T_{v_j}C_\sigma C_i C_\sigma ^{-1}$\\
\hline

$A_2 p_{4,2}$ & $p_2$ & $A_2^{-1}T_{v_k}^{-1}A_2A_{4,2}$ & \hfill $C_\sigma^{-1} R_iC_\sigma  T_{v_i}T_{v_j}T_{v_k}^{-1}C_\sigma ^{-1} C_i C_\sigma $ \hfill (2)\\
\hline

$A_2 p_{4,3}$ & $p_2$ & $A_2^{-1}T_\sigma^{-1}A_2A_{4,3}$ & $T_1^{-1} T_\sigma T_k^{-1}  T_{v_i} C_i$\\
\hline





\hline

\end{tabular}
}
\end{center}
\caption{Cycles coming from $A_2$}
\end{table}

\begin{table}[h]\label{hurwitztableA3}
  \begin{center}
  {\renewcommand{\arraystretch}{1.2}%
\begin{tabular}{|c|c|c|c|}
\hline 
$\mathbf{A.p}$ & $\mathbf{p'}$ & $\mathbf{A_{p'}^{-1}W^{-1}A A_p}$ & $\mathbf{W'}$  \\
\hline

$A_{3,1} p_0$ & $p_0$ & $IA_{3,1}I$ & $T_{v_i}^{-1}T_{v_j}^{-1}T_{v_k}^{-1}$\\
\hline

$A_{3,2} p_0$ & $p_0$ & $IA_{3,2}I$ & $T_{v_i}T_{v_j}^{-1}T_{v_k}^{-1}$\\
\hline


$A_{3,2} p_{4,1}$ & $p_{3,1}$ & $A_{3,1}^{-1}A_{3,2}A_{4,1}$ & $R_iC_iT_{v_i}^{-1}$\\
\hline

$A_{3,3} p_0$ & $p_0$ & $IA_{3,3}I$ & $T_{v_j}^{-1}T_{v_k}^{-1}T_iT_\sigma^{-1}$\\
\hline

\end{tabular}
}
\end{center}
\caption{Cycles coming from $A_{3,1}, A_{3,2}, A_{3,3}$}
\end{table}

\begin{table}[h]\label{hurwitztableA4}
  \begin{center}
  {\renewcommand{\arraystretch}{1.2}%
\begin{tabular}{|c|c|c|c|}
\hline 
$\mathbf{A.p}$ & $\mathbf{p'}$ & $\mathbf{A_{p'}^{-1}W^{-1}A A_p}$ & $\mathbf{W'}$  \\
\hline

$A_{4,1} p_0$ & $p_0$ & $IT_{v_i}^{-1}A_{4,1}I$ & $T_{v_i}^2R_i^{-1}C_i$\\
\hline

$A_{4,1} p_2$ & $p_2$ & $A_2^{-1}T_{v_j}A_{4,1}A_2$& $T_{v_i}^{-1}T_{v_j}R_iC_i$\\
\hline

$A_{4,1} p_{3,1}$ & $p_{3,2}$ & $A_{3,2}^{-1}C_iT_{v_i}^{-1}A_{4,1}A_{3,1}$ & $R_i^{-1}$\\
\hline


$A_{4,1} p_{4,1}$ & $\infty$ & $A_{4,1}^2$ & ${\rm Id}$\\
\hline

$A_{4,1} p_{4,3}$ & $p_{4,4}$ & $A_{4,4}^{-1}C_iT_{v_i}^{-1}R_i^{-1}R_\sigma^2A_{4,1}A_{4,3}$ & $R_\sigma^{-1}$\\
\hline


$A_{4,2} p_0$ & $p_{3,2}$ & $A_{3,2}^{-1}C_iT_{v_i}^{-1}T_{v_j}^{-1}T_{v_k}^{-1}A_{4,2}I$ & $T_{v_i}T_{v_j}T_{v_k}R_i^2C_i $\\
\hline

$A_{4,2} p_2$ & $p_2$ & $A_2^{-1}T_{v_k}^{-1}A_{4,2}A_2$ & \hfill $T_{v_i}T_{v_j}T_{v_k}^{-1}C_\sigma^{-1}R_iC_iC_\sigma$ \hfill (3)\\
\hline

$A_{4,2} p_{3,1}$ & $p_0$ & $IT_{v_i}^{-1}T_{v_j}^{-1}T_{v_k}^{-1}A_{4,2}A_{3,1}$ & $T_{v_i}T_{v_j}T_{v_k}R_i^2$\\
\hline

$A_{4,2} p_{3,3}$ & $p_{3,3}$ & $A_{3,3}^{-1}C_\sigma^{-1}C_iC_\sigma T_{v_i}^{-1}T_{v_j}^{-1}R_i^2A_{4,2}A_{3,3}$ & $C_iT_j^{-1}T_k^{-1}T_{v_j}T_{v_k}^2$\\
\hline

$A_{4,2} p_{4,2}$ & $\infty$ & $A_{4,2}^2$ & \hfill ${\rm Id}$ \hfill (4)\\
\hline

$A_{4,2} p_{4,3}$ & $p_{4,4}$ & $A_{4,4}^{-1}C_iT_{v_i}^{-1}T_{v_j}^{-1}T_{v_k}^{-1}R_\sigma^{-1}R_iA_{4,2}A_{4,3}$ & $C_iT_{v_i}^{-1}T_{v_k}^{-1}R_\sigma R_iR_\sigma^{-1}R_i^{-1}R_\sigma^{-1}T_iT_k$\\
\hline



$A_{4,3} p_0$ & $p_0$ & $IA_{4,3}I$ & $T_iT_kT_{v_i}^{-2}T_{v_j}^{-1}T_{v_k}^{-2}$\\
\hline

$A_{4,4} p_0$ & $p_0$ & $IA_{4,4}I$ & $T_1T_kT_{v_j}^{-2}T_{v_k}^{-1}$\\
\hline

$A_{4,4} p_{4,1}$ & $p_{4,3}$ & $A_{4,3}^{-1}R_\sigma A_{4,4}A_{4,1}$ & $C_iR_\sigma R_i^{-1}R_\sigma R_i^{-1}R_\sigma^{-1}T_{v_i}^{-1}$\\
\hline

\end{tabular}
}
\end{center}
\caption{Cycles coming from $A_{4,1}, A_{4,2}, A_{4,3}, A_{4,4}$}
\end{table}

\begin{table}[h]\label{hurwitztableRC}
  \begin{center}
  {\renewcommand{\arraystretch}{1.2}%
\begin{tabular}{|c|c|c|}
\hline 
Image of a point & Cycle of points & Relation \\
\hline

$R_i p_0=p_0$ & $\infty \xrightarrow{I_0} p_0 \xrightarrow{R_i} p_0 \xrightarrow{I_0}  \infty$& \hfill $[R_i,I_0]={\rm Id}$ \hfill (5)\\
\hline

$R_\sigma p_0=p_0$ & $\infty \xrightarrow{I_0} p_0 \xrightarrow{R_\sigma} p_0 \xrightarrow{I_0}  \infty$& $[R_\sigma,I_0]={\rm Id}$\\
\hline

$C_i p_0=p_0$ & $\infty \xrightarrow{I_0} p_0 \xrightarrow{C_i} p_0 \xrightarrow{I_0}  \infty$& \hfill $[C_i,I_0]={\rm Id}$ \hfill (6)\\
\hline

$C_\sigma p_0=p_0$ & $\infty \xrightarrow{I_0} p_0 \xrightarrow{C_\sigma} p_0 \xrightarrow{I_0}  \infty$& $[C_\sigma,I_0]={\rm Id}$\\
\hline

$R_i p_2=p_2$ & $\infty \xrightarrow{A_2} p_2 \xrightarrow{R_i} p_2 \xrightarrow{A_2^{-1}}  \infty$& $[R_i,A_2]={\rm Id}$\\
\hline

$R_\sigma p_2=p_2$ & $\infty \xrightarrow{A_2} p_2 \xrightarrow{R_\sigma} p_2 \xrightarrow{A_2}^{-1}  \infty$& $[R_\sigma,A_2]={\rm Id}$\\
\hline

$R_i p_{4,1}=p_{4,1}$ & $\infty \xrightarrow{A_{4,1}} p_{4,1} \xrightarrow{R_i} p_{4,1} \xrightarrow{A_{4,1}^{-1}}  \infty$& $[R_i,A_{4,1}]={\rm Id}$\\
\hline

$R_\sigma p_{4,1}=p_{4,1}$ & $\infty \xrightarrow{A_{4,1}} p_{4,1} \xrightarrow{R_\sigma} p_{4,1} \xrightarrow{A_{4,1}^{-1}}  \infty$& $[R_\sigma,A_{4,1}]={\rm Id}$\\
\hline

$R_i p_{4,2}=p_{4,2}$ & $\infty \xrightarrow{A_{4,2}} p_{4,2} \xrightarrow{R_i} p_{4,2} \xrightarrow{A_{4,2}^{-1}}  \infty$& $[R_i,A_{4,2}]={\rm Id}$\\
\hline

$R_\sigma p_{4,2}=p_{4,2}$ & $\infty \xrightarrow{A_{4,2}} p_{4,2} \xrightarrow{R_\sigma} p_{4,2} \xrightarrow{A_{4,2}^{-1}}  \infty$& $[R_\sigma,A_{4,2}]={\rm Id}$\\
\hline



$C_i p_{4,1}=p_{4,1}$ & $\infty \xrightarrow{A_{4,1}} p_{4,1} \xrightarrow{C_i} p_{4,1} \xrightarrow{A_{4,1}^{-1}}  \infty$& $[C_i,A_{4,1}]={\rm Id}$\\
\hline

$C_\sigma p_{4,2}=p_{4,2}$ & $\infty \xrightarrow{A_{4,2}} p_{4,2} \xrightarrow{C_\sigma} p_{4,2} \xrightarrow{A_{4,2}^{-1}}  \infty$& $[C_\sigma,A_{4,2}]={\rm Id}$\\
\hline

\end{tabular}
}
\end{center}
\caption{Action of generators on vertices for $\mathcal{H}$, degenerate cycles coming from $R$'s and $C$'s}
\end{table}

\section{Appendix: presentation for the Eisenstein-Picard modular group ${\rm PU}(2,1,\mathcal{O}_3)$}

The presentation for the Eisenstein-Picard modular group ${\rm PU}(2,1,\mathcal{O}_3)$ which we obtain is $\la T_1, R, I \, | \, \mathcal{R}_3 \ra$, where $\mathcal{R}_3$ is the following set of 19 relations:
$$\begin{array}{l}
I^2 = {\rm Id} \\
    R^{-1}  I  R  I = {\rm Id} \\
    R^6 = {\rm Id} \\
    R  T_1  R^{-1}  T_1^{-1}  R^{-1}  T_1  R = {\rm Id} \\
    I  T_1  I  R  T_1  I  T_1^{-1}  I  R^{-1}  T_1^{-1} = {\rm Id} \\
    T_1^{-1}  I  T_1  I  R  T_1  I  T_1^{-1}  R^{-1}  I = {\rm Id} \\
    R^{-2}  T_1^{-1}  I  R^2  I  R^{-1}  T_1  R  T_1^{-1} = {\rm Id} \\
    I  R^2  T_1  I  R^2  T_1  R^2  I  T_1 = {\rm Id} \\
    I  R  T_1^{-1}  R^{-1}  T_1  I  T_1^{-1}  I  T_1  R^{-1} 
    T_1^{-1}  R = {\rm Id} \\
    T_1^{-1}  R  T_1^2  R^{-1}  T_1^{-1}  R  T_1^{-1}  R^{-1}  T_1 
    R^{-1}  T_1  R = {\rm Id} \\
    I  T_1^{-1}  R  T_1  I  R^{-3}  T_1  R  I  R  T_1 =
    {\rm Id} \\
    I  R^2  T_1^{-1}  R^{-1}  I  R^{-1}  T_1^{-1}  R^{-2}  I 
    T_1^{-1}  R^2 = {\rm Id} \\
    I  R^{-1}  T_1  I  R^{-1}  I  T_1^{-1}  R^{-2}  T_1  I 
    T_1^{-1}  I  R^{-2}  T_1^{-1} = {\rm Id} \\
    I  T_1^{-1}  R^{-2}  T_1  I  T_1^{-1}  R^{-2}  T_1  I  T_1^{-1} 
    I  T_1^{-1}  R^{-1}  I  T_1  R^{-1} = {\rm Id} \\
    T_1  I  T_1^2  I  R  T_1^{-1}  I  R  T_1  I  R^{-2} 
    T_1  I  T_1^{-2}  I  T_1^{-1}  I = {\rm Id} \\
    T_1  R^{-1}  I  T_1^{-1}  R^2  I  T_1^{-1}  I  T_1  I 
    R^{-1}  T_1  I  T_1^{-1}  I  R  T_1^{-1}  R^{-1}  I = {\rm Id} \\
    T_1  I  T_1^{-1}  R^2  T_1^{-1}  R^{-1}  T_1  R^{-1}  T_1  I 
    T_1^{-1}  R^{-2}  T_1  I  T_1^{-1}  R  I  T_1^{-1}  R^2  T_1^{-1} 
    R^{-1} = {\rm Id} \\
    T_1^{-2}  R  T_1  I  T_1  I  T_1^{-1}  R^2  T_1  R  I 
    T_1^{-1}  R  T_1^2  R^{-1}  I  T_1^{-1}  R^{-2}  T_1  R^{-3}  T_1 
    I  T_1  R = {\rm Id} \\
    (T_1^{-1}  R  T_1  R^{-2}  T_1  I  T_1  I  T_1^{-1}  R^{-2}  T_1
     I T_1^{-1})^2 = {\rm Id} 
\end{array}$$

\raggedright

\begin{flushleft}

 \textsc{Alice Mark\\
   Department of Mathematics, Vanderbilt University}\\
       \verb|alice.h.mark@vanderbilt.edu|
\\

  \textsc{Julien Paupert\\
   School of Mathematical and Statistical Sciences, Arizona State University}\\
       \verb|paupert@asu.edu|
\end{flushleft}

\end{document}